\documentclass[hidelinks,journal]{IEEEtran}

\IEEEoverridecommandlockouts                              

\usepackage{xr}

\newcommand{\omitted}[1]{}

\title{\fontsize{22}{15} \selectfont Robust and Adaptive Sequential Submodular Optimization
}

\author{Vasileios Tzoumas,$^{\star}$~\IEEEmembership{Member,~IEEE,} Ali Jadbabaie,{$^{\dagger}$}~\IEEEmembership{Fellow,~IEEE,} George J.~Pappas,{$^{\ddagger}$}~\IEEEmembership{Fellow,~IEEE.}
	\thanks{$^\star$Department of Aerospace Engineering, University of Michigan, Ann Arbor, MI 48109, USA.  At the time the paper was accepted for publication: Laboratory for Information and Decision Systems, Massachusetts Institute of Technology, Cambridge, MA 02139, USA. {\tt\footnotesize vtzoumas@umich.edu}}
	\thanks{$^{\dagger}$Institute for Data, Systems and Society, Massachusetts Institute of Technology, Cambridge, MA 02139 USA. {\tt\footnotesize jadbabai@mit.edu}}
	\thanks{$^{\ddagger}$Department of Electrical and Systems Engineering, University of Pennsylvania, Philadelphia, PA 19104 USA. {\tt\footnotesize pappasg@seas.upenn.edu}}
	\thanks{The work was supported  
		by the AFOSR Complex Networks program, and by the ARL CRA DCIST W911NF-17-2-0181 program.}
}

\usepackage{cite}

\usepackage{comment}
\usepackage{siunitx}
\usepackage{relsize}
\usepackage{ifthen}
\usepackage[colorinlistoftodos]{todonotes}

\usepackage[caption=false]{subfig}

\usepackage[vlined,ruled,linesnumbered]{algorithm2e}
\usepackage{graphics} 
\usepackage{rotating}
\usepackage{color}
\usepackage{enumerate}
\usepackage[T1]{fontenc}
\usepackage{psfrag}
\usepackage{epsfig} 
\usepackage{booktabs}
\usepackage{graphicx,url}
\usepackage{multirow}
\usepackage{array}
\usepackage{latexsym}
\usepackage{amsfonts}
\usepackage{amsmath}
\usepackage{amssymb}
\usepackage{xstring}
\usepackage{algorithmic}
\usepackage{multirow}
\usepackage{xcolor}
\usepackage{prettyref}
\usepackage{flexisym}
\usepackage{bigdelim}
\usepackage{breqn} 
\usepackage{listings}
\let\labelindent\relax
\usepackage{xspace}
\usepackage{bm}
\graphicspath{{./figures/}}
\usepackage{tikz}
\usetikzlibrary{matrix,calc}

\usepackage{mdwlist}

\let\itemize\undefined
\makecompactlist{itemize}{stditemize}
\usepackage{enumitem}

\usepackage{amsthm}

\newtheorem{theorem}{Theorem}

\newtheorem{corollary}[theorem]{Corollary}

\newtheorem{lemma}[theorem]{Lemma}

\newtheorem{definition}[theorem]{Definition}
\newtheorem{proposition}[theorem]{Proposition}

\newtheorem{remark}[theorem]{Remark}

\newcommand{\bdmath}{\begin{dmath}}
\newcommand{\edmath}{\end{dmath}}
\newcommand{\beq}{\begin{equation}}
\newcommand{\eeq}{\end{equation}}
\newcommand{\bdm}{\begin{displaymath}}
\newcommand{\edm}{\end{displaymath}}
\newcommand{\bea}{\begin{eqnarray}}
\newcommand{\eea}{\end{eqnarray}}
\newcommand{\beal}{\beq \begin{array}{lll}}
\newcommand{\eeal}{\end{array} \eeq}
\newcommand{\beas}{\begin{eqnarray*}}
\newcommand{\eeas}{\end{eqnarray*}}
\newcommand{\ba}{\begin{array}}
\newcommand{\ea}{\end{array}}
\newcommand{\bit}{\begin{itemize}}
\newcommand{\eit}{\end{itemize}}
\newcommand{\ben}{\begin{enumerate}}
\newcommand{\een}{\end{enumerate}}

\newcommand{\calA}{{\cal A}}
\newcommand{\calB}{{\cal B}}

\newcommand{\calI}{{\cal I}}

\newcommand{\calK}{{\cal K}}

\newcommand{\calM}{{\cal M}}

\newcommand{\calO}{{\cal O}}
\newcommand{\calP}{{\cal P}}

\newcommand{\calR}{{\cal R}}
\newcommand{\calS}{{\cal S}}

\newcommand{\calV}{{\cal V}}

\newcommand{\calX}{{\cal X}}

\definecolor{myblue}{RGB}{65 105 225}

\newcommand{\hide}[1]{}

\newcommand{\hiddenText}{{\color{gray} hidden text.}}
\newcommand{\hideWithText}[1]{\hiddenText}

\newcommand{\eye}{{\mathbf I}}

\newcommand{\Real}[1]{ { {\mathbb R}^{#1} } }

\newcommand{\scenario}[1]{{\smaller \sf#1}\xspace}

\newcommand{\rsm}{\texttt{RSM}\xspace}

\newcommand{\alg}{\scenario{RAM}}
\newcommand{\post}{post-failure\xspace}

\newcommand{\alglongCap}{Robust and Adaptive Maximization\xspace}

\newcommand{\elem}{v}

\newcommand{\myParagraph}[1]{{\bf #1.}\xspace}

\newcommand{\random}{\scenario{random}}
\newcommand{\greedy}{\scenario{greedy}}

\newcommand{\optimal}{\scenario{optimal}}

\PassOptionsToPackage{end}{algorithmic}

\begin{document}
	
\maketitle

\thispagestyle{empty}
\pagestyle{empty}

\begin{abstract}
Emerging applications of control, estimation, and machine learning, from target tracking to decentralized model fitting, pose resource constraints that limit which of the available sensors, actuators, or data can be simultaneously used across time.  Therefore, many researchers have proposed solutions within discrete optimization frameworks where the optimization is performed over finite sets.
By exploiting notions of discrete convexity, such as submodularity, the researchers have been able to provide scalable algorithms with provable suboptimality bounds.  In this paper, we consider such problems but in adversarial environments, where in every step a number of the chosen elements in the optimization is removed due to failures/attacks.   
Specifically, we consider for the first time a sequential version of the problem that allows us to observe the failures and adapt, while the attacker also adapts to our response. We call the novel problem \textit{Robust Sequential submodular Maximization} (\rsm). Generally, the problem is computationally hard and no scalable algorithm is known for its solution.  However, in this paper we propose \textit{Robust and Adaptive Maximization} (\alg), the first scalable algorithm.  
\alg  runs in an online fashion, adapting in every step to the history of failures. Also, it guarantees a near-optimal performance, even against any number of failures among the used elements. Particularly, \alg has both provable per-instance a priori bounds and tight and/or optimal a posteriori bounds. Finally, we demonstrate \alg's near-optimality in simulations across various application scenarios, along with its robustness \mbox{against several failure types, from worst-case to random.}
\end{abstract}

\section{Introduction}\label{sec:Intro}

Control, estimation, and machine learning applications of the Internet of Things (IoT) and autonomous robots~\cite{abdelzaher2018will}
require the sequential optimization of systems in scenarios such as:

\paragraph*{{\small\textbullet}~Sensor scheduling} An unmanned aerial vehicle (UAV) is assisted for its navigation by on-board and on-ground sensors.  Ideally, the UAV would use all available sensors for navigation. However, limited on-board capacity for measurement-processing necessitates a sequential sensor scheduling problem~\cite{carlone2018attention}: at each time step, which few sensors should be used for the UAV to effectively navigate itself?

\paragraph*{{\small\textbullet}~Target tracking} A wireless sensor network (WSN) is designated to monitor a mobile target.  Limited battery power necessitates a sequential sensor activation problem~\cite{he2006achieving}: at each time step, which few sensors should be activated for the WSN to effectively track the target?

\paragraph*{{\small\textbullet}~Decentralized model fitting} A team of mobile robots collects data to learn the model of an unknown environmental process.  The data are transmitted to a fusion center, performing the statistical analysis.  Ideally, all robots would transmit their data to the center at the same time.  But instead, communication bandwidth constraints necessitate a sequential transmission problem~\cite{yang2019scheduling}:  at each time step, which few robots should transmit their  data for the \mbox{center to effectively learn the model?} 

Similar applications of sensor and data scheduling, but also of actuator scheduling as well as infrastructure design are studied in~\cite{gupta2006stochastic,jawaid2015submodularity,clark2016submodularity,tzoumas2016scheduling,zhang2017kalman,tzoumas2018codesign,zhao2016scheduling,nozari2017time,ikeda2018sparsity,leskovec2007cost,poularakis2017one,summers2014submodularity}.
Particularly, all applications above require the {sequential} selection of {a few} elements, among a finite set of available ones, to optimize performance across multiple steps subject to resource constraints. For example, the target tracking application above requires the sequential activation of a few sensors across the WSN, to optimize an estimation error 
subject to power constraints. Importantly, the activated sensors may vary in time, since each sensor may measure different parts of the target's state (e.g., some sensors may measure only position, others only speed). 
Formally,  all above applications motivate the sequential optimization problem\footnote{ 
Calligraphic fonts denote finite discrete sets (e.g., $\calA$).  $2^\calA$ denotes $\calA$'s power set. $|\calA|$ its cardinality. $\calA\setminus\calB$ denotes set difference: the elements in $\calA$ not in~$\calB$. 
{Given a set function {$f: 2^{\calV_1}\cup \cdots \cup 2^{\calV_T}\mapsto \mathbb{R}$}, and $\calA_1\subseteq \calV_1,\ldots$, $ \calA_t\subseteq \calV_t$, for some positive integer $t \leq T$, the  $f(\calA_1,\ldots,\calA_t)$ denotes $f(\calA_1\cup\cdots\cup\calA_t\cup\emptyset\cup\cdots\cup\emptyset)$ where the $\emptyset$ is repeated $T-t$ times, and $\emptyset$ denotes the empty set.}
$\mathbb{R}$ denotes the set of real numbers.
}
\begin{align}\label{eq:non_res}
\begin{split}
\max_{\calA_1 \subseteq \calV_1}&\cdots
\max_{\calA_T \subseteq \calV_T} \; f(\calA_{1},\ldots,\calA_T), \\
&\;\text{s.t. } \;\;\;|\calA_t|\;= \alpha_t, \;\;\; t=1,\ldots,T,
\end{split}
\end{align}
where $T$ is a {given} horizon; $\calV_t$ is a given finite set of available elements to choose from at $t$ such that $\calV_{t}\cap \calV_{t'}=\emptyset$ for all $t,t'=1,\ldots,T$;\footnote{Even if the elements in $\calV_1,\ldots,\calV_T$ correspond to the same system modules, e.g., sensors, the elements among different $\calV_t$ are differentiated because they are chosen at different times.  For example, consider the case where $T=2$, and  two sensors $s_1$ and $s_2$ are available  to be chosen  at each $t$; then, by denoting with $s_{i,t}$ that sensor $i$ is available to be chosen at $t$, it is $\calV_1=\{s_{1,1}, s_{2,1}\}$ and $\calV_2=\{s_{1,2},s_{2,2}\}$, and, naturally, $\calV_1\cap\calV_2=\emptyset$.
} {$f: 2^{\calV_1}\cup \cdots \cup 2^{\calV_T}\mapsto \mathbb{R}$} is a given objective function;
$\alpha_t$ is a {given} cardinality constraint, capturing the resource constraints at~$t$; and $\calA_t$ are the chosen elements at $t$, resulting from the solution of eq.~\eqref{eq:non_res}.
Notably, in all above applications, and~\cite{gupta2006stochastic,jawaid2015submodularity,clark2016submodularity,tzoumas2016scheduling,zhang2017kalman,tzoumas2018codesign,zhao2016scheduling,nozari2017time,ikeda2018sparsity,leskovec2007cost,poularakis2017one,summers2014submodularity}, $f$ is non-decreasing, and without loss of generality one may consider $f(\emptyset)=0$.
For example, in~\cite{zhao2016scheduling}, $f$ is the trace of the inverse of the controllability Gramian, which captures the average control effort for driving the system; and in~\cite{tzoumas2016scheduling}, $f$ is the logdet of the error covariance of the minimum mean square batch-state estimator. Specifically, in~\cite{tzoumas2016scheduling}, $f$ is also submodular, a diminishing returns property that captures the intuition that a sensor's contribution to $f$'s \mbox{value diminishes when more sensors are activated already.}

Although the problem in eq.~\eqref{eq:non_res} is computationally hard, efficient algorithms have been proposed for its solution: when~$f$ is monotone and submodular, then eq.~\eqref{eq:non_res} is NP-hard~\cite{Feige:1998:TLN:285055.285059} and the greedy algorithm in~\cite[Section~4]{fisher1978analysis} guarantees a constant suboptimality bound across all problem instances; and when $f$ is only monotone, then eq.~\eqref{eq:non_res} is inapproximable (no polynomial time algorithm guarantees a constant bound across all instances)~\cite{foster2015variable,ye2018complexity} 
but the greedy algorithm in~\cite{fisher1978analysis} guarantees {per-instance} bounds instead~\cite{das2011spectral,wang2016approximation,sviridenko2017optimal}.

In this paper, however, we shift focus to a novel reformulation of eq.~\eqref{eq:non_res} that is robust against failures/attacks. Particu- larly, in all above applications, at any time $t$, actuators can be cyber-attacked~\cite{abdelzaher2018toward}, sensors can malfunction~\cite{wood2002dos}, and communication channels can be blocked~\cite{yang2019scheduling}, all resulting to denial-of-service (DoS) failures, {in the sense that the actuators, sensors, channels, etc. will shut down (stop working), at least temporarily}.  Hence, in such failure-prone and adversarial scenarios, eq.~\eqref{eq:non_res} may fail to protect any of the above applications, since it ignores the possibility of DoS failures. Thus, towards guaranteed protection, a robust reformulation becomes necessary \textit{that can both adapt to the history of incured failures and account for future ones}.

Therefore, in this paper we introduce a novel robust optimization framework, named \textit{Robust Sequential submodular Maximization} (\rsm), that goes beyond  the failure-free eq.~\eqref{eq:non_res} and accounts for  DoS failures/attacks. Specifically, we define \rsm 
as the following robust reformulation of eq.~\eqref{eq:non_res}:
\begin{align}\label{pr:robust_sub_max} 
\begin{split}
\text{\textbf{\rsm problem:}}&\\
\max_{\calA_1 \subseteq \calV_1}\min_{\;\calB_1 \subseteq \calA_1}&\!\cdots\!
\max_{\calA_T \subseteq \calV_T}\min_{\;\calB_T \subseteq \calA_T} \; \;f(\calA_{1}\setminus \calB_1,\ldots,\calA_T\setminus \calB_T), \\
&\;\text{s.t. } \;\;\;|\calA_t|\;= \alpha_t, \;\; |\calB_t|\;\leq \beta_t, \;\;\; t=1,\ldots,T.
\end{split}
\end{align}
where $\beta_t$ is a given number of possible failures (generally, $\beta_t\in [0, \alpha_t]$); and $\calB_t$ is the failure against~$\calA_t$. 

By solving \rsm, our goal is to 
maximize $f$ despite {worst-case} failures that occur at each maximization step, as captured by the intermediate/subsequent minimization steps.  Evidently, since \rsm considers worst-case failures, it is suitable when there is no prior on the failure mechanism, or when protection against worst-case failures is essential, such as in safety-cri- tical target tracking and costly experiment designs.

\rsm  can be interpreted as a $T$-stage {perfect information} sequential game between a ``maximization'' player (defender) and a ``minimization'' player (attacker)~\cite[Chapter~4]{myerson2013game}.  The defender starts the game and the players act sequentially, {having perfect knowledge of each others' actions}: at each $t$, the defender selects an $\calA_t$, and then the attacker responds with a worst-case removal $\calB_t$ from $\calA_t$, while both players account for the history of all actions up to $t-1$. In this context, the defender finds an optimal sequence $\calA_1,\ldots,\calA_{T}$ by accounting at each $t$ (i) for the history of responses $\calB_1,\ldots,\calB_{t-1}$, (ii) for the subsequent response $\calB_t$, and (iii) for all remaining future responses $\calB_{t+1},\ldots,\calB_{T}$. 
This is an additional computational challenge in comparison to the failure-free eq.~\eqref{eq:non_res}, which is already computationally hard. 

No scalable algorithms exists for \rsm.
\textit{In this paper, to provide the first scalable algorithm, we  develop an adaptive algorithm} that at each $t$ accounts only (i) for the history of responses up to $t-1$ and (ii) for the subsequent response $\calB_t$ (but not for the remaining future responses up to $t=T$), and as a result is scalable, but which still can guarantee a performance close to the optimal.

\myParagraph{Related work in combinatorial optimization} 
The majority of the related work has focused on the failure-free  eq.~\eqref{eq:non_res}, when~$f$ is either monotone and submodular
or  only monotone.
In more detail, Fisher et al.~\cite{fisher1978analysis} focused on $f$ being monotone and submodular, and proposed offline and online greedy algorithms that both guarantee the {constant} $1/2$ suboptimality bound. 
Similarly, Conforti and Cornu\'{e}jols~\cite{conforti1984curvature}, Iyer et al.~\cite{iyer2013curvature}, and Sviridenko et al.~\cite{sviridenko2017optimal} focused again on $f$ being monotone and sumodular but provided instead {per-instance},  curvature-depended  bounds.  The bounds generally tighten the ones in~\cite{fisher1978analysis}.   Finally,  Krause et al.~\cite{krause2010submodular}, Das and Kempe~\cite{das2011spectral}, Wang et al.~\cite{wang2016approximation}, and Sviridenko et al.~\cite{sviridenko2017optimal} (see also the earlier~\cite{sviridenko2013optimal}) focused on $f$ being only monotone, and proved per-instance, curvature-depended bounds for the greedy algorithms in~\cite{fisher1978analysis}, using notions of curvature ---also referred to as ``submodularity ratio''--- they introduced.

Recent work has also studied failure-robust reformulations of eq.~\eqref{eq:non_res}, typically per \rsm's framework \textit{but only for $T=1$, where no adaptiveness is required}.
Specifically, when $f$ is monotone and submodular, Orlin et al.~\cite{orlin2016robust} and Bogunovic et al.~\cite{bogunovic2017robust} provided greedy algorithms with constant suboptimality bounds. However, the algorithms are valid only for limited numbers of failures (for $\beta_1\leq \sqrt{\alpha_1}$ in~\cite{orlin2016robust} and $\beta_1\leq \alpha_1/(\log\alpha_1)^3$ in~\cite{bogunovic2017robust}).  In contrast, Tzoumas et al.~\cite{tzoumas2017resilient} provided a greedy algorithm with per-instance bounds for any number of failures ($\beta_1$ can take any value in $[0,\alpha_1]$).  Also, Rahmattalabi et al.~\cite{rahmattalabi2018robust} developed a mixed-integer linear program approach for a locations monitoring problem. More recently, Tzoumas et al.~\cite{tzoumas2018resilient-non-sub} and Bogunovic et al.~\cite{bogunovic2018robust} extended the previous works on the $T=1$ case by focusing on $f$ being only monotone, and proved per-instance, curvature-dependent bounds for the algorithm introduced in~\cite{tzoumas2017resilient}.  In more detail, Bogunovic et al.~\cite{bogunovic2018robust} focuses on cardinality constraints, whereas Tzoumas et al.~\cite{tzoumas2018resilient-non-sub} on the more general matroid constraints, {but, still, for the case where $T=1$}.  The latter framework enabled applications of failure-robust multi-robot robot planning, and particularly of active information gathering~\cite{schlotfeldt2018resilient} and target tracking~\cite{zhou2018resilient}.  Other relevant work is that of Mitrovic et al.~\cite{mitrovic2017streaming}, where a memoryless failure-robust reformulation of eq.~\eqref{eq:non_res} is considered, instead of the sequential framework of \rsm, which takes into account the history of past selections/failures.  Finally, Mirzasoleiman et al.~\cite{mirzasoleiman2017deletion} and Kazemi et al.~\cite{kazemi2017deletion} adopted a robust optimization framework against non worst-case failures, in contrast to \rsm which is against worst-case failures.\\
\indent All in all, in comparison to all prior research, in this paper we analyze \rsm's multistep case $T>1$ for the first time, and consider adaptive algorithms. 

\myParagraph{Related work in control} In the robust/secure control literature, various approaches have been proposed towards fault-tolerant control, secure control, as well as secure state estimation, against random failures, data injection and DoS failures/attacks~\cite{blanke2006diagnosis,8758093,liu2011false,jin2018power,clark2018linear,zhu2015game,jin2017adaptive,pasqualetti2013attack,perelman2016sensor,shoukry2017secure,pajic2014robustness,8710309,mo2015physical,8795564,de2014resilient,su2018finite,chen2018resilient,schenato2007foundations,amin2009safe,lu2017input}. 
In contrast to \rsm's resource-constrained  framework,~\cite{blanke2006diagnosis,8758093,liu2011false,jin2018power,clark2018linear,zhu2015game,jin2017adaptive,pasqualetti2013attack,perelman2016sensor,shoukry2017secure,pajic2014robustness,8710309,mo2015physical,8795564,de2014resilient,su2018finite,chen2018resilient,schenato2007foundations,amin2009safe,lu2017input} focus in resource abundant environments where all sensors and actuators stay always active under normal operation. 
For example,~\cite{schenato2007foundations,amin2009safe,lu2017input} focus on DoS failures/attacks from the perspective of packet loss and intermittent network connectivity, which can result to system destabilization. Generally,~\cite{blanke2006diagnosis,8758093,liu2011false,jin2018power,clark2018linear,zhu2015game,jin2017adaptive,pasqualetti2013attack,perelman2016sensor,shoukry2017secure,pajic2014robustness,8710309,mo2015physical,8795564,de2014resilient,su2018finite,chen2018resilient,schenato2007foundations,amin2009safe,lu2017input} focus on failure/attack detection and identification, and/or secure estimator/controller design, \textit{instead of the adaptive activation of a few sensors/actuators against worst-case DoS failures/attacks per \rsm.}

\myParagraph{Contributions} We introduce the novel \rsm problem of robust sequential maximization against DoS failures/attacks. 
We develop the first scalable algorithm, named \textit{\alglongCap} (\alg), that has the properties:

\setcounter{paragraph}{0}
\paragraph*{{\small\textbullet}~Adaptiveness} At each time $t=1,2,\ldots$, \alg selects a robust solution $\calA_t$ in an online fashion, accounting for the history of failures $\calB_1,\ldots,\calB_{t-1}$ and of actions $\calA_1,\ldots,\calA_{t-1}$, as well as, for all possible failures at $t$ from $\calA_t$.
\paragraph*{{\small\textbullet}~System-wide robustness} \alg is valid for any number of failures; that is, for any $\beta_t \in [0,\alpha_t]$, $t=1,2,\ldots$.
\paragraph*{{\small\textbullet}~Polynomial running time} \alg has the same order of running time as the polynomial time greedy algorithm proposed in~\cite[Section~4]{fisher1978analysis} for the failure-free eq.~\eqref{eq:non_res}.
\paragraph*{{\small\textbullet}~Provable approximation performance} \alg has provable {per-instance} suboptimality bounds that quantify \alg's near-optimality at each problem instance at hand.\footnote{Similarly to eq.~\eqref{eq:non_res}, \rsm is generally inapproximable: no polynomial time algorithm guarantees a constant suboptimality bound across all problem instances.  
	For example, it is inapproximable for fundamental applications in control and machine learning 
	such as \textit{sensor selection for optimal Kalman filtering}~\cite{ye2018complexity}, and \textit{feature selection for sparse model fitting}~\cite{foster2015variable}. Thus, in this paper we focus our analysis in per-instance suboptimality bounds.} 
Particularly, we provide both {a priori} and {a posteriori} per-instance bounds.  The a priori bounds quantify \alg's near-optimality {before} \alg has run.
In contrast,
the a posteriori bounds are computable {online} (as \alg runs), {once} the failures at each current step have been observed.  
The a posteriori bounds are tight and/or optimal.\footnote{{A suboptimality bound is called \emph{optimal} when it is the tightest achievable bound among all polynomial time algorithms, given a worst-case family of~$f$.}}  Finally, we present approximations of the a posteriori bounds
 that are computable {before} each failure occurs.
To quantify the bounds, we use curvature notions by Conforti and Conru\'{e}jols~\cite{conforti1984curvature}, for monotone and submodular functions, and  {Sviridenko et al.~\cite{sviridenko2017optimal}, for monotone functions.}

We demonstrate \alg's effectiveness in applications of {sensor scheduling}, and of {target tracking with wireless sensor networks}.
We present a Monte Carlo analysis, where we vary the failure types from worst-case to greedily and randomly selected failures, and compare \alg against a brute-force optimal algorithm (viable only for small-scale instances), the greedy algorithm in~\cite{fisher1978analysis}, and a random algorithm.
In the results, we observe \alg's near-optimality against worst-case failures, its robustness against non worst-case failures, and its superior performance against the compared algorithms.

{\textbf{Comparison with the preliminary results in~\cite{tzoumas2018resilient}, which coincides with preprint~\cite{tzoumas2018arxiv-resilient}:} 
	This paper extents the results in~\cite{tzoumas2018resilient}, considers new simulations, and includes the proofs that were all omitted from~\cite{tzoumas2018resilient}.  Particularly, most of the technical results herein, including Theorem~\ref{th:aposteriori}, Theorem~\ref{th:tightness_optimality}, Corollary~\ref{cor:ineq_from_lemmata}, and Algorithm~\ref{alg:bisection}, are novel and have not been previously published. 
	Also, the simulation scenarios are new and include a sensitivity analysis of \alg against various failure types (in~\cite{tzoumas2018resilient} we tested \alg only against worst-case failures, and in different scenarios). 
	Finally, all proofs in~\cite{tzoumas2018resilient} were omitted and are now included here.  
}

\myParagraph{Organization of the rest of the paper} Section~\ref{sec:algorithm} presents \alg, and quantifies its minimal running time.  Section~\ref{sec:performance} presents \alg's suboptimality bounds.  Section~\ref{sec:simulations} presents \alg's numerical evaluations. Section~\ref{sec:con} concludes the paper. {All proofs are found in the appendix.}

\section{An Adaptive Algorithm: \alg} \label{sec:algorithm}

We present \alg, the first scalable algorithm for \rsm, formulated in eq.~\eqref{pr:robust_sub_max}. 
\alg's pseudo-code is given in Algorithm~\ref{alg:rob_sub_max}.  
Below, we first give an intuitive description of \alg, and then a step-by-step description.  Also, we quantify its running time.  \alg's suboptimality bounds are given in Section~\ref{sec:performance}.

\begin{algorithm}[t]
	\caption{Robust adaptive maximization (\alg). 
	}
	\begin{algorithmic}[1]
		\REQUIRE \alg receives the inputs:
		\begin{itemize}[leftmargin=-.1cm]
			\item \textit{Offline}:~integer $T$; function {$f\!:\!2^{\mathcal{V}_1}\!\cup\! \cdots\!\cup\! 2^{\mathcal{V}_T}\! \mapsto \mathbb{R}$ }such that~$f$ is non-decreasing and $f(\emptyset)=0$; integers $\alpha_t$, $\beta_t$ \\ such that~$0\leq\beta_t \leq \alpha_t \leq |\mathcal{V}_t|$, for all $t=1,\ldots,T$;
			\item \textit{Online}:~at each $t=2,\ldots,T$, {observed} removal $\calB_{t-1}$ from \alg's output $\mathcal{A}_{t-1}$.
		\end{itemize}
		\ENSURE  At each step $t=1,\ldots,T$\!, set $\mathcal{A}_{t}$.
		\medskip
		
		\FORALL {$t=1,\ldots,T$}\label{line:begin_for_1}
		\STATE $\mathcal{S}_{t,1}\leftarrow\emptyset$;~~~$\mathcal{S}_{t,2}\leftarrow\emptyset$;\label{line:initiliaze}
		\STATE Sort elements in $\mathcal{V}_t$ s.t.~$\mathcal{V}_t\equiv\{v_{t,1}, \ldots, v_{t,|\calV_t|}\}$
		and $f(v_{t,1})\geq \ldots \geq f(v_{t,|\calV_t|})$;\label{line:sort}
		\STATE $\mathcal{S}_{t,1}\leftarrow\{v_{t,1}, \ldots, v_{t,\beta_t}\}$; \label{line:bait}
		\WHILE {$|\mathcal{S}_{t,2}|\; < \alpha_t-\beta_t$} \label{line:begin_while} 
		\STATE $x\!\in\! \arg\max_{y \in \mathcal{V}_t\setminus (\mathcal{S}_{t,1}\cup\mathcal{S}_{t,2})}f(\calA_{1}\!\setminus\!\calB_{1},\ldots,\calA_{t-1}\!\setminus\calB_{t-1}, \mathcal{S}_{t,2}\cup \{y\})$; \label{line:greedy_1}
		\STATE $\mathcal{S}_{t,2} \leftarrow \{x\}\cup \mathcal{S}_{t,2}$;\label{line:greedy_2}
		\ENDWHILE \label{line:end_while}
		\STATE $\mathcal{A}_{t}\leftarrow \mathcal{S}_{t,1} \cup \mathcal{S}_{t,2}$. \label{line:selection}
		\ENDFOR \label{line:end_for_3}
		
	\end{algorithmic}\label{alg:rob_sub_max}
\end{algorithm}

\subsection{Intuitive description}\label{subsec:intuition}

\rsm aims to maximize $f$ through a sequence of steps despite compromises to each step.  Specifically, at each $t=1,2,\ldots$, \rsm selects an $\calA_t$ towards a maximal $f$ despite the fact that $\calA_t$ will be compromised by a worst-case removal $\calB_t$, resulting to $f$ being evaluated at $\calA_1\setminus \calB_1,\ldots, \calA_T\setminus \calB_T$ instead of $\calA_1,\ldots,\calA_T$.
In~this~context, \alg aims to achieve \rsm's goal by selecting $\calA_t$ as the union of two sets $\calS_{t,1}$, and $\calS_{t,2}$ (\alg's line~\ref{line:selection}), whose role we describe intuitively below:

\setcounter{paragraph}{0} 

\paragraph{$\calS_{t,1}$ approximates ({aims to guess the}) worst-case removal from~$\calA_{t}$}  With $\calS_{t,1}$, \alg aims to capture the worst-case removal of $\beta_t$ elements from $\calA_t$.  Intuitively, $\calS_{t,1}$ is aimed to act as a ``bait'' to a worst-case attacker that selects the {best}~$\beta_t$ elements to remove from~$\calA_{t}$ at time $t$ (\textit{best} with respect to their contribution towards \rsm's goal). 
\alg  { aims to approximate} them by letting $\calS_{t,1}$ be the set of~$\beta_t$ elements with the largest marginal contributions to $f$ (\alg's lines~\ref{line:sort}-\ref{line:bait}).  {As such, each $\calS_{t,1}$ is independent of the history of actual removals $\calB_1,\ldots, \calB_{t-1}$ and can be computed offline, before any of the $\calB_1,\ldots, \calB_T$ has been realized.  In contrast, $\calS_{t,2}$ can only be computed online, as we describe below.}

\paragraph{$\calS_{t,1}\cup \calS_{t,2}$ approximates optimal solution to  \rsm's $t$-th  step}
To complete $\calA_t$'s construction, \alg needs to select a set $\calS_{t,2}$ of $\alpha_t-\beta_t$ elements (since $|\calA_t|=\alpha_t$ and $|\calS_{t,1}|=\beta_t$), and return $\calA_{t}=\calS_{t,1}\cup \calS_{t,2}$ (\alg's line~\ref{line:selection}).  Assuming $\calS_{t,1}$'s removal from $\calA_t$, for $\calA_{t}$ to be an optimal solution to \rsm's $t$-th maximization step, \alg needs to select $\calS_{t,2}$ as a \textit{best} set of $\alpha_t-\beta_t$ elements from $\calV_t\setminus\calS_{t,1}$.  
Nevertheless, this problem
is NP-hard~\cite{Feige:1998:TLN:285055.285059}.  Thereby, \alg  { approximates} such a best set, 
using the greedy procedure in \alg's lines~\ref{line:begin_while}-\ref{line:end_while}.  {Particularly, \alg's line~6 adapts $\calS_{t,2}$ to the history of removals $\calB_1,\ldots, \calB_{t-1}$ and  selections $\calA_1,\ldots, \calA_{t-1}$, since it constructs $\calS_{t,2}$ given $\calA_1\setminus\calB_1,\ldots, \calA_{t-1}\setminus\calB_{t-1}$.  As such, each $\calS_{t,2}$, in contrast to $\calS_{t,1}$, can be computed only online, only once the history of removals $\calB_1,\ldots, \calB_{t-1}$ has been~realized.}   

Overall, \alg { adaptively} constructs an $\calA_t$ to approximate an optimal solution to \rsm's $t$-th maximization step. 

{\begin{remark}[Further intuition on why $\calS_{t,1}$ and $\calS_{t,2}$ are selected as in \alg]
We first discuss why \alg  (i) selects $\calA_{t}$ as the union of $\calS_{t,1}$ and $\calS_{t,2}$, and (ii) selects $\calS_{t,2}$ as a greedily picked subset of $\calV_t\setminus\calS_{t,1}$.  We then focus on $\calS_{t,1}$.

If $\calS_{t,1}$ guesses correctly the removal $\calB_t$ from $\calA_{t}=\calS_{t,1}\cup \calS_{t,2}$, then all elements in $\calS_{t,2}$ remain intact \emph{($\calA_t\setminus \calB_t = \calS_{t,2}$)}.  Therefore, since $\calS_{t,2}$ has been selected using the greedy algorithm in \alg's lines~\ref{line:begin_while}-\ref{line:end_while}, which is an  {optimal approximation} algorithm for maximizing monotone functions subject to cardinality constraints~\cite{sviridenko2017optimal},\footnote{{An approximation algorithm is called \emph{optimal} when it achieves the tightest possible achievable suboptimality bound among all polynomial time algorithms, given a worst-case family of functions $f$.}} $\calA_{t}$ is an {optimal} approximation to \rsm's $t$-th maximization step.  This explains why \alg  selects $\calA_{t}$ as the union of two sets ($\calS_{t,1}$ and $\calS_{t,2}$), and why \alg selects $\calS_{t,2}$ greedily from $\calV_t\setminus\calS_{t,1}$ given $\calS_{t,1}$. 
Generally, if  $\calS_{1,1},\ldots, \calS_{T,1}$ guess  $\calB_1,\ldots,\calB_T$ correctly, then \rsm becomes equivalent to the attack-free eq.~\eqref{eq:non_res}, and \alg becomes equivalent to the optimal  greedy algorithm for eq.~\eqref{eq:non_res}.

But if $\calS_{t,1}$ guesses incorrectly $\calB_t$,  then some of $\calS_{t,1}$'s elements will survive, and, instead, some of $\calS_{t,2}$'s elements will be removed.  The question  arising now is: Can the elements that survived in $\calS_{t,1}$ compensate for the removed elements  from $\calS_{t,2}$?  In this paper, we develop tools to prove that if the elements of $\calS_{t,1}$  are chosen as in \alg's lines~\ref{line:sort}-\ref{line:bait}, this is indeed the case (proof of Theorem~\ref{th:alg_rob_sub_max_performance}), providing the first provable approximation guarantees for \rsm via \alg.
\end{remark}
}

\subsection{Step-by-step description}
 
\alg executes four steps for each $t=1,\ldots,T$: 
\setcounter{paragraph}{0}
\paragraph{Initialization (\alg's line~\ref{line:initiliaze})} \alg defines two auxiliary sets, namely, $\calS_{t,1}$ and $\mathcal{S}_{t,2}$, and initializes them with the empty set (\alg's line~\ref{line:initiliaze}). 

\paragraph{Construction of set $\calS_{t,1}$ (\alg's lines~\ref{line:sort}-\ref{line:bait})} \alg constructs $\calS_{t,1}$
by selecting $\beta_t$ elements, among all $s \in \calV_t$, with the highest values $f(s)$.  In detail, $\calS_{t,1}$ is constructed by first indexing the elements in~$\calV_t$ such that  $\mathcal{V}_t\equiv\{v_{t,1},$ $\ldots, v_{t,|\calV_t|}\}$
and $f(v_{t,1})\geq \ldots \geq f(v_{t,|\calV_t|})$ (\alg's line~\ref{line:sort}), and then by including in $\calS_{t,1}$ the fist $\beta_t$ elements (\alg's line~\ref{line:bait}).

\paragraph{Construction of set $\calS_{t,2}$ (\alg's lines~\ref{line:begin_while}-\ref{line:end_while})} \alg constructs  $\calS_{t,2}$ by picking greedily $\alpha_t-\beta_t$ elements from $\calV_t\setminus \calS_{t,1}$,
taking also into account the history of selections and removals, that is, $\calA_{1}\setminus\calB_1,\ldots, \calA_{t-1}\setminus\calB_{t-1}$.
Specifically, the ``while loop''  (\alg's lines~\ref{line:begin_while}-\ref{line:end_while}) selects an element $y\in\mathcal{V}_t\setminus (\mathcal{S}_{t,1}\cup\mathcal{S}_{t,2})$ to add in $\calS_{t,2}$ only if $y$ maximizes the value of  $f(\calA_{1}\setminus\calB_1,\ldots, \calA_{t-1}\setminus\calB_{t-1}, \mathcal{S}_{t,2}\cup \{y\})$.

\paragraph{Construction of set $\calA_{t}$ (\alg's line~\ref{line:selection})} 
\alg constructs $\calA_t$ as the union of $\calS_{t,1}$ and~$\mathcal{S}_{t,2}$.  

The above steps are valid for any number of failures $\beta_t$.

\subsection{Running time}

We now analyze the computational complexity of \alg.

\begin{proposition}\label{prop:runtime}
	At each $t=1,2,\ldots$, \alg runs in $O[|\calV_t|(\alpha_t-\beta_t)\tau_f]$ time, where $\tau_f$ is $f$'s evaluation time. 
\end{proposition}

\begin{remark}[Minimal running time]
Even though \alg robustifies the traditional, failure-free sequential optimization in eq.~\eqref{eq:non_res}, \alg has the same order of running time as the state-of-the-art algorithms for eq.~\eqref{eq:non_res}~\cite[Section~4]{fisher1978analysis}~\cite[Section~8]{sviridenko2017optimal}. 
\end{remark}

\smallskip

In summary, \alg selects adaptively a solution for \rsm, in minimal running time, and is valid for any number of failures.  We quantify its approximation performance next.

\section{Suboptimality Guarantees
 }\label{sec:performance}

We present \alg's suboptimality bounds.
We first present \alg's 
a priori bounds, 
and, then, the {a posteriori} bounds.  
Finally, we present the latter's {pre-failure} approximations.

\subsection{Curvature and total curvature}\label{sec:total_curvature}
 
To present \alg's suboptimality bounds we use the notions of \emph{curvature} and \emph{total curvature}. 
To this end, we start by recalling the definitions of \textit{modularity} and 
\textit{submodularity}, where we consider the notation:
\begin{itemize}
	\item $\calV\triangleq \bigcup_{i=1}^T\calV_t$; i.e., $\calV$ is the union   across the horizon $T$ of all the available elements to choose from;
\end{itemize}

\begin{definition}[Modularity~{\cite{nemhauser78analysis}}]\label{def:modular}
$f:2^\calV\mapsto \mathbb{R}$ is modular if and only if $f(\calA)=\sum_{\elem\in \calA}f(\elem)$, for any $\calA\subseteq \calV$.
\end{definition}

Therefore, if $f$ is modular, then $\calV$'s elements complement each other through $f$. Particularly, Definition~\ref{def:modular} implies $f(\{\elem\}\cup\calA)-f(\calA)= f(\elem)$, for any $\calA\subseteq\calV$ and $\elem\in \calV\setminus\calA$.

\begin{definition}[Submodularity~{\cite{nemhauser78analysis}}]\label{def:sub} $f:2^\calV\mapsto \mathbb{R}$ is \emph{submodular} if and only if 
$f(\calA\cup \{\elem\})-f(\calA)\geq f({\mathcal{B}}\cup \{\elem\})-f({\mathcal{B}})$, for any $\calA\subseteq {\mathcal{B}}\subseteq\calV$ and $\elem\in \calV$.
\end{definition}

The definition implies $f$ is submodular if and only if  the return $f(\calA\cup \{\elem\})-f(\calA)$ diminishes as $\calA$ grows, for any $\elem$. 
In contrast to $f$ being modular, if $f$ is submodular, then $\calV$'s elements substitute each other. Specifically, without loss of generality, consider $f$ to be non-negative: then, Definition~\ref{def:sub} implies $f(\{\elem\}\cup\calA)-f(\calA)\leq f(\elem)$. That is, in the presence of $\calA$, $\elem$'s contribution to $f(\{\elem\}\cup\calA)$'s~value is diminished.

\begin{definition}\label{def:curvature}
\emph{{(Curvature~\cite{conforti1984curvature})}}
Consider a non-decreasing submodular $f:2^\calV\mapsto\mathbb{R}$ such that $f(\elem)\neq 0$, for any $\elem \in \calV$,  without loss of generality.  Then, $f$'s \emph{curvature} is defined as
 \begin{equation}\label{eq:curvature}
\kappa_f\triangleq 1-\min_{\elem\in\calV}\frac{f(\calV)-f(\calV\setminus\{\elem\})}{f(\elem)}.
\end{equation}
\end{definition}

Definition~\ref{def:curvature} implies $\kappa_f \in [0,1]$.   Particularly, $\kappa_f$ measures how far~$f$ is from modularity: if $\kappa_f=0$, then  $f(\calV)-f(\calV\setminus\{v\})=f(v)$, for all $v\in\calV$; that is, $f$ is modular. In~contrast, if $\kappa_f=1$, then there exist $v\in\calV$ such that $f(\calV)=f(\calV\setminus\{v\})$; that is, $v$~has no contribution to $f(\calV)$ in the presence of $\calV\setminus\{v\}$.  Therefore, $\kappa_f$ can also been interpreted as a measure of how much $\calV$'s elements complement/substitute each other.

\begin{definition}[Total curvature~\cite{sviridenko2017optimal,sviridenko2013optimal}] \label{def:total_curvature}
Consider a monotone $f:2^\calV\mapsto\mathbb{R}$.  Then, $f$'s total curvature is defined as 
\begin{equation}\label{eq:total_curvature}
c_f\triangleq 1-\min_{v\in\calV}\min_{\calA, \mathcal{B}\subseteq \calV\setminus \{v\}}\frac{f(\{v\}\cup\calA)-f(\calA)}{f(\{v\}\cup\mathcal{B})-f(\mathcal{B})}.
\end{equation}
\end{definition}
Similarly to $\kappa_f$, it also is $c_f\in [0,1]$. Remarkably, when $f$ is  submodular, then $c_f=\kappa_f$. 
Generally, if  $c_f=0$, then $f$ is modular, while if $c_f=1$, then eq.~\eqref{eq:total_curvature} implies the assumption that $f$ is non-decreasing.
In~\cite{lehmann2006combinatorial}, any monotone $f$ with total curvature $c_f$ is called $c_f$-submodular, as repeated below.\footnote{Lehmann et al.~\cite{lehmann2006combinatorial}	defined $c_f$-submodularity by considering in eq.~\eqref{eq:total_curvature}  $\calA\subseteq\calB$ instead of $\calA\subseteq \calV$.  Generally, non submodular but monotone functions have been referred to as \textit{approximately} or \textit{weakly} submodular~\cite{krause2010submodular,elenberg2016restricted}, names that have also been adopted for the definition of $c_f$ in~\cite{lehmann2006combinatorial}, e.g., in~\cite{chamon2016near,guo2019actuator}.} 

\begin{definition}[$c_f$-submodularity~\cite{lehmann2006combinatorial}]\label{def:approx_sub} Any monotone function $f:2^\calV\mapsto\mathbb{R}$  with total curvature $c_f$ is called \emph{$c_f$-submodular}.
\end{definition}

\begin{remark}[Dependence on the {size of $\calV$} and length of horizon $T$]\label{rem:dependence}
{Evidently, both $\kappa_f$ and $c_f$  are non-decreasing as $\calV$ grows (cf.~Definition~\ref{def:curvature} and Definition~\ref{def:total_curvature}).} Therefore, $\kappa_f$ and $c_f$ are also non-decreasing as $T$ increases, since $\calV\equiv \bigcup_{i=1}^T\calV_t$.
\end{remark}

\begin{algorithm}[t]
	\caption{Online greedy algorithm~\cite[Section~4]{fisher1978analysis}.}
	\begin{algorithmic}[1]
		\REQUIRE 
		Integer $T>0$; {$f:2^{\mathcal{K}_1}\cup \cdots\cup 2^{\mathcal{K}_T} \mapsto \mathbb{R}$} such that \\ $f$ is non-decreasing and $f(\emptyset)=0$; integers \!$\delta_t$ \!such that $0\leq\delta_t \leq |\mathcal{K}_t|$, for all $t=1,\ldots,T$\!.
		\ENSURE  At each step $t=1,\ldots,T$\!, set $\mathcal{M}_{t}$.
		\medskip
		
		\FORALL {$t=1,\ldots,T$}
		\STATE $\mathcal{M}_{t}\leftarrow\emptyset$;
		\WHILE {$|\mathcal{M}_{t}| < \delta_t$} \label{line2:begin_while} 
		\STATE $x\!\in\!\arg\max_{y \in \mathcal{K}_t\setminus \mathcal{M}_{t}}\!f({\calM_{1},\ldots,\calM_{t-1},} \mathcal{M}_{t}\!\cup\!\{y\}\!)$; \label{line2:greedy_1}
		\STATE $\mathcal{M}_{t} \leftarrow \{x\}\cup \mathcal{M}_{t}$;\label{line2:greedy_2}
		\ENDWHILE \label{line2:end_while}
		\ENDFOR 
		
	\end{algorithmic}\label{alg:local}
\end{algorithm}

\subsection{A priori suboptimality bounds}

We present \alg's a priori suboptimality bounds, using the above notions of curvature.  We use also the notation:
\begin{itemize}
	\item $f^\star$ is the optimal value of \rsm;
	\item $\calA_{1:t}\triangleq (\calA_1,\ldots,\calA_t)$, where $\calA_{t}$ is the selected set by \alg at $t=1,2,\ldots$; 
	\item {$(\calB^\star_1,\ldots, \calB^\star_T)$ is an optimal removal from  $\calA_{1:T}$;}
	\item $\mathcal{B}^\star_{1:t}\triangleq (\calB^\star_1,\ldots, \calB^\star_t)$;
	\item $\calA_{1:t}\setminus \mathcal{B}^\star_{1:t}\triangleq (\calA_1\setminus\calB^\star_1,\ldots, \calA_t\setminus\calB^\star_t)$.
\end{itemize}

\begin{theorem}[A priori bounds]\label{th:alg_rob_sub_max_performance}
\alg selects $\calA_{1:T}$ such that  $|\calA_t|\;\leq \alpha_t$, and
if $f$ is submodular, then
\begin{equation}\label{ineq:bound_sub}
\frac{f(\calA_{1:T}\setminus \mathcal{B}^\star_{1:T})}{f^\star}\geq \left\{\begin{array}{lr}
\frac{1-{e^{-\kappa_f}}}{\kappa_f}(1-\kappa_f), & T=1;\\
(1-\kappa_f)^4, & T>1;
\end{array} \right. 
 \end{equation}
whereas, if $f$ is $c_f$-submodular, then
\begin{equation}\label{ineq:bound_non_sub}
\hspace*{-11mm}\frac{f(\calA_{1:T}\setminus \mathcal{B}^\star_{1:T})}{f^\star}\geq 
\left\{\begin{array}{lr}
 (1-c_f)^3, & T=1;\\
 (1-c_f)^5, & T>1.
\end{array} \right. 
\end{equation}
\end{theorem}

Evidently, Theorem~\ref{th:alg_rob_sub_max_performance}'s bounds are {a priori}, since ineqs.~\eqref{ineq:bound_sub}'s and~\eqref{ineq:bound_non_sub}'s right-hand-sides are independent of the selected $\calA_{1:T}$ by \alg, and the incurred failures~$\calB^\star_{1:T}$. 

Importantly, the bounds compare \alg's selection $\calA_{1:T}$  against an {optimal one that knows a priori all future failures} (and achieves that way  the value $f^\star$).  Instead, \alg's has {no knowledge} of the future failures.  Within this challenging setting, Theorem~\ref{th:alg_rob_sub_max_performance} nonetheless implies: for functions $f$ with $\kappa_f<1$ or  $c_f<1$, \alg's selection $\calA_{1:T}$ is finitely close to the optimal, instead of arbitrarily suboptimal.
Indeed, then Theorem~\ref{th:alg_rob_sub_max_performance}'s bounds are non-zero.
We discuss functions with $\kappa_f<1$ or $c_f<1$ below, along with relevant~applications.

\begin{remark}[Functions with $\kappa_f<1, c_f<1$, and applications]\label{rem:functionsiwthboundedcurv}
Functions with $\kappa_f<1$ are the concave over modular functions~\cite[Section~2.1]{iyer2013curvature} and the $\log\det$ of~\mbox{positive-definite} matrices~\cite{sharma2015greedy}. Also, functions with $c_f<1$ are the support selection functions~\cite{elenberg2016restricted}, 
the average minimum square error of the Kalman filter (trace of error covariance)~\cite[Section~IV]{chamon2017mean}, and the LQG cost as a function of the active sensors~\cite[Theorem~4]{tzoumas2018codesign}.  
The aforementioned functions appear in control and machine learning applications such as feature selection~\cite{das2011spectral,khanna2017scalable}, and actuator and sensor scheduling~\cite{gupta2006stochastic,jawaid2015submodularity,clark2016submodularity,tzoumas2016scheduling,zhang2017kalman,zhao2016scheduling,nozari2017time,ikeda2018sparsity,chamon2017mean,tzoumas2018codesign}.  
\end{remark}

Evidently, when $\kappa_f$ and  $c_f$ tend to 0, then \alg becomes optimal, 
since all bounds in Theorem~\ref{th:alg_rob_sub_max_performance} tend to $1$; for example, {{$1/\kappa_f(1-e^{-\kappa_f})(1-\kappa_f)$ increases as $\kappa_f$ decreases, and its limit is equal to $1$ for $\kappa_f\rightarrow0$}}.  Application examples of this sort involve the regression of Gaussian processes with RBF kernels~\cite[Theorem~5]{sharma2015greedy}, such as in sensor selection for temperature monitoring~\cite{krause2008near}.

{Finally, since both $\kappa_f$ and $c_f$ are non-decreasing in $T$ {and $\calV$} (Remark~\ref{rem:dependence}), the bounds are non-increasing in $T$ {and~$\calV$}.}

\paragraph*{Tightness and optimality (towards a posteriori bounds)} \alg's curvature-dependent bounds are the first suboptimality bounds for \rsm, and make a first step towards separating
the classes of monotone functions into
functions for which \rsm
can be approximated well (low curvature functions), and functions for which it cannot \mbox{(high curvature functions).} 
Moreover, although for the failure-free eq.~\eqref{eq:non_res} the {a priori} bounds $1/\kappa_f(1-e^{-\kappa_f})$ and $1/(1+\kappa_f)$ (where $f$ is submodular) are known to be tight~\cite[Theorem~2.12, Theorem~5.4]{conforti1984curvature}, the tightness of ineq.~\eqref{ineq:bound_sub} is an open problem.  Similarly, although for eq.~\eqref{eq:non_res} the {a priori} bound $1-c_f$ (where $f$ is $c_f$-submodular) is known to be optimal ({the tightest possible in polynomial time in a worst-case})~\cite[Theorem~8.6]{sviridenko2017optimal},
the optimality of  ineq.~\eqref{ineq:bound_non_sub} is an open problem. Notably, in the latter case ($f$ is $c_f$-submodular) both $1-c_f$ and the bound in ineq.~\eqref{ineq:bound_non_sub}  are $0$ for $c_f=1$, which is in agreement with the inapproximability of both eq.~\eqref{eq:non_res} and \rsm in the worst-case.  

In contrast to Theorem~\ref{th:alg_rob_sub_max_performance}'s  {a priori} bounds, we next present tight and/or optimal {a posteriori}~bounds.

\subsection{A posteriori suboptimality bounds
 }\label{sec:post_performance}

We now present \alg's a posteriori bounds, which are computable once all failures up to step $t$ have been observed.  
Henceforth, we use the~notation:
\begin{itemize}
	\item $f^\star_t$ is the optimal value of \rsm for $T=t$;
	\item $\calM_t$ is the set returned by the online, failure-free greedy Algorithm~\ref{alg:local} at $t=1,\ldots,T$,\footnote{{We refer to Algorithm~2 as ``online'' since each $\calM_t$ can be chosen in real time (at time $t$) sequentially, i.e., given the history of past selections $\calM_1,\ldots,\calM_{t-1}$.  Observe, however, that if one wishes so one  can also execute all steps of Algorithm~2 offline at time $t=0$.}} when we consider therein $\delta_t=\alpha_t-\beta_t$ and $\calK_t=\calV_t\setminus \calS_{1,t}$;

	\item $\calM_{1:t}\triangleq \{\calM_1, \ldots, \calM_t\}$.
\end{itemize}

\begin{remark}[Interpretation of $\calM_{1:t}$]
	Since each $\calS_{1,t}$ is the expected future failures (``baits'') selected in \alg's lines~\ref{line:sort}-\ref{line:bait} (see Section~\ref{sec:algorithm}), $\calM_{1:t}$ are the sets one would greedily select per Algorithm~\ref{alg:local} if it was known a priori that indeed the future failures are the $\calS_{1,t}$, $t=1,\ldots,T$.
\end{remark}

\begin{theorem}[A posteriori  bounds]\label{th:aposteriori}
	For all $t=1,\ldots, T$, {given the observed history $\calB^\star_{1:t}$,} \alg selects $\calA_{t}$ such that  $|\calA_t|\;\leq \alpha_t$, and
		if $f$ is submodular, then
		\begin{equation}\label{ineq:a_posteriori_bound_sub}
		\hspace*{-.70mm}\frac{f(\calA_{1:t}\setminus \mathcal{B}^\star_{1:t})}{f^\star_t}\geq\left\{\begin{array}{lr}
		\frac{1-{e^{-\kappa_f}}}{\kappa_f}
		\frac{f(\calA_{1}\setminus \mathcal{B}^\star_{1})}{f(\calM_{1})}, & t=1;\\
		\frac{1}{1+\kappa_f}
		\frac{f(\calA_{1:t}\setminus \mathcal{B}^\star_{1:t})}{f(\calM_{1:t})}, & t>1;
		\end{array} \right. 
		\end{equation}
		whereas, if $f$ is $c_f$-submodular, then
		\begin{equation}\label{ineq:a_posteriori_bound_non_sub}
		\hspace*{-2mm}\frac{f(\calA_{1:t}\setminus \mathcal{B}^\star_{1:t})}{f^\star_t}\geq(1-c_f)
		\frac{f(\calA_{1:t}\setminus \mathcal{B}^\star_{1:t})}{f(\calM_{1:t})}.
		\end{equation}
\end{theorem}

\begin{theorem}[Tightness and optimality]\label{th:tightness_optimality}
There exist families of $f$ such that the suboptimality bounds in ineq.~\eqref{ineq:a_posteriori_bound_sub} are tight. Also, there exist families of $f$ such that the suboptimality bounds in eq.~\eqref{ineq:a_posteriori_bound_non_sub} are optimal {(the tightest possible)} across all algorithms that evaluate $f$ a polynomial number of times.\footnote{Theorem~\ref{th:tightness_optimality}'s function families are the same as those in the proofs of~\cite[Theorem~2.12, Theorem~5.4]{conforti1984curvature} and~\cite[Theorem~8.6]{sviridenko2017optimal}, which prove the tightness and optimality of $1/\kappa_f(1-e^{-\kappa_f})$, $1/(1+\kappa_f)$  and $(1-c_f)$ for eq.~\eqref{eq:non_res}.}
\end{theorem}

The bounds break down into the a priori $\kappa_f$- and $c_f$-depended parts, and the a posteriori $f(\calA_{1:t}\setminus \mathcal{B}^\star_{1:t})/f(\calM_{1:t})$.  We refer to the latter as {a posteriori} since it is computable \textit{after} $\mathcal{B}^\star_{t}$ has been observed.  Intuitively, the a posteriori part captures how successful the  ``bait'' $\calS_{1,t}$ has been in approximating the anticipated worst-case failure $\calB^\star_t$. Indeed, if $\calB^\star_t=\calS_{1,t}$ for all $t=1,2,\ldots$, then $f(\calA_{1:t}\setminus \mathcal{B}^\star_{1:t})/f(\calM_{1:t})=1$ and  Theorem~\ref{th:aposteriori}'s bounds become the tight/optimal a priori bounds $1/\kappa_f(1-e^{-\kappa_f})$, $1/(1+\kappa_f)$ and $1-c_f$,\footnote{Theorem~\ref{th:tightness_optimality} is proved based on this observation.} {and, as such, they are also tighter than Theorem~\ref{th:alg_rob_sub_max_performance}'s a priori bounds.}

{In general, Theorem~\ref{th:aposteriori}'s a posteriori bounds may be looser than
 Theorem~\ref{th:alg_rob_sub_max_performance}'s a priori bounds; yet, they are tighter when $f(\calA_{1:t}\setminus \mathcal{B}^\star_{1:t})/f(\calM_{1:t})$ is close enough to 1: e.g., for $f$ being $c_f$-submodular and $T>1$, if  $f(\calA_{1:t}\setminus \mathcal{B}^\star_{1:t})/f(\calM_{1:t})>(1-c_f)^4$, then the a posteriori bound in eq.~\eqref{ineq:a_posteriori_bound_non_sub} is tighter than the a priori in eq.~\eqref{ineq:bound_non_sub}.  Indeed, in the numerical evaluations of Section~\ref{sec:simulations}, $f(\calA_{1:t}\setminus \mathcal{B}^\star_{1:t})/f(\calM_{1:t})$ is nearly 1, whereas $(1-c_f)^4\leq.0001$; thus, eq.~\eqref{ineq:a_posteriori_bound_non_sub} is $3$ orders tighter than eq.~\eqref{ineq:bound_non_sub}.
}

Notably, the a priori parts $1/\kappa_f(1-e^{-\kappa_f})$, $1/(1+\kappa_f)$ are non-zero for any values of $\kappa_f$. In more detail, $1/\kappa_f(1-e^{-\kappa_f})\geq 1-1/e$ and $1/(1+\kappa_f)\geq 1/2$ for all $\kappa_f \in [0,1]$; {particularly, $1/\kappa_f(1-e^{-\kappa_f})$ increases as $\kappa_f$ decreases, and its limit is equal to $1$ for $\kappa_f\rightarrow0$}. Therefore,
in contrast to the a priori bound in eq.~\eqref{ineq:bound_sub}, which for $\kappa_f=1$ becomes 0, eq.~\eqref{ineq:a_posteriori_bound_sub} for $\kappa_f=1$ becomes instead
\begin{equation}
\hspace*{-.70mm}\frac{f(\calA_{1:t}\setminus \mathcal{B}^\star_{1:t})}{f^\star_t}\geq\left\{\begin{array}{lr}
(1-1/e)
\frac{f(\calA_{1}\setminus \mathcal{B}^\star_{1})}{f(\calM_{1})}, & t=1;\\

\frac{f(\calA_{1:t}\setminus \mathcal{B}^\star_{1:t})}{2f(\calM_{1:t})}, & t>1.
\end{array} \right. 
\end{equation}

Nevertheless, such simplification for eq.~\eqref{ineq:a_posteriori_bound_non_sub} is not evident, a fact that is in agreement with both (i) \rsm's inapproximability when $f$ is \textit{not} submodular, necessitating~per-instance suboptimality bounds for any polynomial time algorithm, and (ii)  eq.~\eqref{ineq:a_posteriori_bound_non_sub}'s optimality per Theorem~\ref{th:tightness_optimality}.


\smallskip
Overall, Theorem~\ref{th:aposteriori}'s bounds are computable online, at each $t=1,2,\ldots$, \textit{after} failure $\calB^\star_t$ has been observed. 
%
%
%
We next approximate the bounds \textit{before}  $\calB^\star_t$ occurs. 

\subsection{Pre-failure approximations of \post bounds}\label{sec:pre_failure}

We present pre-failure approximations to Theorem~\ref{th:aposteriori}'s post-failure bounds. 
In particular, we propose a method to lower bound $f(\calA_{1:t}\setminus \mathcal{B}^\star_{1:t})$ by a value $\hat{f}_t$, at each $t=1,\ldots,T$,
{
given $\mathcal{B}^\star_{1:t-1}$ (but {before}  $\calB^\star_t$ occurs).
}

In more detail, we recall $f(\calA_{1:t}\setminus \mathcal{B}^\star_{1:t})$ is the value of the constrained optimization problem
\begin{equation}\label{eq:def_opt_attacks}
 f(\calA_{1:t}\setminus \mathcal{B}^\star_{1:t})\equiv \min_{\mathcal{B}_t \subseteq \calA_t,\; |\mathcal{B}_t| \leq \beta_t}\; f(\calA_{1:t-1}\setminus\calB_{1:t-1}^\star,\calA_t\setminus \mathcal{B}_t).
\end{equation}
Computing $f(\calA_{1:t}\setminus \mathcal{B}^\star_{1:t})$ is NP-hard, even if $f$ is submodular~\cite{schrijver2000combinatorial}.
But lower bounding $f(\calA_{1:t}\setminus \mathcal{B}^\star_{1:t})$ can be efficient.  Specifically, the non-constrained reformulation of eq.~\eqref{eq:def_opt_attacks} in eq.~\eqref{eq:regul_def_opt_attacks} below is efficiently solvable (see \cite{schrijver2000combinatorial,iwata2009simple,lee2015faster,chakrabarty2017subquadratic} for $f$ being submodular; and \cite{halabi2019minimizing} for $f$ being $c_f$-submodular):
\begin{equation}\label{eq:regul_def_opt_attacks}
\min_{\mathcal{B}_t \subseteq \calA_t} f(\calA_{1:t-1}\setminus\calB_{1:t-1}^\star,\calA_t\setminus \mathcal{B}_t)+\lambda_t |\mathcal{B}_t|,
\end{equation}
where $\lambda_t\geq 0$ and constant {($\lambda_t$ acts {similarly to} a Lagrange multiplier~\cite{boyd2004convex})}. Evidently, Lemma~\ref{lem:bisection} below holds true, where $\hat{\calB}_t(\lambda_t)$ denotes an optimal solution to eq.~\eqref{eq:regul_def_opt_attacks}, i.e.,
\begin{equation}\label{eq:b_regul_def_opt_attacks}
\hat{\calB}_t(\lambda_t)\in\arg \min_{\mathcal{B}_t \subseteq \calA_t} f(\calA_{1:t-1}\setminus\calB_{1:t-1}^\star,\calA_t\setminus \mathcal{B}_t)+\lambda_t |\mathcal{B}_t|,
\end{equation}
{and where $\hat{f}_t(\lambda_t)$ denotes the value of $f(\calA_{1:t-1}\setminus\calB_{1:t-1}^\star,\calA_t\setminus \mathcal{B}_t)$ when $ \mathcal{B}_t=\hat{\calB}_t(\lambda_t)$, i.e.,}
\begin{equation}\label{eq:def_f_hat}
{\hat{f}_t(\lambda_t) \triangleq f(\calA_{1:t-1}\setminus\calB_{1:t-1}^\star,\calA_t\setminus \hat{\calB}_t(\lambda_t)).}
\end{equation}

\begin{lemma}\label{lem:bisection}
There exists $\lambda_t^\star$ such that $\hat{f}_t(\lambda_t)\leq f(\calA_{1:t}\setminus \mathcal{B}^\star_{1:t})$ and $|\hat{\calB}_t(\lambda_t)|\; > \beta_t$  for $\lambda_t< \lambda_t^\star$; whereas, $\hat{f}_t(\lambda_t)\geq {f}(\calA_{1:t}\setminus \mathcal{B}^\star_{1:t})$ and $|\hat{\calB}_t(\lambda_t)|\; \leq  \beta_t$  for $\lambda_t\geq \lambda_t^\star$.
\end{lemma}

\begin{algorithm}[t]
	\caption{Bisection.}
	\begin{algorithmic}[1]
		\REQUIRE 
		{Integer $\beta_t$ per \rsm; function $f$ per \rsm; histories $\calA_{1:t}$ \!and\! $\calB^\star_{1:t-1}$; $u_0\!>\!0$ such that $|\hat{\calB}_t(u_0)| \;< \beta_t$, \!where $\hat{\calB}_t(\cdot)$ is defined in eq.~\eqref{eq:b_regul_def_opt_attacks}; {$\epsilon>0$, which defines bise- ction's stopping condition (accuracy level).}}
		\ENSURE $\lambda_t\geq 0$ such that  
		$\hat{f}_t(\lambda_t)\leq$ $f(\calA_{1:t}\setminus \mathcal{B}^\star_{1:t})$, where \\ $\lambda_t$ and $\hat{f}_t(\lambda_t)$ are defined in eq.~\eqref{eq:regul_def_opt_attacks} and  eq.~\eqref{eq:def_f_hat}.
		\medskip
		
		\STATE {$l \leftarrow 0$; \quad$u\leftarrow u_0$; \quad$\lambda_t\leftarrow (l+u)/2$;}
		\WHILE {{$u-l>\epsilon$}}\label{line3:begin-while}
		\STATE {Find $\hat{\calB}_t(\lambda_t)$ by solving eq.~\eqref{eq:b_regul_def_opt_attacks};}
		\IF {$|\hat{\calB}_t(\lambda_t)|\; < \beta_t$}
		\STATE{$u \leftarrow \lambda_t$;}  \COMMENT{\color{gray}$u$ always satisfies $|\hat{\calB}_t(u)|\; < \beta_t$\color{black}}\label{line3:u-update}
		\ELSE
		\STATE{$l \leftarrow \lambda_t$;} \COMMENT{\color{gray}$l$ always satisfies $|\hat{\calB}_t(l)|\; \geq \beta_t$\color{black}}\label{line3:l-update}
		\ENDIF
		\STATE {{$\lambda_t\leftarrow (l+u)/2$;}}
		\ENDWHILE\label{line3:end-while}
		\STATE $\lambda_t \leftarrow l$; \COMMENT{\color{gray}$l$ is $\epsilon$-close to $\lambda_t^\star$ \!($\lambda_t^\star$\! is\! defined \!in \!Lemma \ref{lem:bisection}) and satisfies $|\hat{\calB}_t(l)|\; \geq \beta_t$\color{black}}\label{line3:lambda-update1}
		\RETURN  $\lambda_t$.\label{line3:lambda-update2}
	\end{algorithmic}\label{alg:bisection}
\end{algorithm}

{To observe such a value $\lambda_t^\star$ exists, it suffices to observe: (i) for $\lambda_t=0$, the cardinality of $\hat{\calB}_t$ in eq.~\eqref{eq:regul_def_opt_attacks} is unconstrained, and, thus, the optimal solution in eq.~\eqref{eq:regul_def_opt_attacks} is to remove all $\calA_t$, which implies $|\hat{\calB}_t(\lambda_t)|\;=\alpha_t \geq \beta_t$; (ii)  more generally, for $\lambda_t>0$, the cardinality of $\hat{\calB}_t$ in eq.~\eqref{eq:regul_def_opt_attacks} is increasingly penalized as $\lambda_t$ increases, and, thus, $|\hat{\calB}_t(\lambda_t)|$ is a decreasing function of $\lambda_t$ (in particular, if $\lambda_t\rightarrow +\infty$, then $|\hat{\calB}_t(\lambda_t)|\;\rightarrow 0$).
	Now, given (i)-(ii), denote by $\lambda_t^\star$  the first value of $\lambda_t$ such that $|\hat{\calB}_t(\lambda_t)|\;\leq \beta_t$, when $\lambda_t$   is initially set to $0$ and then continuously increases: then, for  $\lambda_t< \lambda_t^\star$, it is $|\hat{\calB}_t(\lambda_t)|\; >  \beta_t$, and, thus, $\hat{f}_t(\lambda_t)\leq f(\calA_{1:t}\setminus \mathcal{B}^\star_{1:t})$, since $|\hat{\calB}_t(\lambda_t)| \; >\beta_t=|\calB^\star_t|$ and $\hat{\calB}_t(\lambda_t)$ is an optimal solution to eq.~\eqref{eq:regul_def_opt_attacks}; whereas, for  $\lambda_t\geq \lambda_t^\star$, it is $|\hat{\calB}_t(\lambda_t)|\; \leq  \beta_t$, and, thus, $\hat{f}_t(\lambda_t)\geq f(\calA_{1:t}\setminus \mathcal{B}^\star_{1:t})$.}

Although $\lambda^\star_t$ is unknown, it can be approximated by using bisection. For example, Algorithm~\ref{alg:bisection} uses bisection with accuracy level $\epsilon>0$ (Algorithm~\ref{alg:bisection}'s lines \ref{line3:begin-while}-\ref{line3:end-while}) to find a $\lambda_t$ that is $\epsilon$-close to $\lambda^\star_t$ and for which $\hat{f}_t(\lambda_t)\leq f(\calA_{1:t}\setminus \mathcal{B}^\star_{1:t})$. 
	  {To start the bisection, Algorithm~\ref{alg:bisection} assumes a large enough $u_0 \geq 0$ such that $|\hat{\calB}_t(u_0)| \;< \beta_t$; such a $u_0$ can be found since $|\hat{\calB}_t(u_0)|\;\rightarrow 0$ for $u_0\rightarrow +\infty$.  Next, at each ``while loop'' (lines~\ref{line3:begin-while}-\ref{line3:end-while} of Algorithm~\ref{alg:bisection}), $\lambda_t^\star \in [l,u]$, since  $|\hat{\calB}_t(u)|\; < \beta_t$ and $|\hat{\calB}_t(l)|\; \geq \beta_t$ (cf.~line~\ref{line3:u-update} and line~\ref{line3:l-update} of Algorithm~\ref{alg:bisection}). Per line~\ref{line3:begin-while} of the algorithm, $l$ and $u$ are updated until $u-l\leq \epsilon$. Then, after at most $\log_2[(u-l)/\epsilon]$ iterations, the algorithm terminates by setting $\lambda_t$ equal to the latest value of $l$ (lines~\ref{line3:lambda-update1}-\ref{line3:lambda-update2} of the algorithm).  Therefore, $\lambda_t$ is $\epsilon$-close to $\lambda_t^\star$ and satisfies $|\hat{\calB}_t(l)|\; \geq \beta_t$, which in turn implies $\hat{f}_t(\lambda_t)\leq f(\calA_{1:t}\setminus\calB^\star_{1:t})$, as desired.}
All in all, given an approximation $\lambda_t$ to $\lambda^\star_t$, Lemma~\ref{lem:bisection} implies the following approximation of Theorem~\ref{th:aposteriori}'s bounds.

\begin{corollary}[Pre-failure approximation of a posteriori bounds]\label{th:pre_failure_aposteriori}
Let Algorithm~\ref{alg:bisection} return $\lambda_t$, for $t=1,\ldots,T$. \alg selects $\calA_{t}$ such that  $|\calA_t|\;\leq \alpha_t$, and  if $f$ is submodular, then
		\begin{equation}\label{ineq:pre_failure_a_posteriori_bound_sub}
		\frac{f(\calA_{1:t}\setminus \mathcal{B}^\star_{1:t})}{f^\star_t}\geq
		\left\{\begin{array}{lr}
		\frac{1-{e^{-\kappa_f}}}{\kappa_f}
		\frac{\hat{f}_1(\lambda_1) }{f(\calM_{1})}, & t=1;\\
		\frac{1}{1+\kappa_f}
		\frac{\hat{f}_t(\lambda_t) }{f(\calM_{1:t})}, & t>1;
		\end{array} \right. 
		\end{equation}
whereas if $f$ is $c_f$-submodular, then
		\begin{equation}\label{ineq:pre_failure_a_posteriori_bound_non_sub}
		\frac{f(\calA_{1:t}\setminus \mathcal{B}^\star_{1:t})}{f^\star_t}\geq(1-c_f)\frac{\hat{f}_t(\lambda_t) }{f(\calM_{1:t})}.
		\end{equation}
\end{corollary}

Corollary~\ref{th:pre_failure_aposteriori} describes an online mechanism to predict \alg's performance before the upcoming failures, step by step.

{\begin{remark}[Utility of Corollary~\ref{th:pre_failure_aposteriori}'s bounds]\label{rem:utility}
For $T=1$, Corollary~\ref{th:pre_failure_aposteriori}'s bounds are computed before $\calB_1^\star$ occurs, which implies $\hat{f}_1$ can be computed \emph{offline} \emph{(}at any time $t<1$\emph{)}.   Therefore, Corollary~\ref{th:pre_failure_aposteriori}'s bounds, when tighter than Theorem~\ref{th:alg_rob_sub_max_performance}'s a priori bounds, allow for an a priori assessment of \alg's approximation performance \emph{(}\emph{before} \alg is deployed in the real world\emph{)}.  Examples of 1-step design problems where $T=1$, include \emph{robust actuator and sensor placement}~\cite{clark2016submodularity,summers2014submodularity,ye2020resilient}, \emph{robust feature selection}~\cite{krause2008robust,das2011spectral,bogunovic2017robust}, \emph{robust graph covering}~\cite{rahmattalabi2019exploring}, and \emph{robust server placement}~\cite{poularakis2017one,lu2020robust}. 

For $T>1$, Corollary~\ref{th:pre_failure_aposteriori}'s bounds  can only be computed \emph{online}, once $\calB_{1:t-1}^\star$ has been observed (but before $\calB_t^\star$ has occurred), for each  $t=1,\ldots,T$.\footnote{{Corollary~\ref{th:pre_failure_aposteriori}'s bounds  can only be computed \emph{online} since \alg itself is an online algorithm, computing $\calA_t$ only once $\calB_{1:t-1}^\star$ has been observed (yet before $\calB_t^\star$ has occurred), for each  $t=1,\ldots,T$.}}  As such, the bounds can be used to balance the trade-off between (i) computation time requirements (including computation and energy consumption requirements) for solving each step $t$ of \rsm, and (ii) approximation performance requirements for solving \rsm at each $t$. For example, if Corollary~\ref{th:pre_failure_aposteriori}'s bounds indicate a good performance by \alg  at $t$ (e.g., the bounds are above a given threshold), \alg is used to select $\calA_t$, since \alg is computation time inexpensive, being a polynomial time algorithm.  However, if the bounds indicate poor performance by \alg at $t$ (less than the given threshold), then an optimal algorithm can be used instead at $t$ (such as the one proposed in~\cite{rahmattalabi2018robust}), but at the expense of higher computation time, since any optimal algorithm is non-polynomial in the worst-case.\footnote{{Even when \alg is used in combination with another algorithm to choose $\calA_{1:t}$, Corollary~\ref{th:pre_failure_aposteriori}'s bounds are still applicable since they are \emph{algorithm agnostic} (cf.~proof of Theorem~\ref{th:aposteriori}).}}
\end{remark}}
%

\newcommand{\myhspace}{\hspace{-3.5mm}}
\newcommand{\mpw}{4.5cm}
\begin{figure*}[tbp]
	\begin{minipage}{\textwidth}
		\begin{tabular}{cccc}%
			\myhspace
			\begin{minipage}{\mpw}%
				\centering
				\includegraphics[width=1\columnwidth,trim=1cm 6.5cm 2.5cm 6.5cm, clip]{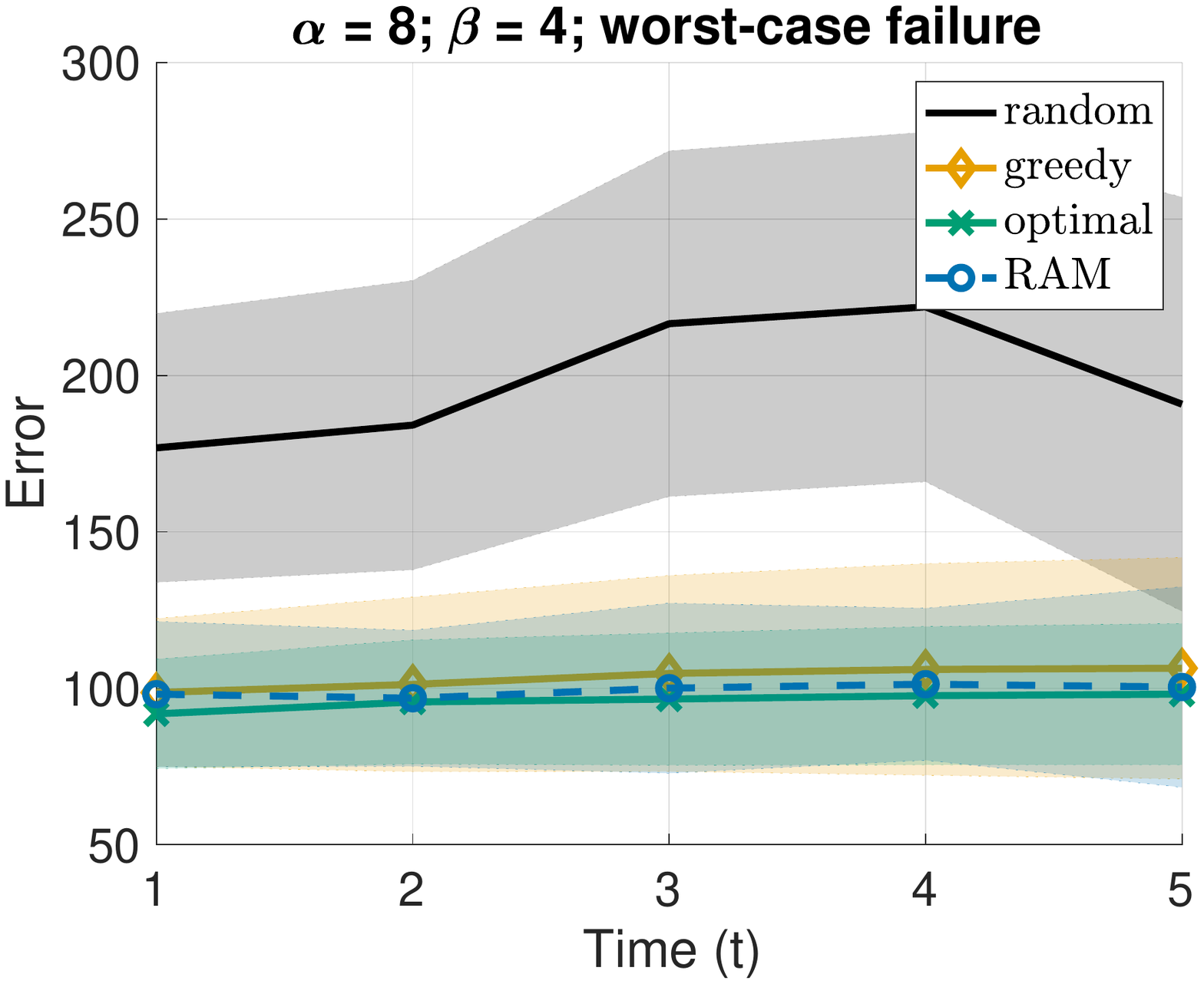}
			\end{minipage}
			&		\myhspace
			\begin{minipage}{\mpw}%
				\centering%
				\includegraphics[width=1\columnwidth,trim=1cm 6.5cm 2.5cm 6.5cm, clip]{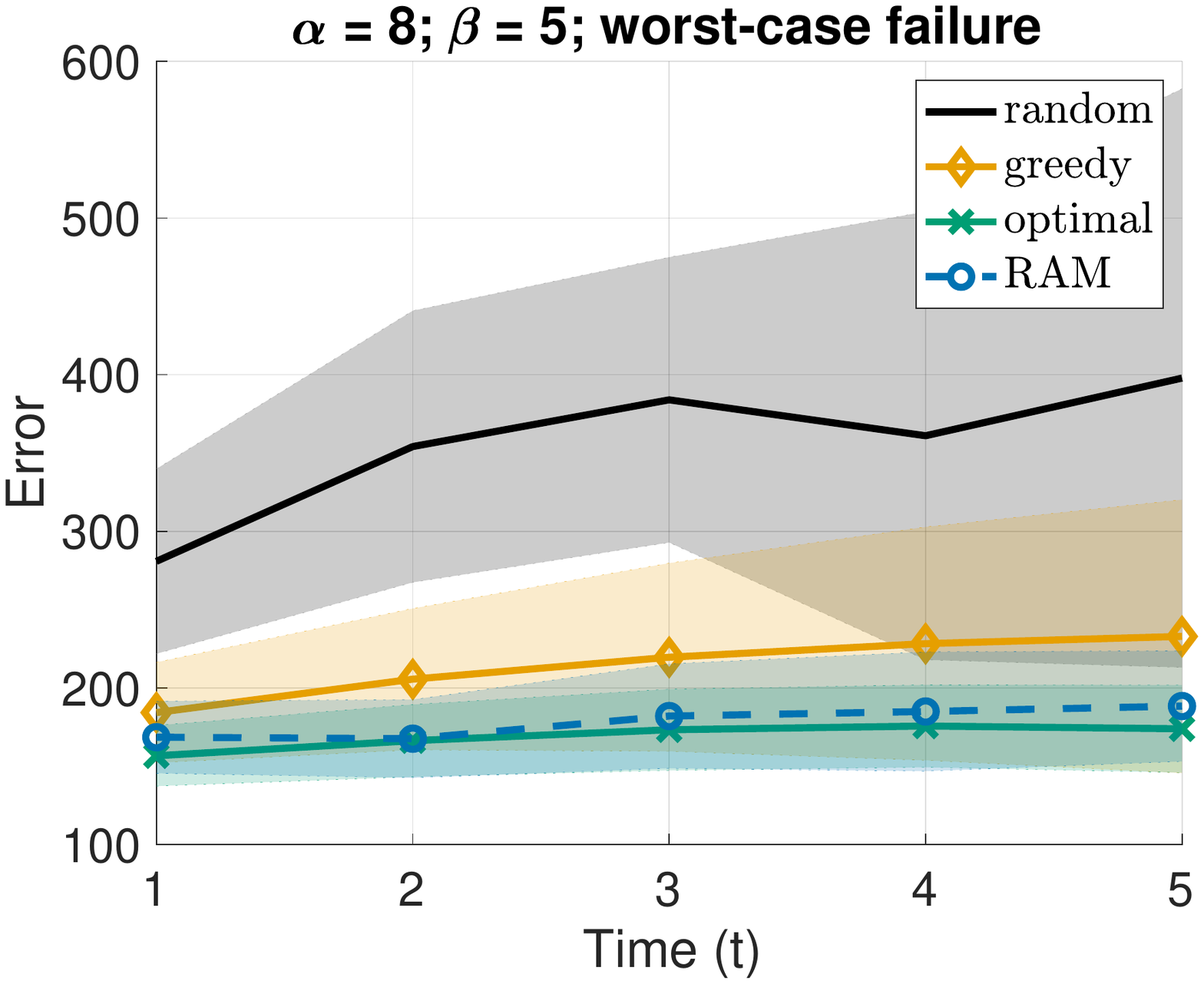}
			\end{minipage}
			&\myhspace
			\begin{minipage}{\mpw}%
				\centering
				\includegraphics[width=1\columnwidth,trim=1cm 6.5cm 2.5cm 6.5cm, clip]{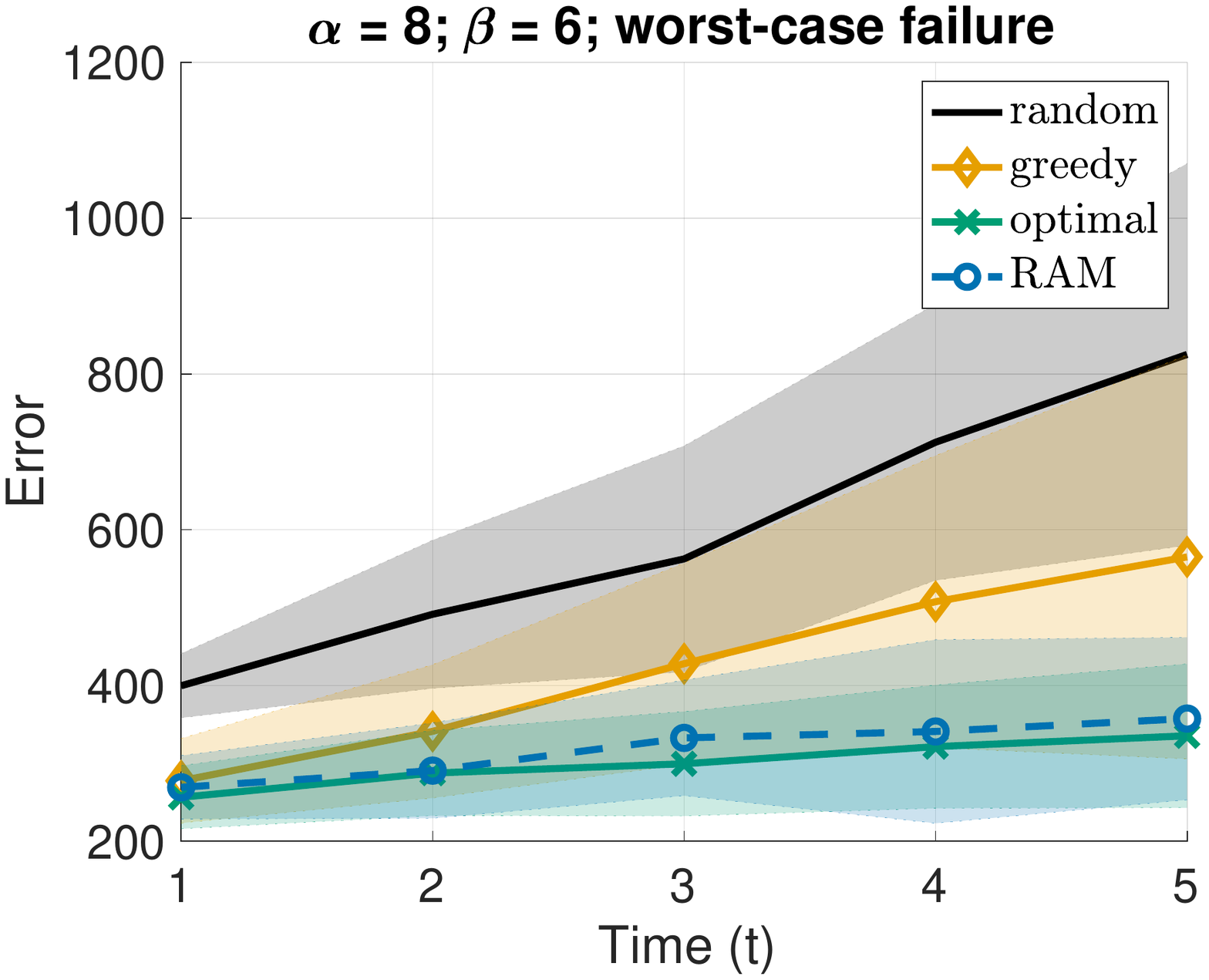}
			\end{minipage}
			& \myhspace
			\begin{minipage}{\mpw}%
				\centering%
				\includegraphics[width=1\columnwidth,trim=1cm 6.5cm 2.5cm 6.5cm, clip]{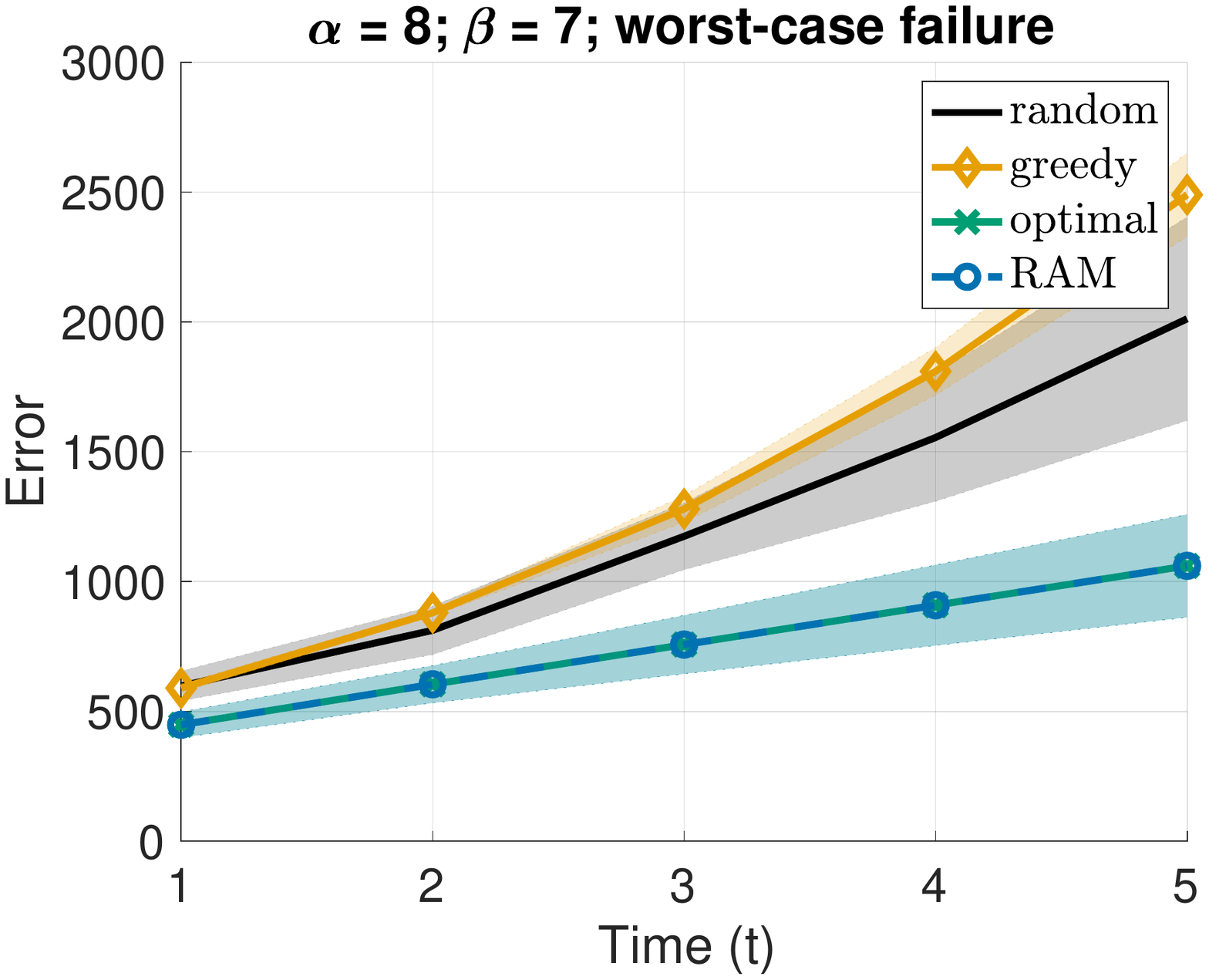} 
			\end{minipage}
			\\
			\myhspace
			\begin{minipage}{\mpw}%
				\centering
				\includegraphics[width=1\columnwidth,trim=1cm 6.5cm 2.5cm 6.5cm, clip]{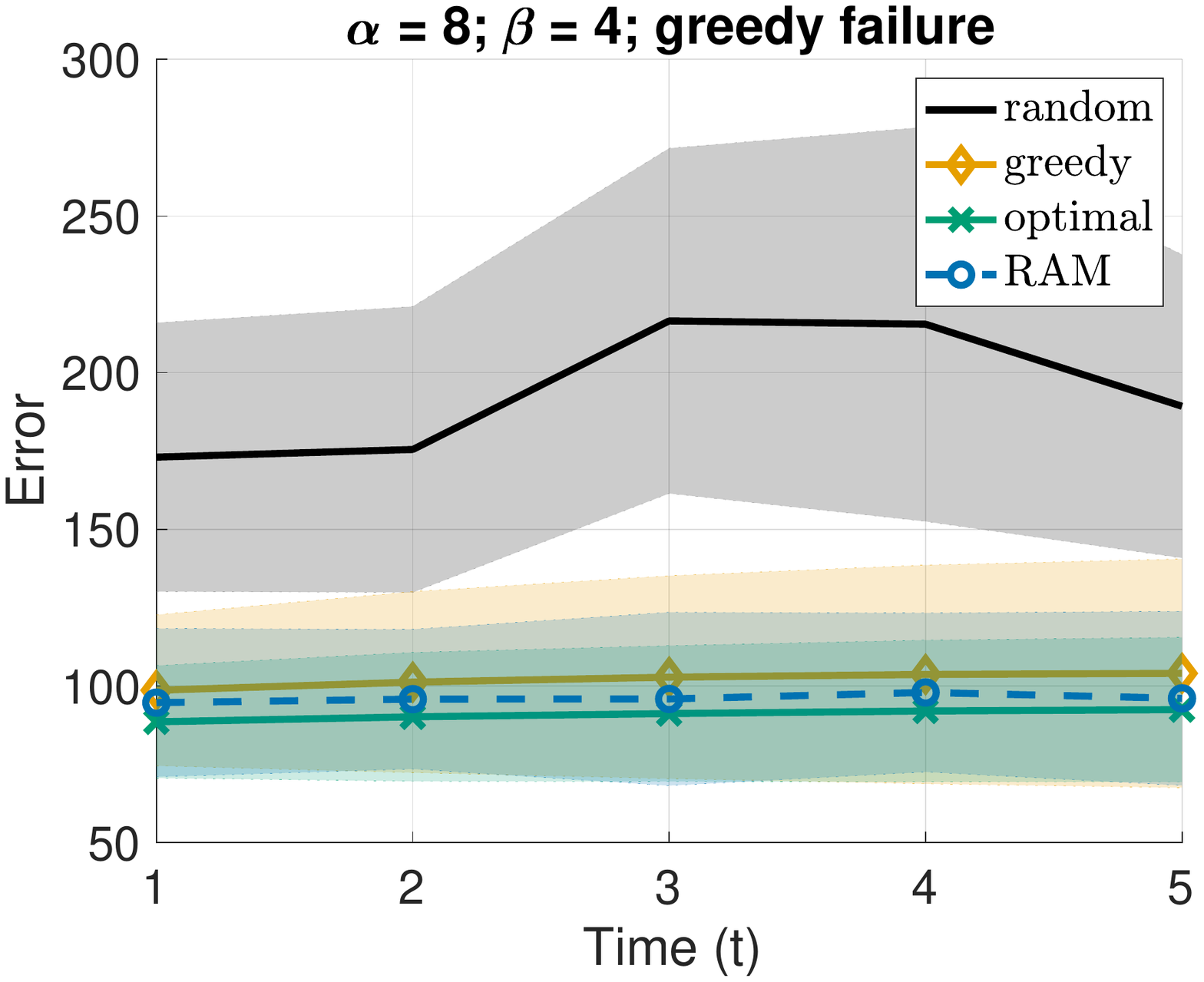}
			\end{minipage}
			&		\myhspace
			\begin{minipage}{\mpw}%
				\centering%
				\includegraphics[width=1\columnwidth,trim=1cm 6.5cm 2.5cm 6.3cm, clip]{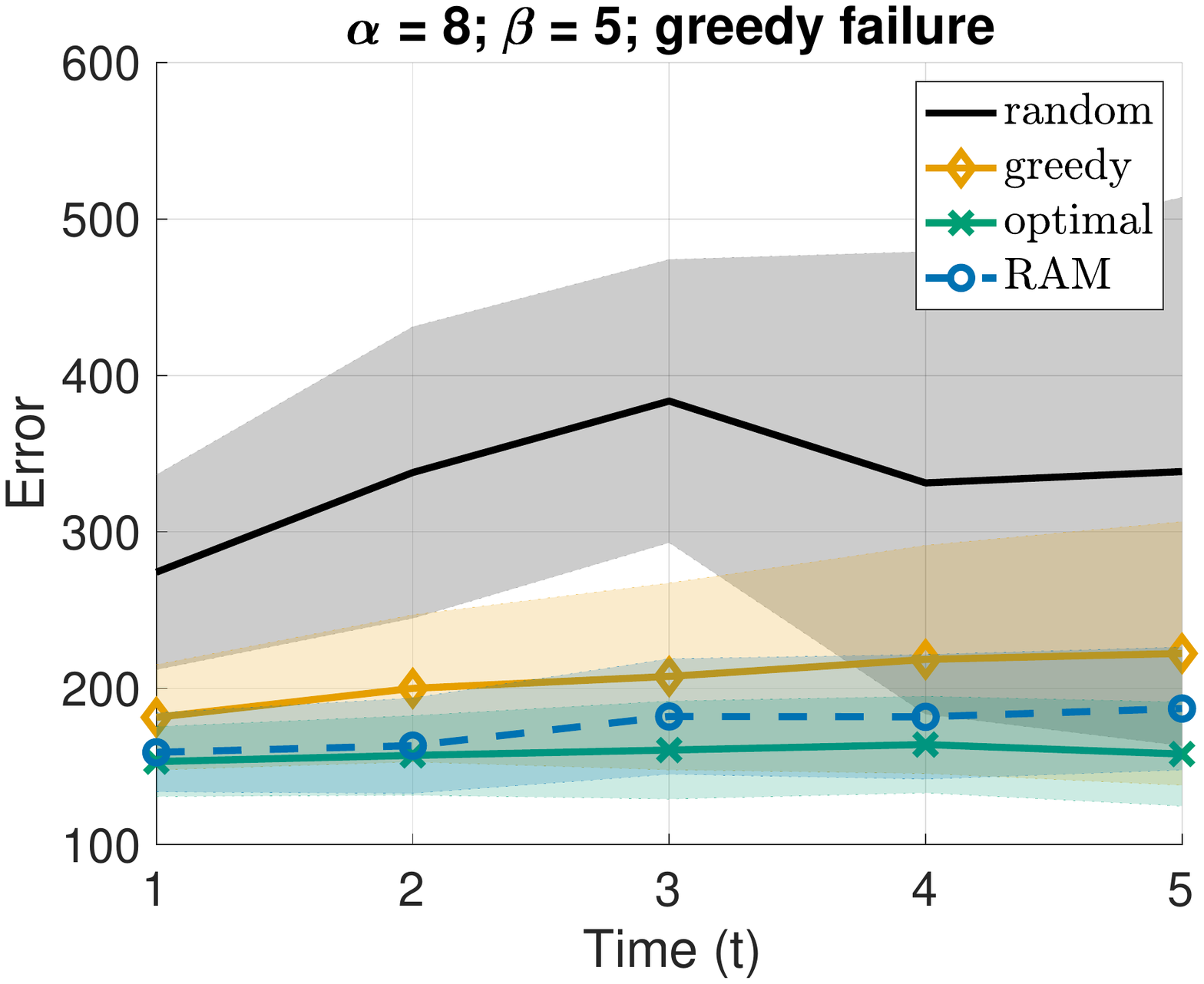}
			\end{minipage}
			&\myhspace
			\begin{minipage}{\mpw}%
				\centering
				\includegraphics[width=1\columnwidth,trim=1cm 6.5cm 2.5cm 6.3cm, clip]{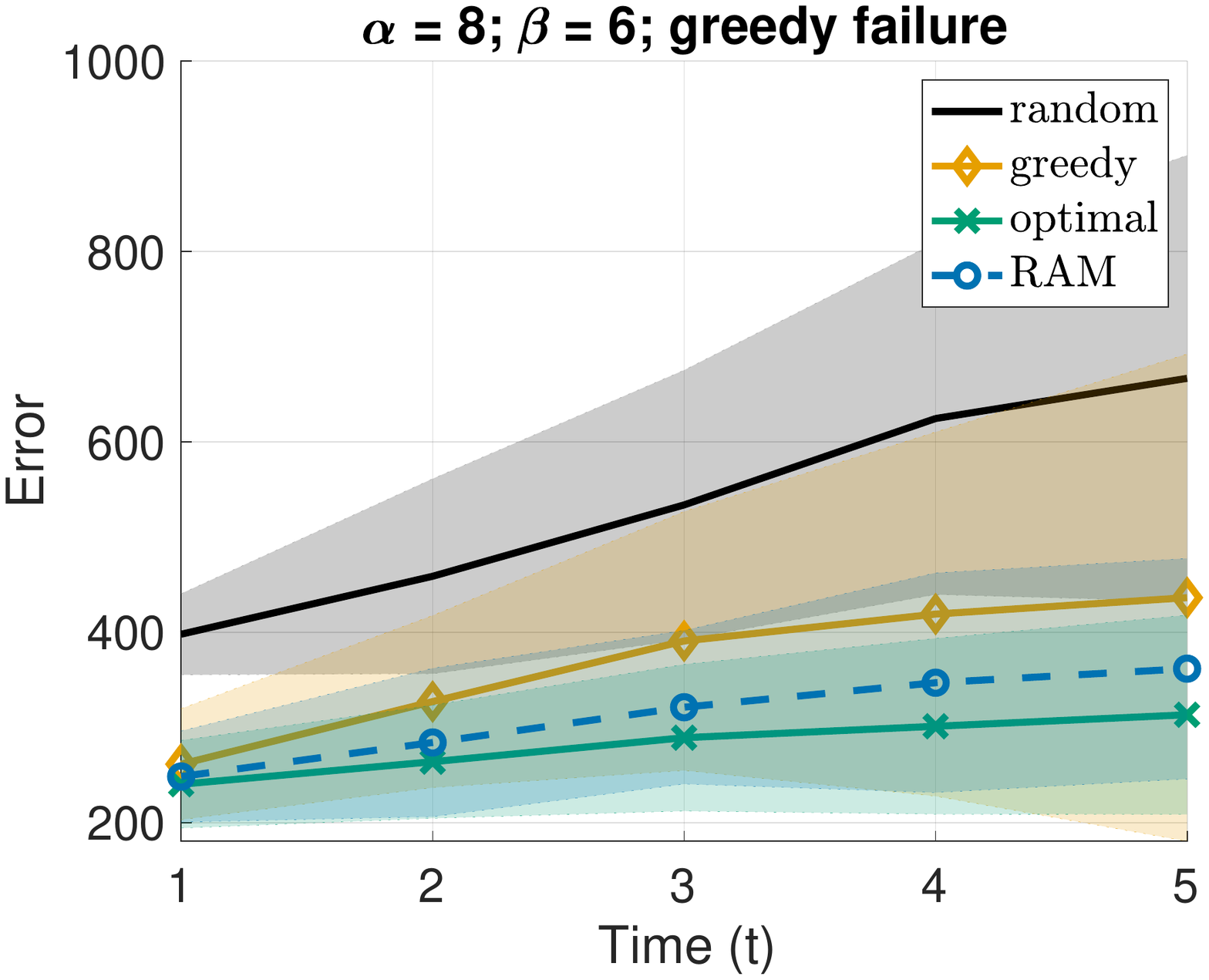}
			\end{minipage}
			& \myhspace
			\begin{minipage}{\mpw}%
				\centering%
				\includegraphics[width=1\columnwidth,trim=1cm 6.5cm 2.5cm 6.3cm, clip]{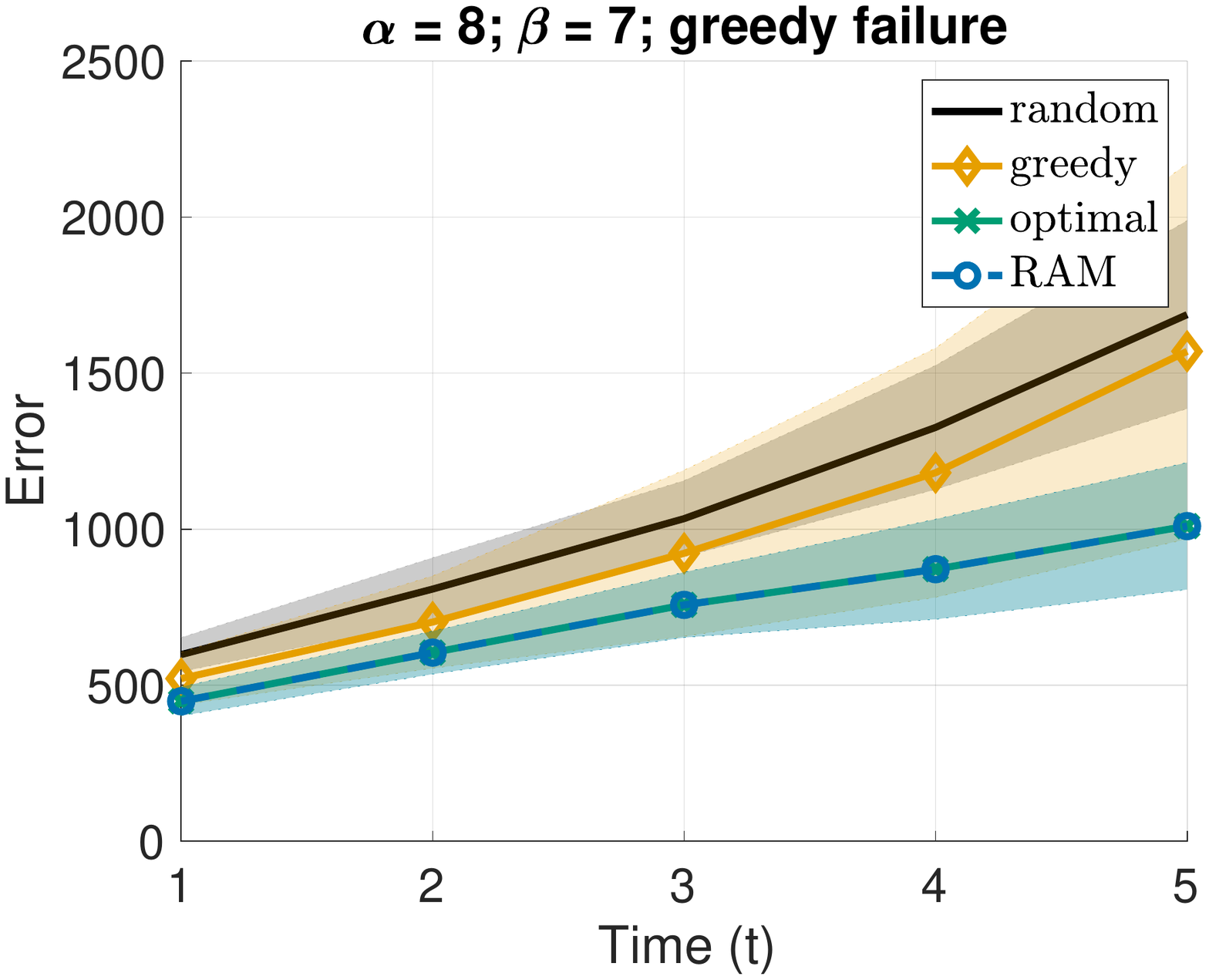} 
			\end{minipage}
			\\
			\myhspace
			\begin{minipage}{\mpw}%
				\centering
				\includegraphics[width=1\columnwidth,trim=1cm 6.5cm 2.5cm 6.3cm, clip]{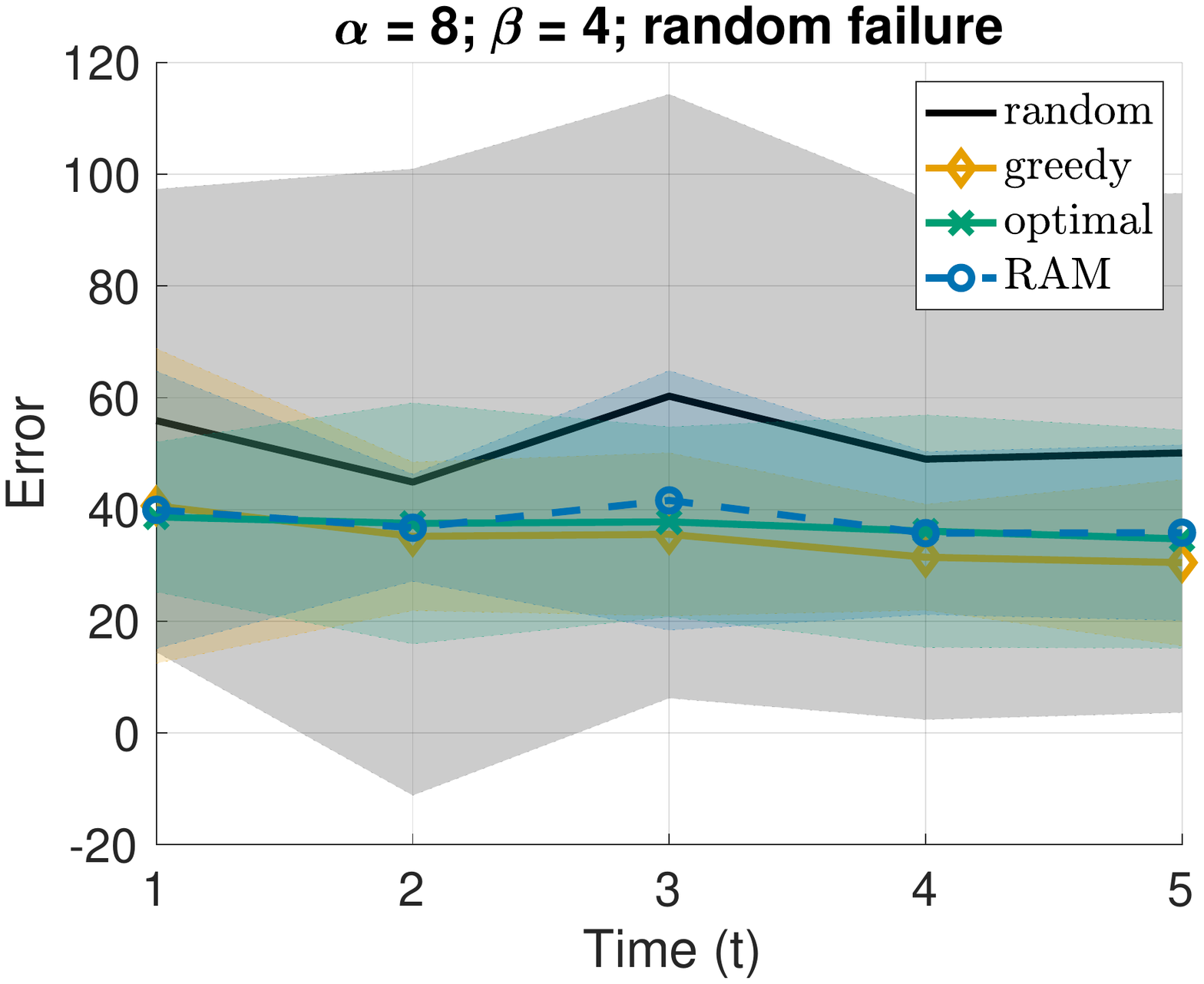}
			\end{minipage}
			&		\myhspace
			\begin{minipage}{\mpw}%
				\centering%
				\includegraphics[width=1\columnwidth,trim=1cm 6.5cm 2.5cm 6.3cm, clip]{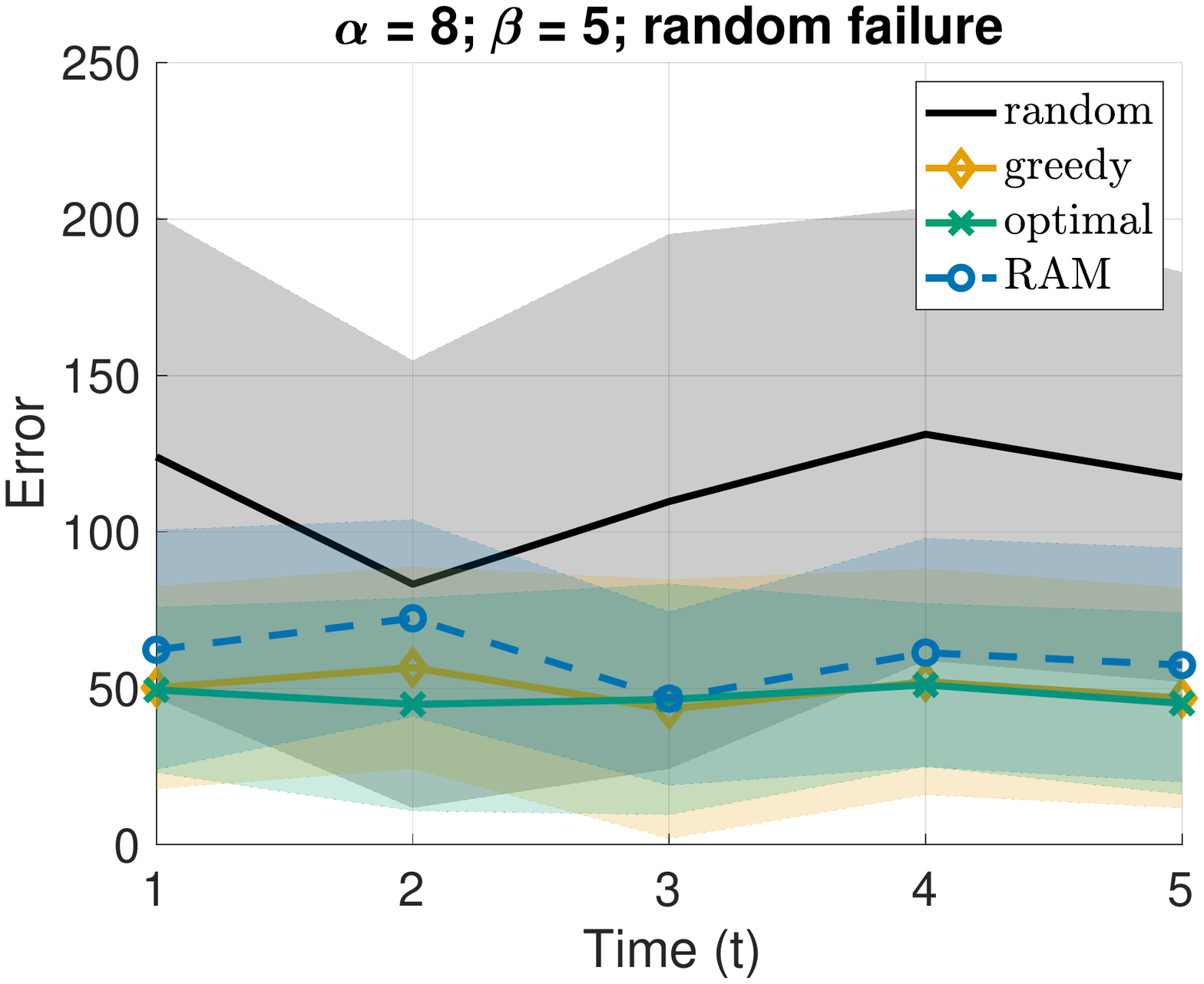}
			\end{minipage}
			&\myhspace
			\begin{minipage}{\mpw}%
				\centering
				\includegraphics[width=1\columnwidth,trim=1cm 6.5cm 2.5cm 6.5cm, clip]{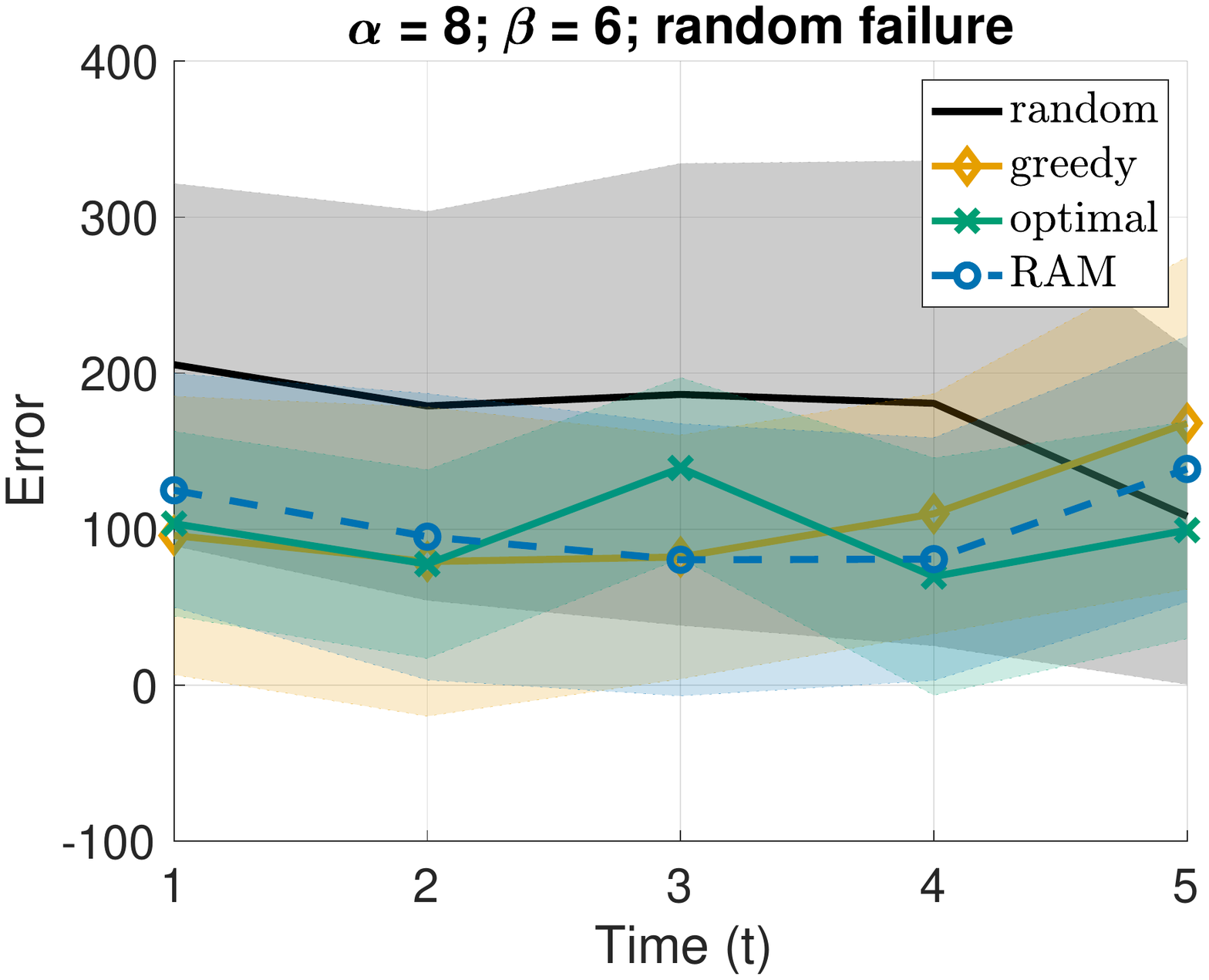}
			\end{minipage}
			& \myhspace
			\begin{minipage}{\mpw}%
				\centering%
				\includegraphics[width=1\columnwidth,trim=1cm 6.5cm 2.5cm 6.5cm, clip]{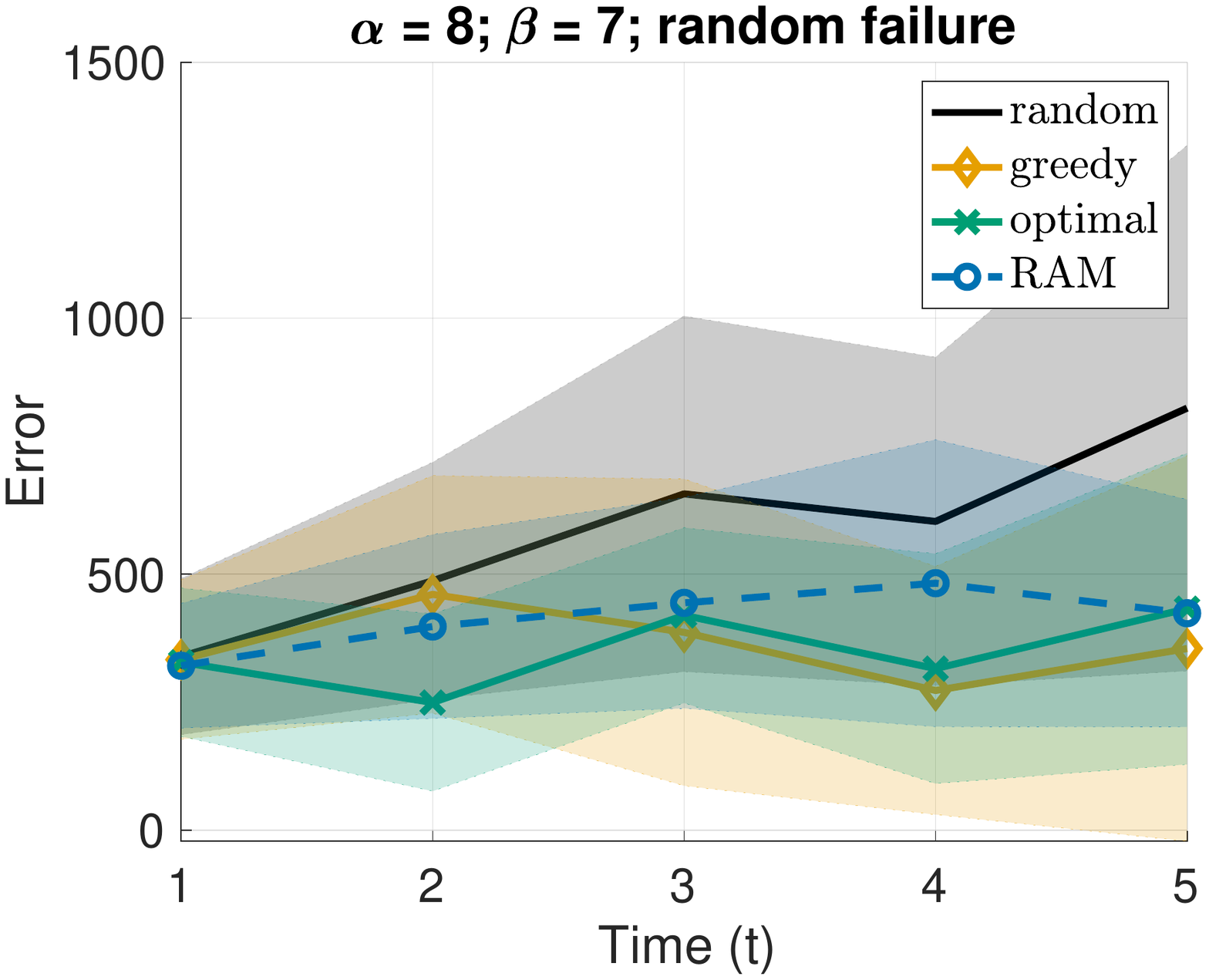} 
			\end{minipage}
		\end{tabular}
	\end{minipage}%
	\vspace{-1mm}
	\caption{\label{fig:flying}\small
		Representative simulation results for the application \textit{sensor scheduling for autonomous navigation}.  Results are averaged across 100 Monte Carlo runs.  Depicted is the estimation error for increasing time $t$, per eq.~\eqref{eq:opt_sensors}, where $\alpha_t=\alpha=8$ across all subfigures, whereas $\beta_t=\beta$ where $\beta$ varies across subfigures column-wise. Finally, the failure type also varies, but row-wise. \textit{Each subfigure has different scale.}
	}
\end{figure*}

\section{Applications}\label{sec:simulations}

We evaluate \alg's performance in applications.  We start by assessing its near-optimality against {worst-case} failures.  We continue by testing its sensitivity against \textit{non} worst-case failures, particularly, random and greedily selected failures. For such failures, one would expect \alg's performance to be the same, or improve, since \alg is designed to withstand the worst-case. 
To these ends, we consider two applications from the introduction: \textit{sensor scheduling for autonomous navigation}, and 
\textit{target tracking with wireless sensor networks}.

\subsection{Sensor scheduling for autonomous navigation}
\label{sec:autonomous}

We demonstrate \alg's performance in autonomous navigation scenarios, in the presence of sensing failures.  
We focus on small-scale instances, to enable \alg's comparison with a brute-force algorithm attaining the optimal to \rsm.  Instead, in Section~\ref{sec:tracking} we consider larger-scale instances.

A UAV moves in a 3D space, starting from a
randomly selected initial location. 
Its objective is to land at $[0,\;0,\;0]$ with zero velocity. 
The UAV is modeled as a double-integrator
with state $x_t = [p_t \; v_t]^\top \in \Real{6}$, where $t=1,2,\ldots$ is the time index, $p_t$ is the UAV's position, and $v_t$ is its velocity.  The UAV controls its acceleration.  The process noise has covariance $\eye_6$. 

The UAV is equipped with two on-board sensors: a GPS, measuring the 
UAV's position $p_t$ with a covariance $2 \cdot\eye_3$, 
and an altimeter, measuring $p_t$'s altitude component with standard deviation $0.5\rm{m}$. 
Also, the UAV 
can communicate with $10$ linear ground sensors.  These sensors are randomly generated at each Monte Carlo run, along with their noise covariance.

The UAV has limited on-board battery power and  measure- ment-processing bandwidth. Hence, it uses only a few sensors at each $t$.
Particularly, among the $12$ available sensors, the UAV uses at most $\alpha$, where $\alpha$  varies from 1 to 12 in the Monte Carlo analysis (per \rsm's notation, $\alpha_t=\alpha$ for all $t=1,2,\ldots$).
The UAV selects the sensors to minimize the cumulative batch-state error over a horizon $T=5$, captured~by
\begin{equation}\label{eq:opt_sensors}
c(\calA_{1:t})=\log\det[\Sigma_{1:t}(\calA_{1:t})],
\end{equation}
where $\Sigma_{1:t}(\calA_{1:t})$ is the error covariance of the minimum variance estimator of $(x_1,\ldots,x_t)$ given the used sensors up to $t$~\cite{anderson2015batch}.  Notably, $f(\calA_{1:t})=-c(\calA_{1:t})$ is non-decreasing and submodular, in congruence to \rsm's framework~\cite{tzoumas2016scheduling}.

Finally, we consider that at most $\beta$ failures are possible at each $t$ (per \rsm's notation, $\beta_t=\beta$ for all $t$). In the Monte Carlo analysis, $\beta$ varies from $0$ to $\alpha-1$.  

\myParagraph{Baseline algorithms}
We compare \alg with three algorithms.
The~first algorithm is a brute-force, optimal algorithm, denoted as {\scenario{optimal}}.  Evidently, {\scenario{optimal}}  is viable only for small-scale problem instances, such as herein where the available sensors are 12.
The second algorithm performs random selection and is denoted as {\scenario{random}}.
The third algorithm, denoted as {\scenario{greedy}}, greedily selects sensors to optimize eq.~\eqref{eq:opt_sensors} per the failure-free optimization setup in eq.~\eqref{eq:non_res}.

\myParagraph{Results} The results are averaged over 100 Monte Carlo runs.
For $\alpha=8$ and $\beta=4,5,6,7$, they are reported in Fig.~\ref{fig:flying}.  For the remaining $\alpha$ and $\beta$ values, the qualitative results are the same. From Fig.~\ref{fig:flying}, the following observations are due: 

\paragraph{Near-optimality against worst-case failures}  We focus on Fig.~\ref{fig:flying}'s first row of subfigures, where $\beta$ varies from 4 to 7 (from left to right).  Across all $\beta$, \alg nearly matches \scenario{optimal}. In contrast, \greedy nearly matches \optimal only for $\beta=4$ (and, generally, for $\beta\leq \alpha/2$, taking into account the simulation results for the remaining values of $\alpha$). Expectedly, \random is always the worst among all compared algorithms.  Importantly, as $\beta$ tends to $\alpha$, \greedy's performance tends to \random's. The observation exemplifies the insufficiency of the traditional optimization paradigm in eq.~\eqref{eq:non_res} against failures.

{Across all values of $\alpha$ and $\beta$ in the Monte Carlo analysis, the suboptimality bound in Theorem~\ref{th:aposteriori}'s eq.~\eqref{ineq:a_posteriori_bound_sub} is at least $.59$, informing \alg performs at least $50\%$ the optimal ($\kappa_f$ remains always less than $.93$, while $f(\calA_{1:t}\setminus \mathcal{B}^\star_{1:t})/f(\calM_{1:t})$ is close to $.95$).} In contrast, in Fig.~\ref{fig:flying} we observe an almost optimal performance. This is an example where the actual performance of the algorithm is significantly closer to the optimal than what is indicated by the algorithm's suboptimality bound. Indeed, this is a common observation for greedy-like algorithms: for the failure-agnostic \mbox{greedy in~\cite{fisher1978analysis} see, e.g.,~\cite{leskovec2007cost}.}

\paragraph{Robustness against non worst-case failures} We compare Fig.~\ref{fig:flying}'s subfigures column-wise, where the failure type varies among worst-case, greedy, and random (from top to bottom).\footnote{We refer to a failure $\calB_t$ as ``greedy,'' when $\calB_t$ is selected greedily towards minimizing $f(\calA_{1:t-1}\setminus\calB_{1:t-1},\calA_t\setminus\calB_t)$, where   $\calA_{1:t}$ and $\calB_{1:t-1}$ are given, as in Algorithm~\ref{alg:local} but now for minimization instead of maximization.}  Particularly, \alg's performance remains the same, or improves, against non worst-case~failures, and the best performance is being observed against random failures, as expected.  For example, if we focus on the rightmost column (where $\alpha=8; \beta=7$), at $t=5$, then we observe: for worst-case failures, \alg achieves error 1061; instead, for greedy failures, \alg achieves the reduced error 1010; while for random failures, \alg achieves even less error (less than 500).  Finally, against greedy failures, \alg is still superior to \greedy, while against random failures, they fare similarly.

Overall, the above numerical simulations demonstrate both the necessity for failure-robust optimization (\rsm), as well as the near-optimality of \alg, even for increasing number of failures (system-wide failures).  Similar conclusions we make over the second application scenario below.

\subsection{Target tracking with wireless sensor networks}\label{sec:tracking}

\begin{figure*}[tbp]
	\begin{minipage}{\textwidth}
		\begin{tabular}{cccc}%
			\myhspace
			\begin{minipage}{\mpw}%
				\centering
				\includegraphics[width=1\columnwidth,trim=1cm 6.5cm 2.5cm 6.5cm, clip]{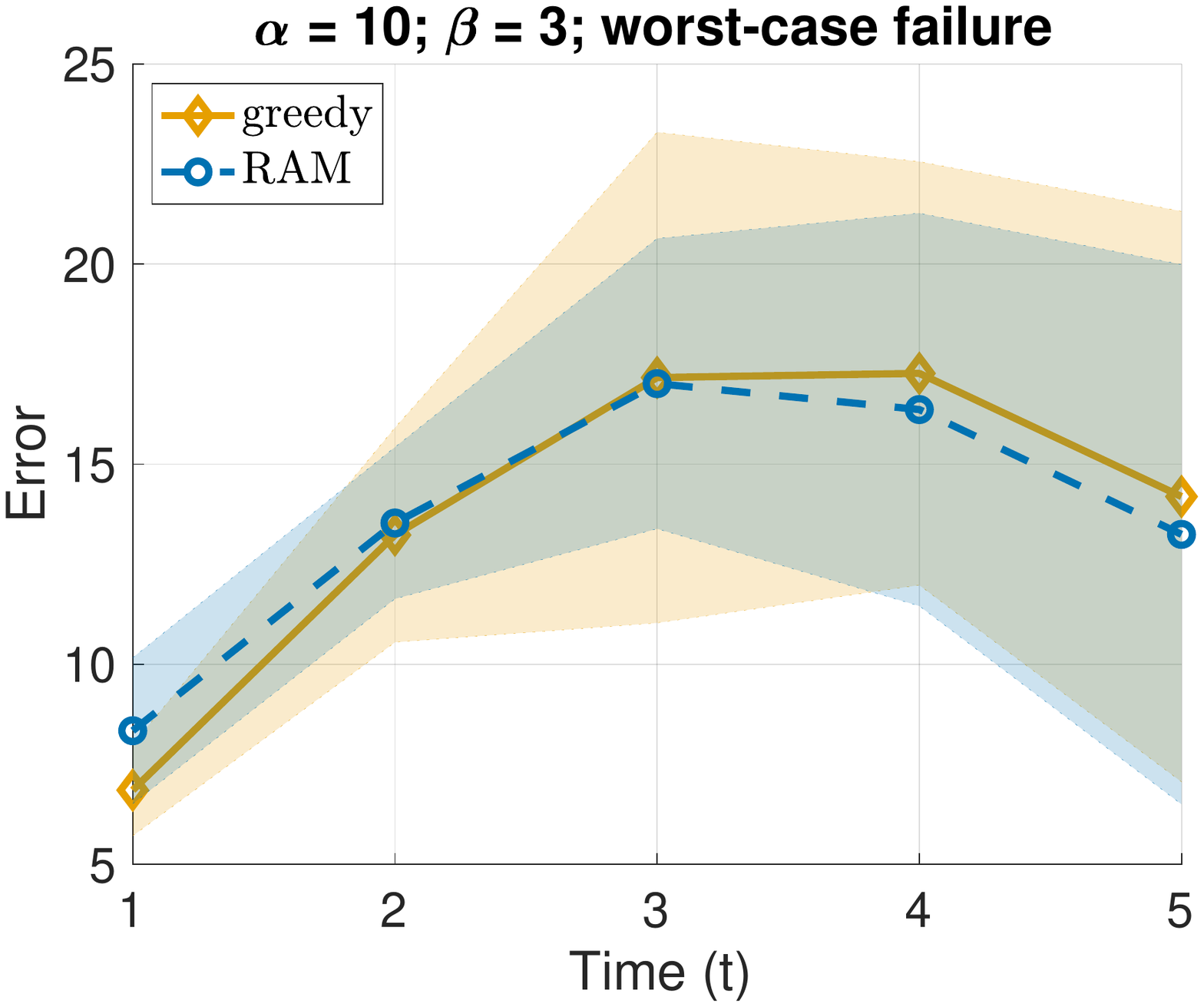}
			\end{minipage}
			&		\myhspace
			\begin{minipage}{\mpw}%
				\centering%
				\includegraphics[width=1\columnwidth,trim=1cm 6.5cm 2.5cm 6.5cm, clip]{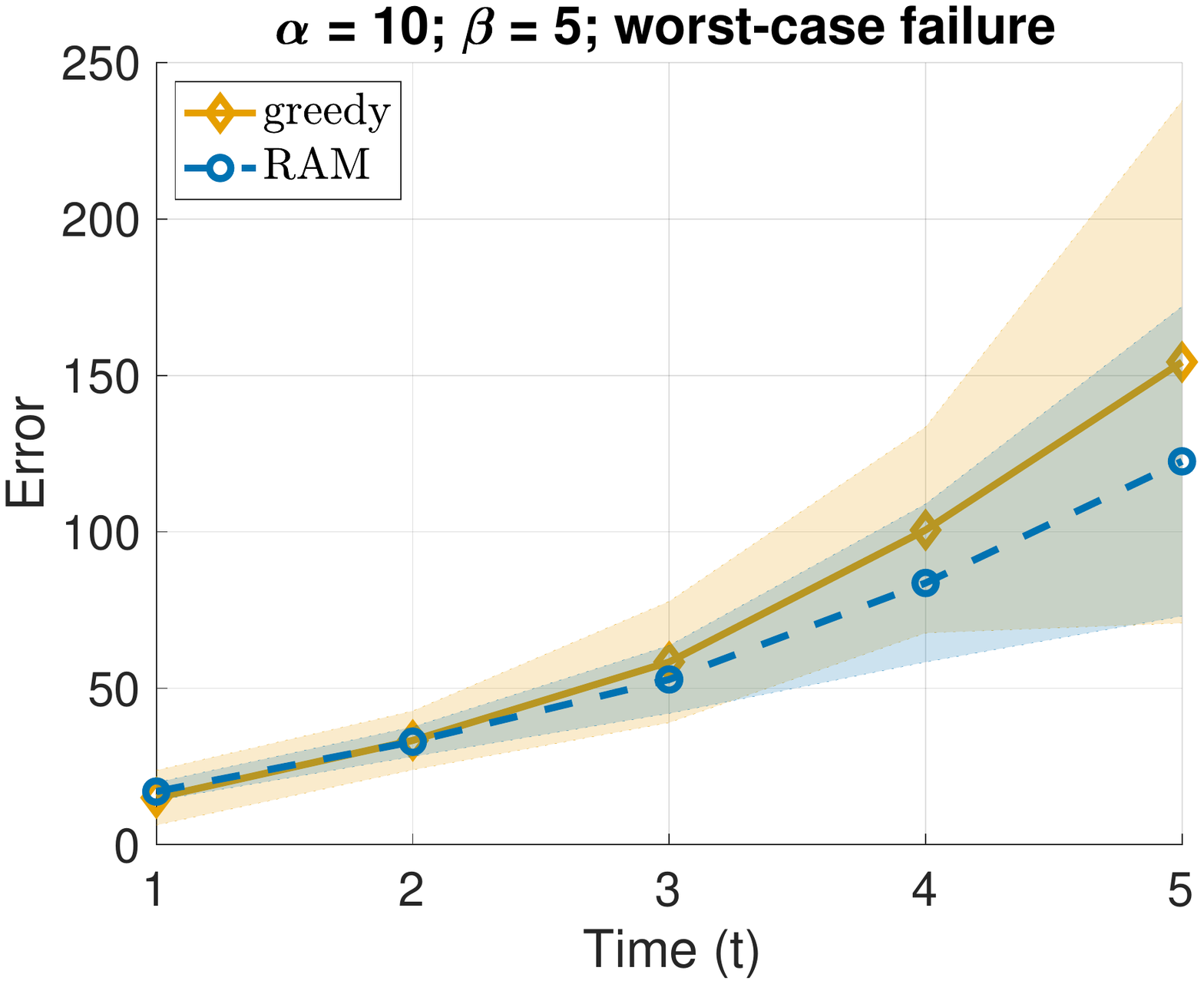}
			\end{minipage}
			&\myhspace
			\begin{minipage}{\mpw}%
				\centering
				\includegraphics[width=1\columnwidth,trim=1cm 6.5cm 2.5cm 6.5cm, clip]{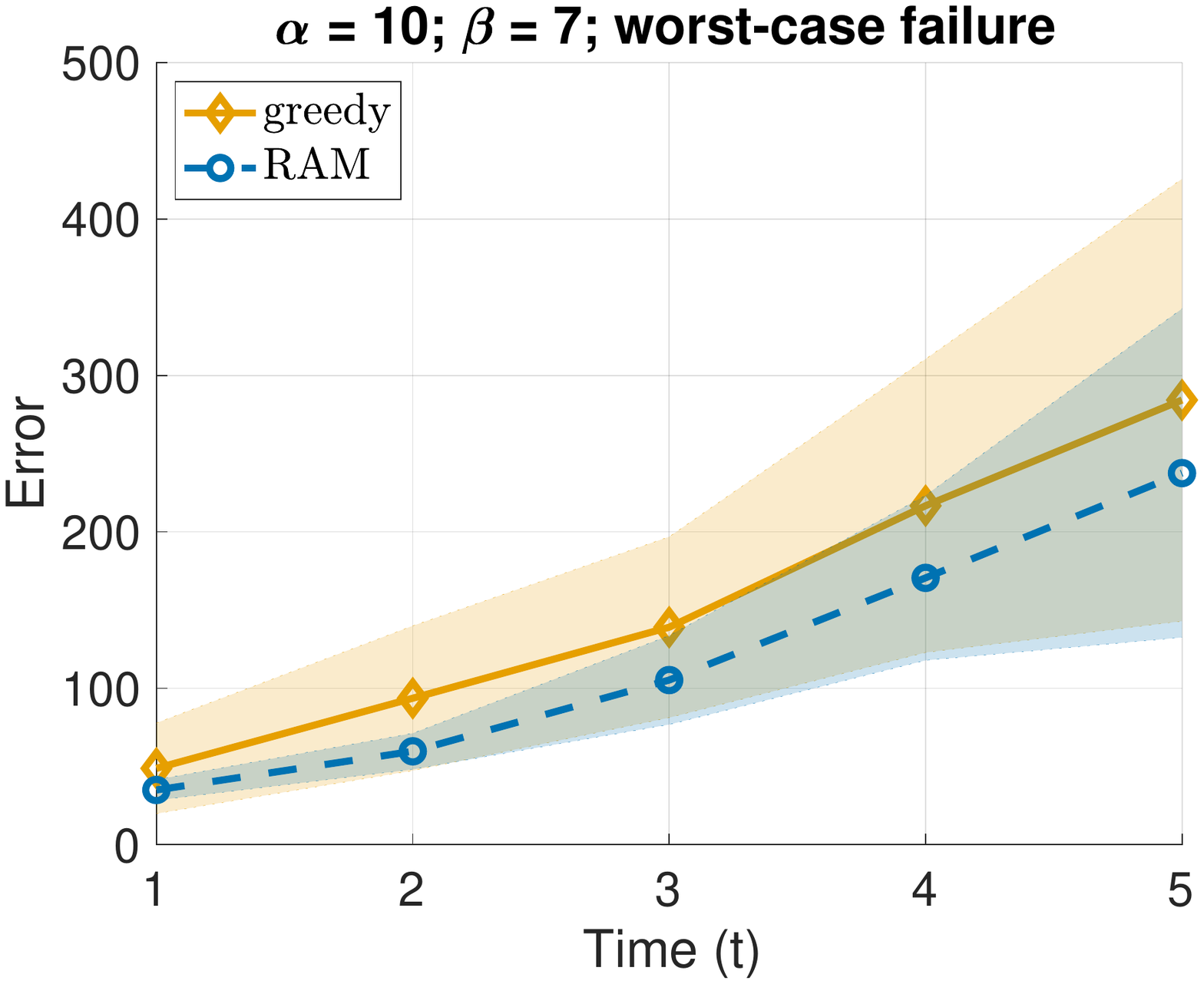}
			\end{minipage}
			& \myhspace
			\begin{minipage}{\mpw}%
				\centering%
				\includegraphics[width=1\columnwidth,trim=1cm 6.5cm 2.5cm 6.5cm, clip]{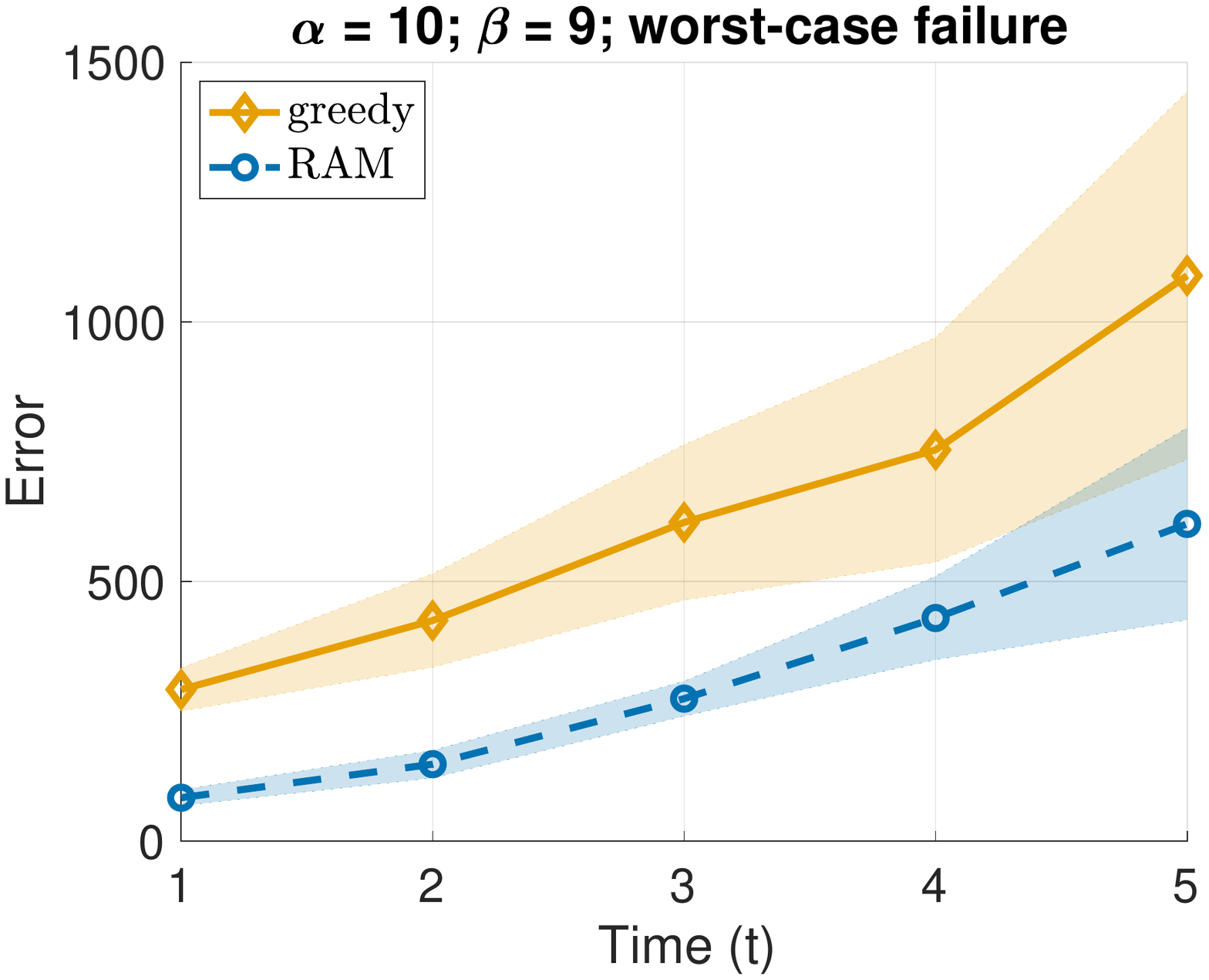} 
			\end{minipage}
			\\
			\myhspace
			\begin{minipage}{\mpw}%
				\centering
				\includegraphics[width=1\columnwidth,trim=1cm 6.5cm 2.5cm 6.5cm, clip]{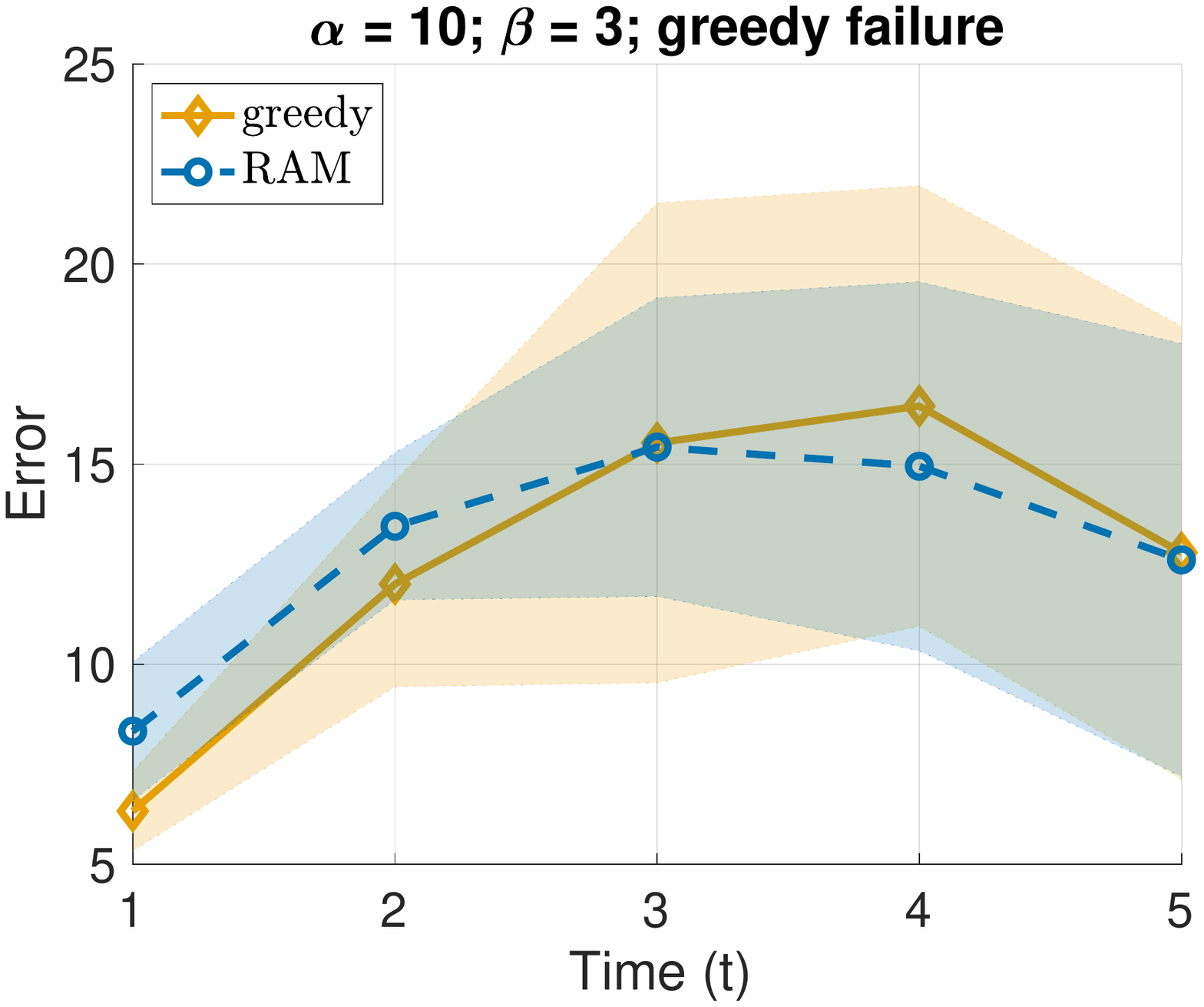}
			\end{minipage}
			&		\myhspace
			\begin{minipage}{\mpw}%
				\centering%
				\includegraphics[width=1\columnwidth,trim=1cm 6.5cm 2.5cm 6.3cm, clip]{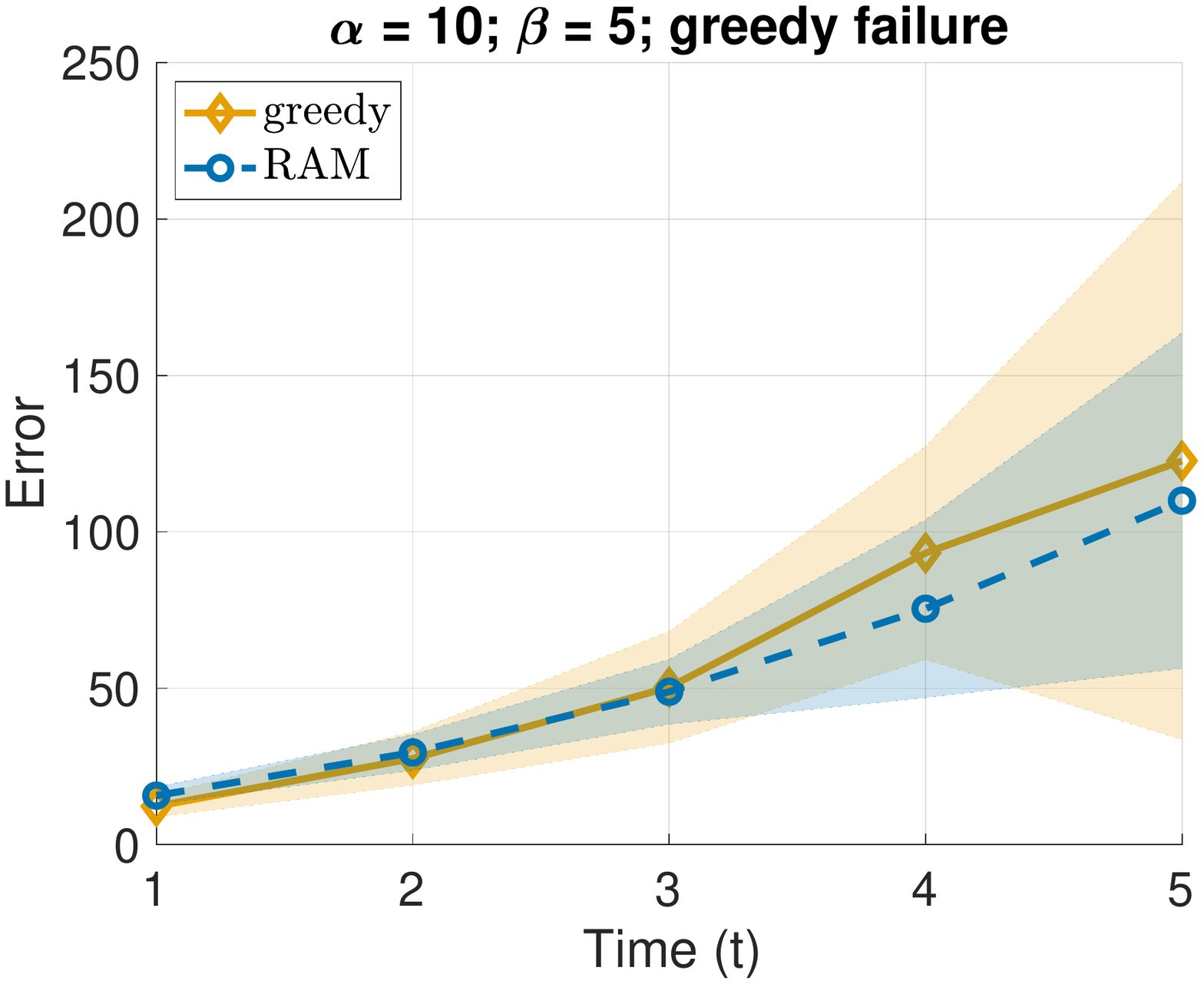}
			\end{minipage}
			&\myhspace
			\begin{minipage}{\mpw}%
				\centering
				\includegraphics[width=1\columnwidth,trim=1cm 6.5cm 2.5cm 6.3cm, clip]{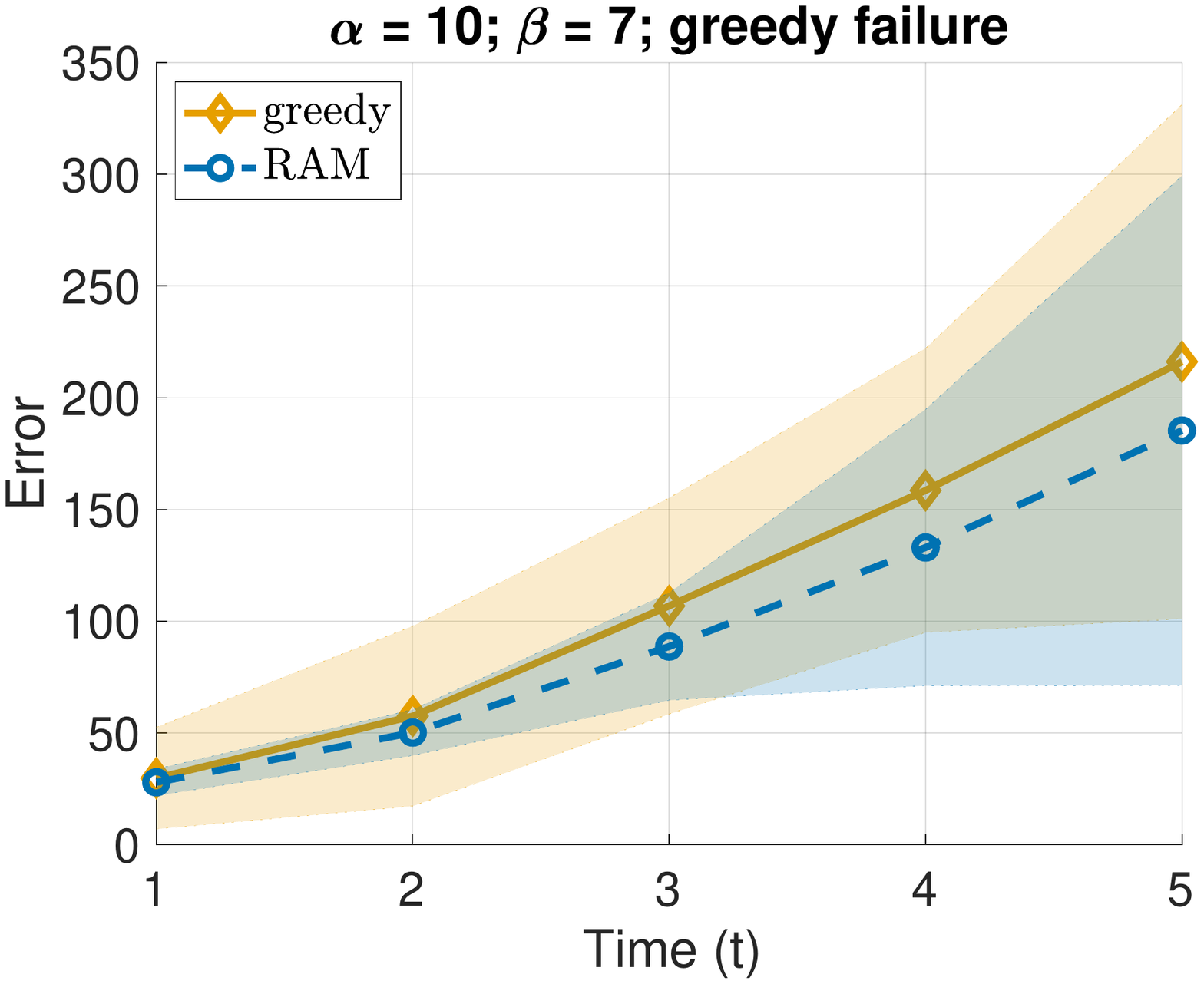}
			\end{minipage}
			& \myhspace
			\begin{minipage}{\mpw}%
				\centering%
				\includegraphics[width=1\columnwidth,trim=1cm 6.5cm 2.5cm 6.3cm, clip]{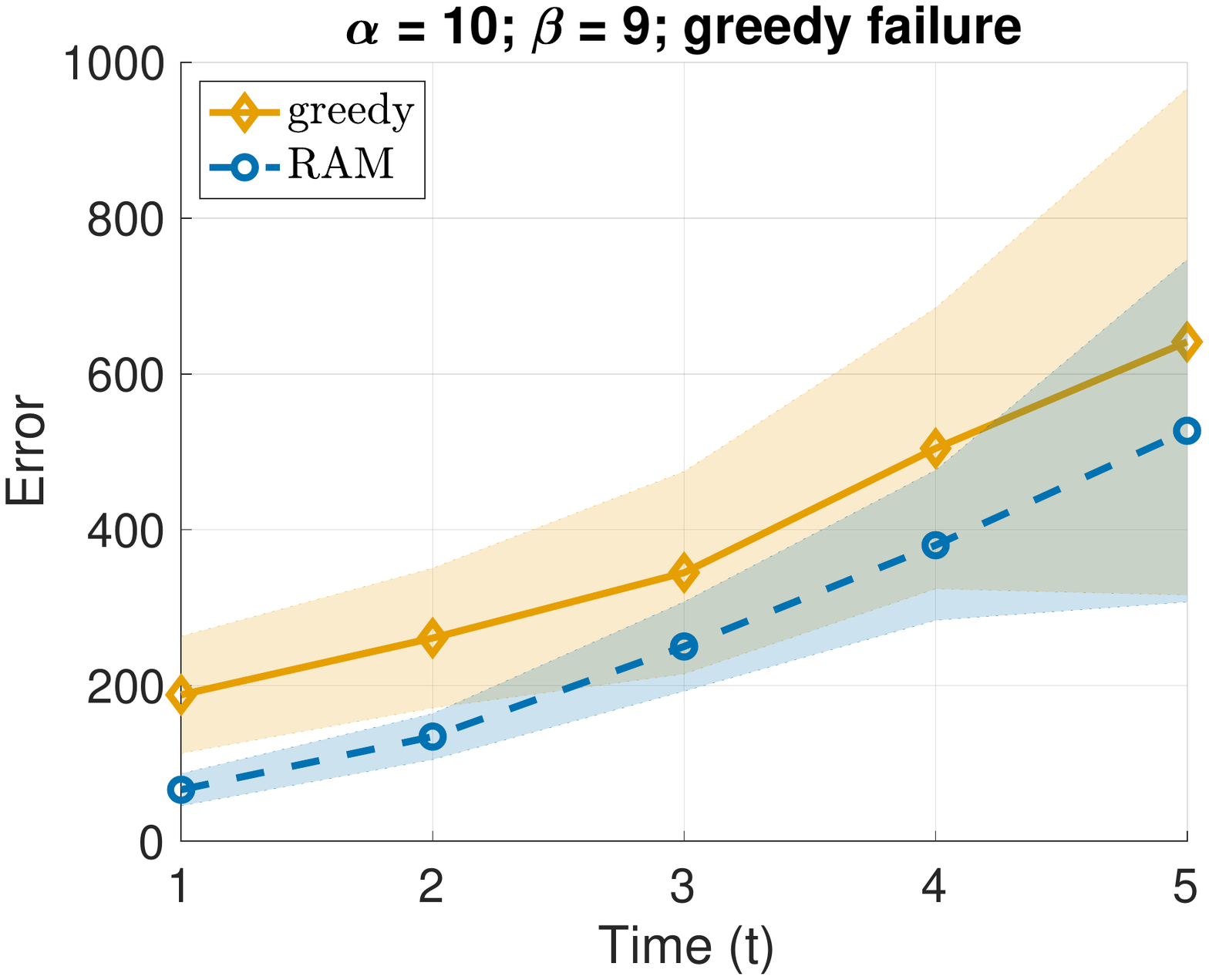} 
			\end{minipage}
			\\
			\myhspace
			\begin{minipage}{\mpw}%
				\centering
				\includegraphics[width=1\columnwidth,trim=1cm 6.5cm 2.5cm 6.3cm, clip]{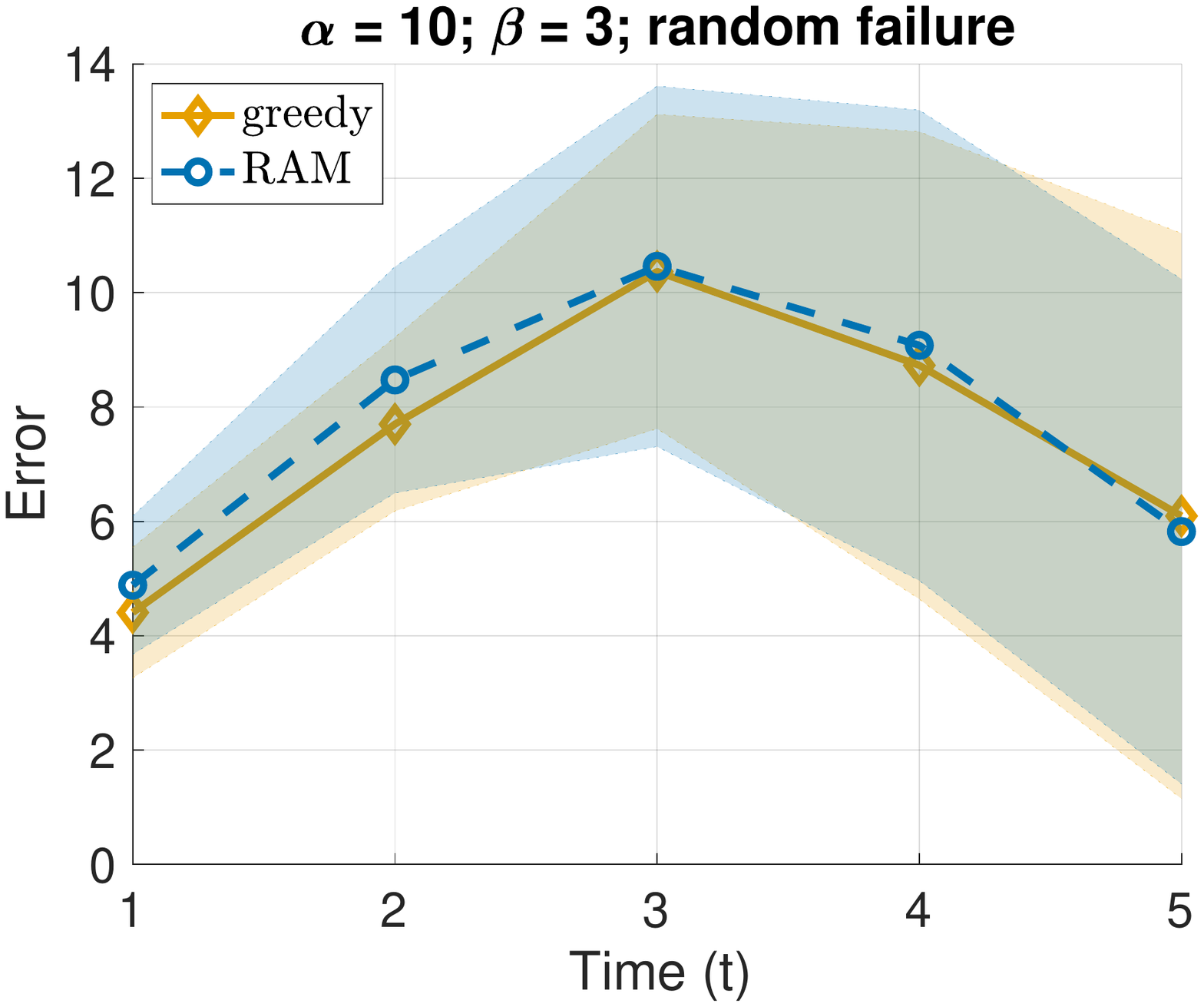}
			\end{minipage}
			&		\myhspace
			\begin{minipage}{\mpw}%
				\centering%
				\includegraphics[width=1\columnwidth,trim=1cm 6.5cm 2.5cm 6.3cm, clip]{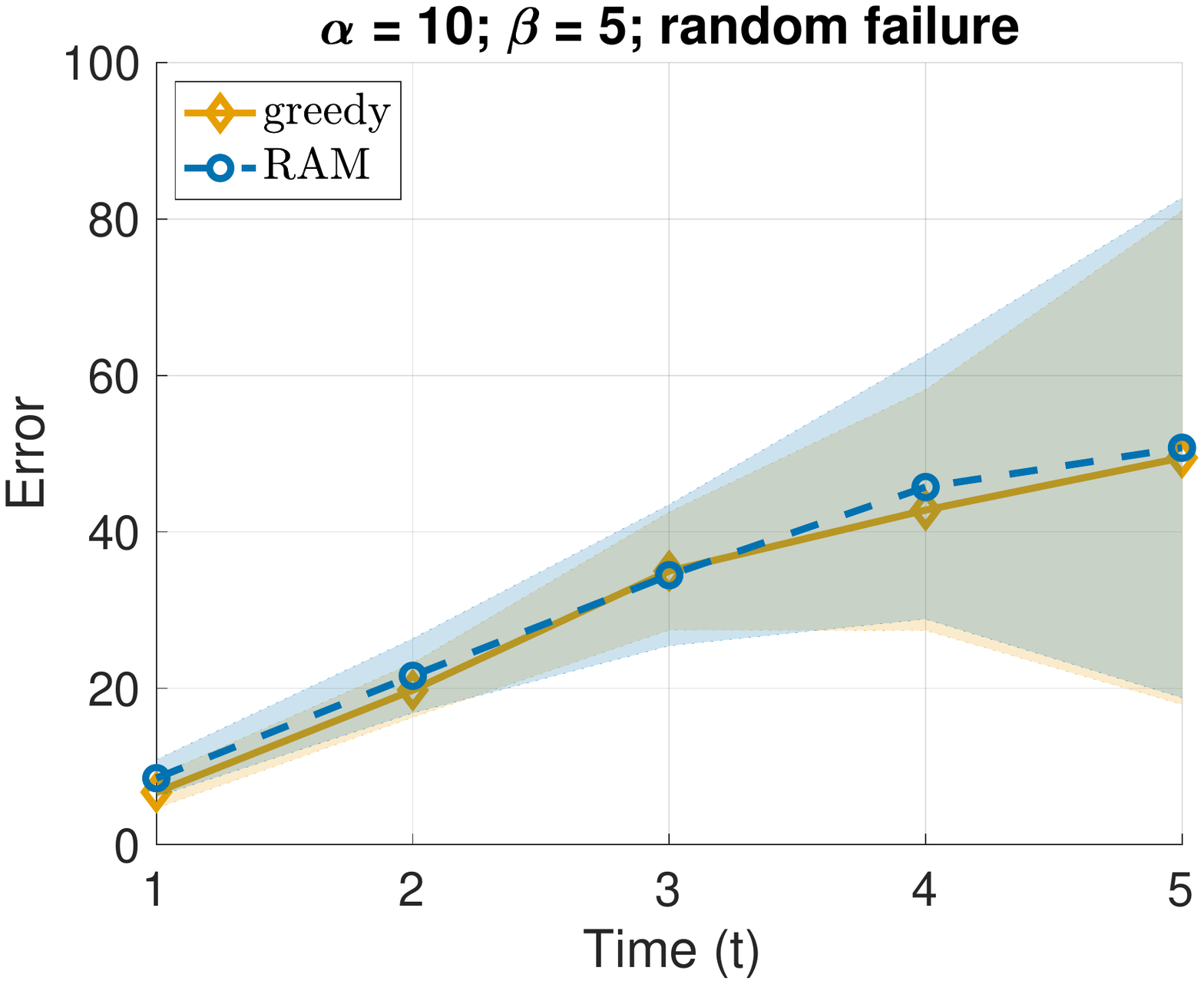}
			\end{minipage}
			&\myhspace
			\begin{minipage}{\mpw}%
				\centering
				\includegraphics[width=1\columnwidth,trim=1cm 6.5cm 2.5cm 6.5cm, clip]{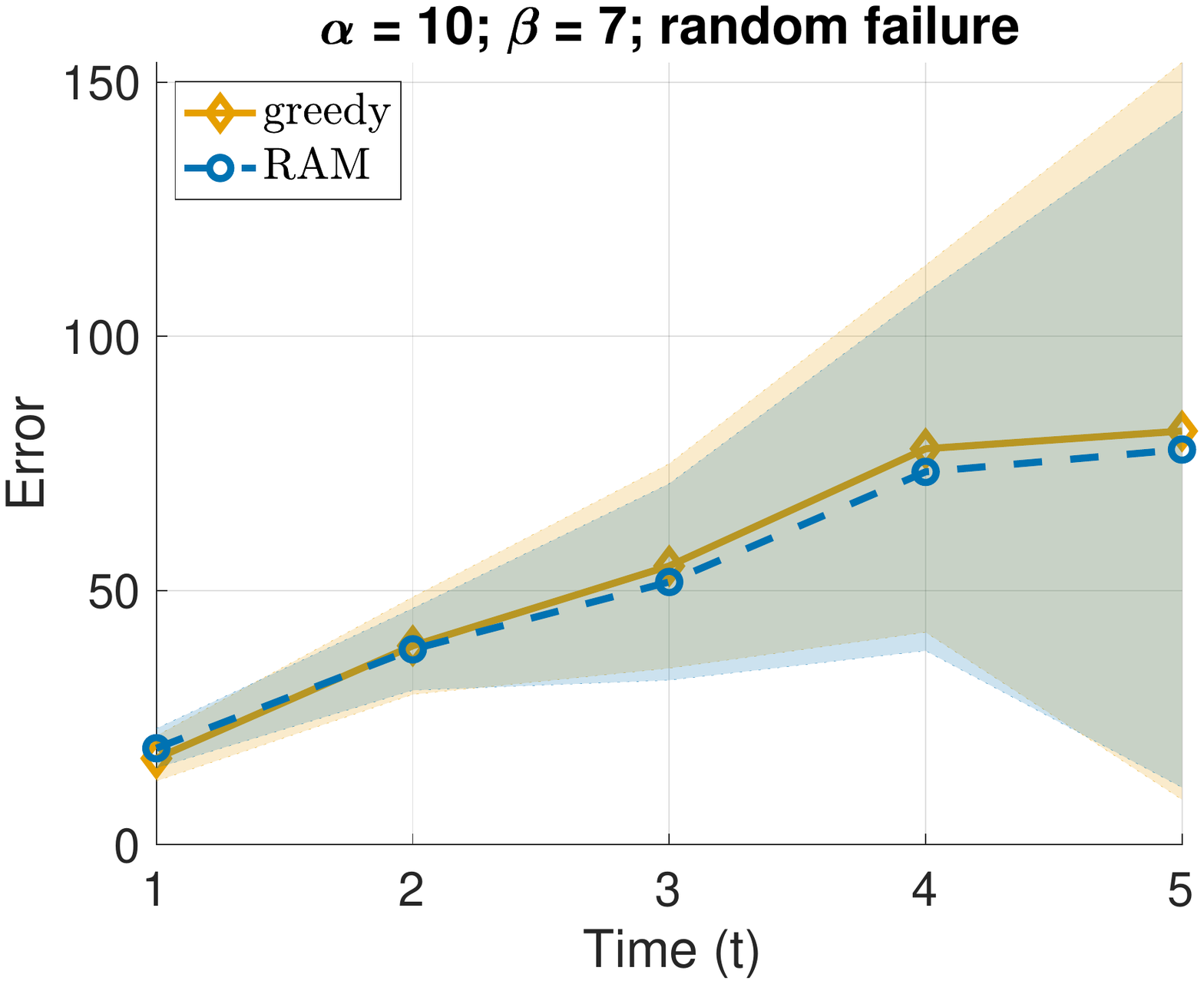}
			\end{minipage}
			& \myhspace
			\begin{minipage}{\mpw}%
				\centering%
				\includegraphics[width=1\columnwidth,trim=1cm 6.5cm 2.5cm 6.5cm, clip]{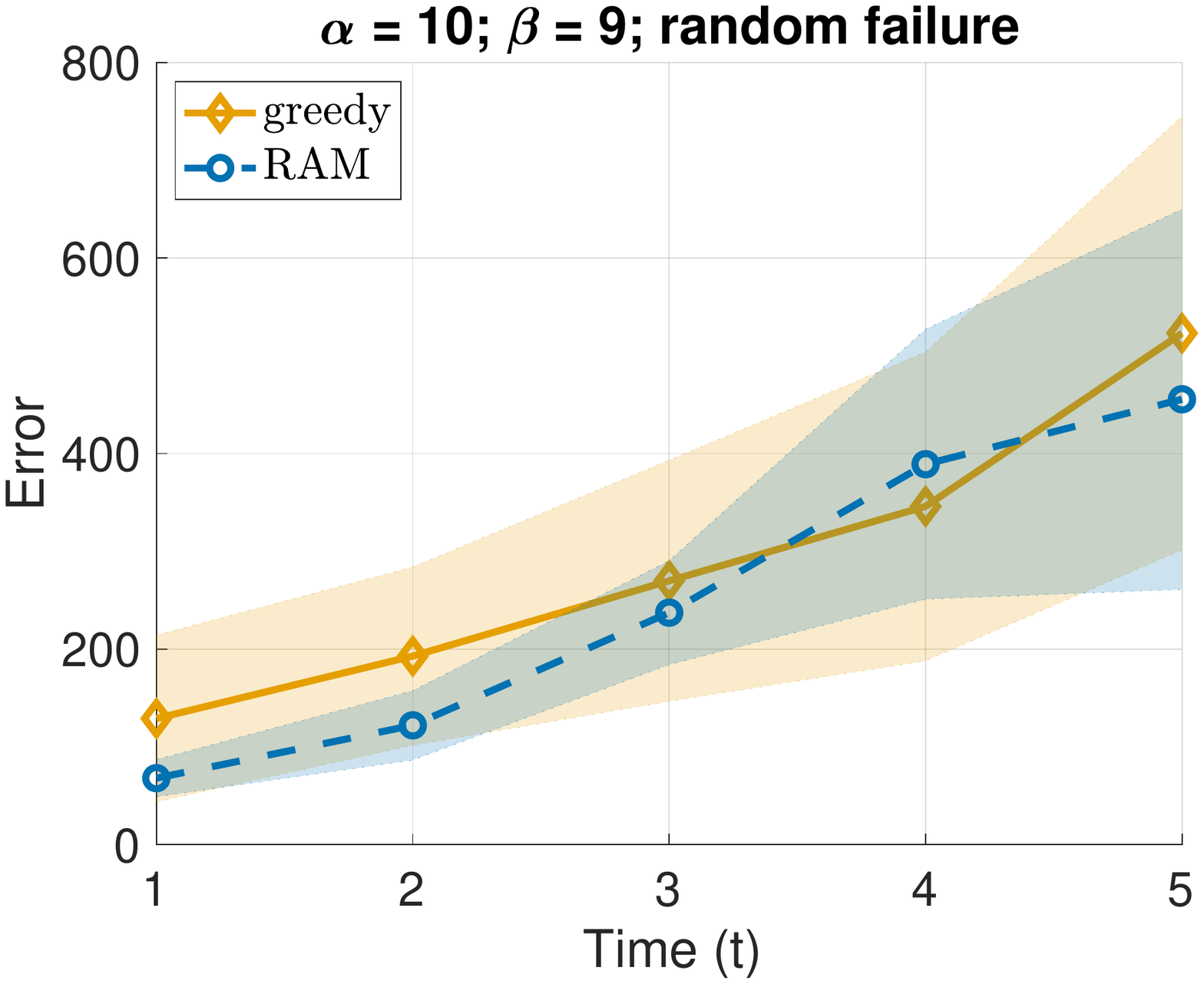} 
			\end{minipage}
		\end{tabular}
	\end{minipage}%
	\vspace{-1mm}
	\caption{\label{fig:sensitivity}\small
 Representative simulation results for the application \textit{target tracking with wireless sensor networks}.  Results are averaged over 100 Monte Carlo runs.  Depicted is the estimation error for increasing $t$, per eq.~\eqref{eq:trace_opt_sensors}, where $\alpha_t=\alpha=10$ across all subfigures, whereas $\beta_t=\beta$ where $\beta$ varies across subfigures column-wise. Finally, the failure type also varies, but row-wise. 
 \textit{Each subfigure has different scale.}
	}
\end{figure*}

We demonstrate \alg's performance in adversarial target tracking scenarios.  Particularly, we consider a mobile target who aims to escape detection from a wireless sensor network (WSN).  To this end, the agent causes failures to the network. 

A UAV (the target) is moving in a 3D, cubic shaped space.  The UAV moves on a straight line, across two opposite boundaries of the cube, keeping constant altitude and speed. The line's start and end points are randomly generated at each Monte Carlo run.  {The UAV's model is as in the autonomous navigation scenario in Section~\ref{sec:autonomous}.}

The WSN is composed of 100 ground sensors.  It is aware of the UAV's model, but can only noisily observe its state.   The sensors are randomly generated at each Monte Carlo run.

Due to power consumption and bandwidth limitations, only a few sensors can be active at each $t=1, 2, \ldots$. Particularly, we assume $\alpha=10$ active sensors at each $t$.  Also, we assume the sensors are activated so the cumulative Kalman filtering error over a horizon $T=5$ is minimized, as prescribed by
\begin{equation}\label{eq:trace_opt_sensors}
c(\calA_{1:t})=\sum_{t=1}^{T}\text{trace}[\Sigma_{t|t}(\calA_{1:t})],
\end{equation}
where $\Sigma_{t|t}(\calA_{1:t})$ is the Kalman filtering error covariance.  Noticeably, $f(\calA_{1:t})= -c(\calA_{1:t})$ is non-decreasing and $c_f$-submodular, in agreement with \rsm's framework~\cite{chamon2017mean}.

Finally, at most $\beta$ failures are possible at each $t$. In the Monte Carlo analysis, $\beta$ varies from $1$ to $\alpha-1=9$.  

\myParagraph{Baseline algorithms}
We compare \alg with {\scenario{random}}, and {\scenario{greedy}}.  We cannot compare with 
{\scenario{optimal}}, since the network's large-scale size makes {\scenario{optimal}} unfeasible.

\myParagraph{Results} The simulation results are averaged over 100 Monte Carlo runs. For $\beta=3,5,7,9$, they are reported in Fig.~\ref{fig:sensitivity}, where \random is excluded since it results to exceedingly larger errors.  For the remaining $\beta$ values, the qualitative results remain the same.  From Fig.~\ref{fig:sensitivity}, we make the observations:

\paragraph{Superiority against worst-case failures} We focus on Fig.~\ref{fig:sensitivity}'s first row, where $\beta$ takes the values 3, 5, 7, and 9 (from left to right).  For $\beta=3$ (also, for $\beta=1,2$, accounting for the remaining, non depicted simulations), \alg fares similar to \greedy.  In contrast, for the remaining values of $\beta$, \alg dominates \scenario{greedy}, achieving significantly lower error (observe the different scales among the subfigures for $\beta=5,7,9$).

{Across all $\beta$ values in the Monte Carlo analysis (including those in Fig.~\ref{fig:sensitivity}), the suboptimality bound in eq.~\eqref{ineq:a_posteriori_bound_non_sub} ranges from $.02$ to $.10$, informing that \alg performs at least $2\%$ to $10\%$ the optimal. Specifically, $c_f$ ranges from $.89$ to $.98$, whereas  $f(\calA_{1:t}\setminus \mathcal{B}^\star_{1:t})/f(\calM_{1:t})$ remains again close to $.95$.}   Hence, the possible conservativeness of the bound stems from the conservativeness of its term~$1-c_f$.

\paragraph{Robustness against non worst-case failures} We compare Fig.~\ref{fig:sensitivity}'s subfigures column-wise.  Similarly to the autonomous navigation scenarios, \alg's performance remains the same, or improves, against non worst-case~failures, and the lowest error is being observed against random failures.  For example, if we focus on the rightmost column (where $\alpha=10; \beta=9$), at $t=5$, then: for worst-case failures, \alg achieves error 611; in contrast, for greedy failures, \alg achieves the lower error 526; and for random failures, \alg achieves the even lower error 456.  Generally, against greedy failures, \alg is again still superior to \greedy; while against random failures, both have similar performance.

In summary, \alg remains superior even against system-wide failures, and even if the failures are non worst-case.

\section{Conclusion} \label{sec:con}

We made the first step to adaptively protect critical control, estimation, and machine learning applications against sequential failures.  Particularly, we focused on scenarios requiring the robust discrete optimization of systems per \rsm.  We provided \alg, the first online algorithm, which adapts to the history of failures, and guarantees a near-optimal performance even against system-wide failures despite its minimal running time. To quantify \alg's performance, we provided  per-instance a priori bounds and tight, optimal a posteriori bounds.  To this end, we used curvature notions, and contributed a first step towards characterizing the curvature's effect on the per-instance approximability of \rsm.  Our curvature-dependent bounds complement the current knowledge on the curvature's effect on the approximability of the failure-free optimization paradigm in eq.~\eqref{eq:non_res}~\cite{conforti1984curvature,lehmann2006combinatorial,sviridenko2017optimal,sviridenko2013optimal}. Finally, we supported our theoretical results with numerical evaluations.\\
\indent The paper opens several avenues for future research. 
One is the decentralized implementation of \alg towards robust multi-agent autonomy and large-scale network design.
And another is the extension of our results to optimization frameworks with  general constraints (instead of cardinality, as in \rsm), such as  observability/controllability requirements, including matroid constraints, towards multi-robot planning. 

\appendices

\section*{Appendix}\label{app:notation}

In this appendix, we provide all proofs.  We use the notation: 
\begin{equation}\label{notation:marginal}
f(\mathcal{X}\;|\;\mathcal{X}')\triangleq f(\mathcal{X}\cup\mathcal{X}')-f(\mathcal{X}'),
\end{equation}
for any $\mathcal{X},\mathcal{X}'$.
Also, $\calX_{1:t}\triangleq (\calX_1,\ldots,\calX_t)$ for any $\calX_1,\ldots,\calX_t$ {(and $(\calX_1,\ldots,\calX_t) \equiv \calX_1\cup \ldots\cup\calX_t$).}

\section*{Appendix A: Preliminary lemmas}\label{app:prelim}

We list lemmas that support the proofs.

\begin{lemma}\label{lem:curvature2}
Consider a non-decreasing $f:2^\calV\mapsto \mathbb{R}$ such that $f(\emptyset)=0$. Then, for any  $\calA, \calB\subseteq \calV$ such that $\calA \cap \calB=\emptyset$, 
\begin{equation*}
f(\calA\cup \calB)\geq \left(1-c_f\right)\left[f(\calA)+f(\calB)\right].
\end{equation*}
\end{lemma}
\paragraph*{Proof of Lemma~\ref{lem:curvature2}}
Let $\calB=\{b_1, b_2, \ldots, b_{|\calB|}\}$. Then, 
\begin{equation}
f(\calA\cup \calB)=f(\calA)+\sum_{i=1}^{|\calB|}f(b_i\;|\; \calA\cup \{b_1, b_2, \ldots, b_{i-1}\}). \label{eq1:lemma_curvature2}
\end{equation} 
The definition of $c_f$ implies
\begin{align}
&f(b_i\;|\;\calA\cup \{b_1, b_2, \ldots, b_{i-1}\})\geq\nonumber\\
& \qquad (1-c_f)\; f(b_i\;|\;\{b_1, b_2, \ldots, b_{i-1}\}). \label{eq2:lemma_curvature2}
\end{align} 
The proof is completed by substituting ineq.~\eqref{eq2:lemma_curvature2} in eq.~\eqref{eq1:lemma_curvature2}, along with $f(\calA)\geq \left(1-c_f\right)f(\calA)$, since $c_f\leq 1$.
\hfill $\blacksquare$

\begin{lemma}\label{lem:curvature}
Consider a non-decreasing $f:2^\calV\mapsto \mathbb{R}$ such that $f(\emptyset)=0$. Then, for any $\calA, \calB\subseteq \calV$ such that $\calA \cap \calB=\emptyset$,
\begin{equation*}
f(\calA\cup \calB)\geq (1-c_f)\left[f(\calA)+\sum_{b \in \calB}f(b)\right].
\end{equation*}
\end{lemma}
\paragraph*{Proof of Lemma~\ref{lem:curvature}}
Let $\calB=\{b_1, b_2, \ldots, b_{|\calB|}\}$. Then, 
\begin{equation}
f(\calA\cup \calB)=f(\calA)+\sum_{i=1}^{|\calB|}f(b_i\;|\;\calA\cup \{b_1, b_2, \ldots, b_{i-1}\}). \label{eq1:lemma_curvature}
\end{equation} 
Now, $c_f$'s Definition~\ref{def:total_curvature} implies
\begin{align}
f(b_i\;|\;\calA\cup \{b_1, b_2, \ldots, b_{i-1}\})&\geq (1-c_f)f(b_i\;|\;\emptyset)\nonumber\\
&=(1-c_f)f(b_i), \label{eq2:lemma_curvature}
\end{align} 
where the latter holds since $f(\emptyset)=0$.
The proof is completed by substituting eq.~\eqref{eq2:lemma_curvature} in eq.~\eqref{eq1:lemma_curvature}, along with $f(\calA)\geq (1-c_f)f(\calA)$, since $c_f\leq 1$. \hfill $\blacksquare$

\begin{lemma}\label{lem:subratio}
Consider a non-decreasing $f:2^\calV\mapsto \mathbb{R}$ such that $f(\emptyset)=0$. Then, for any $\calA, \calB\subseteq \calV$ such that $\calA\setminus\calB\neq \emptyset$,
\begin{equation*}
f(\calA)+(1-c_f) f(\calB)\geq (1-c_f) f(\calA\cup \calB)+f(\calA\cap \calB).
\end{equation*}
\end{lemma}
\paragraph*{Proof of Lemma~\ref{lem:subratio}}
Let $\calA\setminus\calB=\{i_1,i_2,\ldots, i_r\}$, where $r=|\calA-\calB|$.  $c_f$'s Definition~\ref{def:total_curvature} implies $f(i_j\;|\;{(\calA \cap \calB)} \cup \{i_1, i_2, \ldots, i_{j-1}\})\geq (1-c_f) f(i_j\;|\;\calB \cup \{i_1, i_2, \ldots, i_{j-1}\})$,  for any $i=1,\ldots, r$. Summing the $r$ inequalities,
$$f(\calA)-f(\calA\cap \calB)\geq (1-c_f) \left[f(\calA\cup \calB)-f(\calB)\right],$$
which implies the lemma. \hfill $\blacksquare$

\begin{corollary}\label{cor:ineq_from_lemmata}
Consider a non-decreasing $f:2^\calV\mapsto \mathbb{R}$ such that $f(\emptyset)=0$. Then, for any  $\calA, \calB\subseteq \calV$ such that $\calA\cap\calB=\emptyset$, 
\begin{equation*}
f(\calA)+\sum_{b \in \calB}f(b) \geq (1-c_f)  f(\calA\cup \calB).
\end{equation*}
\end{corollary}
\paragraph*{Proof of Corollary~\ref{cor:ineq_from_lemmata}}
 Let $\calB=\{b_1,b_2,\ldots,b_{|\calB|}\}$.   
 
 If $\calA\neq \emptyset$, then
\begin{align}
f(\calA)+\sum_{i=1}^{|\calB|}f(b_i) &\geq (1-c_f) f(\calA)+\sum_{i=1}^{|\calB|}f(b_i)\label{ineq:cor_aux1} \\
& \geq (1-c_f) f(\calA\cup \{b_1\})+\sum_{i=2}^{|\calB|}f(b_i)\nonumber\\
& \geq (1-c_f) f(\calA\cup \{b_1,b_2\})+\sum_{i=3}^{|\calB|}f(b_i)\nonumber\\
& \;\;\vdots \nonumber\\
& \geq (1-c_f) f(\calA\cup \calB),\nonumber
\end{align}
where ineq.~\eqref{ineq:cor_aux1} holds since $0\leq c_f\leq 1$, {and the remaining inequalities are implied by applying Lemma~\ref{lem:subratio} multiple times ($\calA\cap\calB=\emptyset$  implies $\calA\setminus \{b_1\}\neq \emptyset$, $\calA\cup \{b_1\}\setminus \{b_2\}\neq \emptyset$, $\ldots$, $\calA\cup \{b_1,b_2,\ldots, b_{|\calB|-1}\}\setminus \{b_{|\calB|}\}\neq \emptyset$).}

{If $\calA= \emptyset$,  then the proof follows the same reasoning as above but now we need to start from the following inequality, instead of  ineq.~\eqref{ineq:cor_aux1}:
\belowdisplayskip=0pt\begin{equation*}\sum_{i=1}^{|\calB|}f(b_i) \geq (1-c_f) f(\{b_1\})+\sum_{i=2}^{|\calB|}f(b_i).
\end{equation*}
}
\hfill $\blacksquare$

\begin{lemma}\label{lem:adapt_res_greedy_vs_any}
Consider the sets $\calS_{1,1},\ldots,\calS_{T,1}$ selected by \alg's lines~\ref{line:sort}-\ref{line:bait}.  Also, for all $t=1,\ldots,T$, let $\calO_t$ be any subset of $\calV_t\setminus \calS_{t,1}$ such that $|\calO_t|\leq \alpha_t-\beta_t$. Then, 
\begin{equation}\label{eq:adapt_res_greedy_vs_any}
f(\calS_{1,2},\ldots,\calS_{T,2})\geq (1-c_f)^2f(\calO_{1:T}).
\end{equation}
\end{lemma}
\paragraph*{Proof of Lemma~\ref{lem:adapt_res_greedy_vs_any}}
For all $t=1,\ldots,T$\!, let $\calR_t\triangleq \calA_t\setminus \calB_t$; namely, $\calR_t$ is the set that remains after the optimal (worst-case) removal $\calB_t$ from $\calA_t$.  Furthermore, let $s_{t,2}^i\in \calS_{t,2}$ denote the $i$-th element added to $\calS_{t,2}$ per \alg's lines~\ref{line:begin_while}-\ref{line:end_while}; i.e.,  $\calS_{t,2}=\{s_{t,2}^1,\ldots,s_{t,2}^{\alpha_t-\beta_t}\}$.  Additionally,  for all $i=1,\ldots, \alpha_t-\beta_t$, denote $\calS_{t,2}^i\triangleq\{s_{t,2}^1,\ldots,s_{t,2}^{i}\}$, and set $\calS_{t,2}^{0}\triangleq \emptyset$.  Next, order the elements in each $\calO_t$ so that $\calO_t=\{o_t^1,\ldots,o_t^{\alpha_t-\beta_t}\}$  and if $o_t^i\in\calS_{t,2}$, then $o_t^i=s_{t,2}^i$; i.e., order the elements so that the common elements in $\calO_t$ and $\calS_{t,2}$ appear at the same index.
Moreover, for all $i=1,\ldots, \alpha_t-\beta_t$, denote $\calO_{t}^i\triangleq\{o_{t}^1,\ldots,o_{t}^{i}\}$, and also set $\calO_{t}^{0}\triangleq \emptyset$.  Finally, let: {$\calO_{1:t}\triangleq (\calO_1,\ldots, \calO_{t})$}; $\calO_{1:0}\triangleq\emptyset$; {$\calS_{1:t,2}\triangleq (\calS_{1,2},\ldots, \calS_{t,2})$}; and $\calS_{1:0,2}\triangleq\emptyset$. Then,
\begin{align}
f(\calO_{1:T})
&=\sum_{t=1}^{T}\sum_{i=1}^{\alpha_t-\beta_t}f(o_t^i\;|\;\calO_{1:t-1}\cup \calO_t^{i-1})\label{aux10:1}\\
&\leq \frac{1}{1-c_f}\sum_{t=1}^{T}\sum_{i=1}^{\alpha_t-\beta_t}f(o_t^i\;|\;\calR_{1:t-1}\cup \calS_{t,2}^{i-1})\label{aux10:2}\\
&\leq \frac{1}{1-c_f}\sum_{t=1}^{T}\sum_{i=1}^{\alpha_t-\beta_t}f(s_{t,2}^i\;|\;\calR_{1:t-1}\cup \calS_{t,2}^{i-1})\label{aux10:3}\\
&\leq \frac{1}{(1-c_f)^2}\sum_{t=1}^{T}\sum_{i=1}^{\alpha_t-\beta_t}f(s_{t,2}^i\;|\;\calS_{1:t-1,2}\cup \calS_{t,2}^{i-1})\label{aux10:4}\\
%
&=\frac{1}{(1-c_f)^2}f(\calS_{1,2},\ldots,\calS_{T,2})\label{aux10:5}.
\end{align}
where eq.~\eqref{aux10:1} holds due to eq.~\eqref{notation:marginal}; ineq.~\eqref{aux10:2} due to ineq.~\eqref{eq:total_curvature};
ineq.~\eqref{aux10:3} holds since $s_{t,2}^i$ is chosen greedily by the algorithm, given $\calR_{1:t-1}\cup \calS_{t,2}^{i-1}$; ineq.~\eqref{aux10:4} holds for the same reasons as ineq.~\eqref{aux10:2};  eq.~\eqref{aux10:5} holds for the same reasons as eq.~\eqref{aux10:1}.
\hfill $\blacksquare$

\begin{lemma}\label{lem:local_greedy_vs_opt}
Consider the sets $\calS_{1,1},\ldots,\calS_{T,1}$ selected by \alg's lines~\ref{line:sort}-\ref{line:bait}.  Also, for all $t=1,\ldots,T$, let in Algorithm~\ref{alg:local} be $\calK_t=\calV_t\setminus \calS_{t,1}$ and $\delta_t=\alpha_t-\beta_t$.
Finally, consider $\calP_t$ such that $\calP_t\subseteq \calK_t$, $|\calP_t|\leq\delta_t$, and $f({\calP}_{1:T})$ is maximal, that is,
\begin{equation}\label{eq:local_greedy}
\calP_{1:T}\in \arg \max_{\bar{\calP}_1\subseteq \calK_1, |\bar{\calP}_1|\leq \delta_1}\cdots \max_{\bar{\calP}_T\subseteq \calK_T, |\bar{\calP}_T|\leq \delta_T} f(\bar{\calP}_{1:T}).
\end{equation} Then, $f(\calM_{1:T})\geq (1-c_f)f(\calP_{1:T})$.
\end{lemma}
\paragraph*{Proof of Lemma~\ref{lem:local_greedy_vs_opt}}  The proof is the same as that of~\cite[Theorem~6]{sviridenko2017optimal}.
\hfill $\blacksquare$

\begin{corollary}\label{lem:adapt_res_greedy_vs_opt_local_greedy}
Consider the sets $\calS_{1,1},\ldots,\calS_{T,1}$ selected by \alg's lines~\ref{line:sort}-\ref{line:bait}, as well as, the sets $\calS_{1,2},\ldots,\calS_{T,2}$ selected by \alg's lines~\ref{line:begin_while}-\ref{line:end_while}.  Finally, consider $\calK_t=\calV_t\setminus\calS_{t,1}$ and $\delta_t=\alpha_t-\beta_t$, and   $\calP_t$ per eq.~\eqref{eq:local_greedy}. Then, 
\begin{equation*}\label{ineq:main_adapt_res_greedy_vs_opt_local_greedy}
f(\calS_{1,2},\ldots,\calS_{T,2})\geq (1-c_f)^3f(\calP_{1:T}).
\end{equation*}
\end{corollary}
\paragraph*{Proof of Corollary~\ref{lem:adapt_res_greedy_vs_opt_local_greedy}}
Let $\calO_t=\calM_t$ in ineq.~\eqref{eq:adapt_res_greedy_vs_any}.  Then,
\begin{equation}\label{eq:adapt_res_greedy_vs_opt_local_greedy}
f(\calS_{1,2},\ldots,\calS_{T,2})\geq (1-c_f)^2f(\calM_{1:T}).
\end{equation}
Using in ineq.~\eqref{eq:adapt_res_greedy_vs_opt_local_greedy} Lemma~\ref{lem:local_greedy_vs_opt}, the proof is complete.
\hfill $\blacksquare$

\begin{lemma}\label{lem:from_max_to_minmax}
Per the notation in Corollary~\ref{lem:adapt_res_greedy_vs_opt_local_greedy}, 
for all $t=1,\ldots,T$\!, consider $\calK_t=\calV_t\setminus\calS_{t,1}$, $\delta_t=\alpha_t-\beta_t$, and $\calP_t$ per eq.~\eqref{eq:local_greedy}. Then, it holds true that
\begin{equation}\label{eq:toprovefrom_max_to_minmax}
f(\calP_{1:T})\geq f^\star.
\end{equation}
\end{lemma}
\paragraph*{Proof of Lemma~\ref{lem:from_max_to_minmax}}
We use the notation
\begin{align}
\begin{split}\label{def:h_eq}
&\!\!\!h(\calS_{1,1},\ldots,\calS_{T,1})\triangleq\\
&\max_{\bar{\calP}_1\subseteq \calV_1\setminus\calS_{1,1}, |\bar{\calP}_1|\leq \delta_1}\cdots \max_{\bar{\calP}_T\subseteq \calV_T\setminus\calS_{T,1}, |\bar{\calP}_T|\leq \delta_T} f(\bar{\calP}_{1:T}).
\end{split}
\end{align}

{For any  $\hat{\calP}_1,\ldots,\hat{\calP}_T$ such that  $\hat{\calP}_t\subseteq \calV_t\setminus\calS_{t,1}$ and $|\hat{\calP}_t|\;\leq \delta_t$ (for all $t=1,\ldots,T$), $h$'s definition in eq.~\eqref{def:h_eq} implies
\begin{equation}
h(\calS_{1,1},\ldots,\calS_{T,1})\geq f(\hat{\calP}_1,\ldots,\hat{\calP}_T).\label{aux40:0}
\end{equation}
Since ineq.~\eqref{aux40:0} holds for any $\hat{\calP}_T$ such that  $\hat{\calP}_T\subseteq \calV_T\setminus\calS_{T,1}$ and  $|\hat{\calP}_T|\;\leq \delta_T$, then it also holds for the $\hat{\calP}_T$ that maximizes the right-hand-side of ineq.~\eqref{aux40:0}, i.e.,
\begin{equation}
h(\calS_{1,1},\ldots,\calS_{T,1})\geq\max_{\bar{\calP}_T\subseteq \calV_T\setminus\calS_{T,1}, |\bar{\calP}_T|\leq \delta_T}f(\hat{\calP}_{1:T-1},\bar{\calP}_T).\label{aux400:00}
\end{equation}
Treating in ineq.~\eqref{aux400:00} the set $\calS_{T,1}$ as a free variable, since \eqref{aux400:00} holds for any $\calS_{T,1}\subseteq \calV_T$ such that $|\calS_{T,1}|\;\leq \beta_T$, then 
\begin{align}
&\min_{\bar{\calB}_{T}\subseteq \calV_T, |\bar{\calB}_{T}|\leq \beta_T}h(\calS_{1,1},\ldots,\calS_{T-1,1},\bar{\calB}_{T})\geq\nonumber\\
&\hspace{0.5cm}\min_{\bar{\calB}_{T}\subseteq \calV_T, |\bar{\calB}_{T}|\leq \beta_T}\;\;\max_{\bar{\calP}_T\subseteq \calV_T\setminus\bar{\calB}_{T}, |\bar{\calP}_T|\leq \delta_T}f(\hat{\calP}_{1:T-1},\bar{\calP}_T).\label{aux40:1}
\end{align}
Specifically, ineq.~\eqref{aux40:1} holds true for the same reason the following holds true: given a set function $f_1:\calI\mapsto \mathbb{R}$ (representing the left-hand-side of ineq.~\eqref{aux400:00}) and a set function $f_2:\calI\mapsto \mathbb{R}$ (representing the right-hand-side of ineq.~\eqref{aux400:00}), where $\calI=\{\calS: \calS\subseteq \calV_T, |\calS|\;\leq \beta_T\}$, if $f_1(\calS)\geq f_2(\calS)$ for every $\calS \in \calI$ (as ineq.~\eqref{aux400:00} defines), then also the minimum of $f_1$ must be greater or equal to the minimum of $f_2$, i.e., $\min_{\calS \in \calI} f_1(\calS)\geq \min_{\calS \in \calI} f_2(\calS)$.  The reason: if $\calS_1^\star\in \arg\min_{\calS \in \calI} f_1(\calS)$, then $f_1(\calS_1^\star)\geq f_2(\calS_1^\star)$, since $f_1(\calS)\geq f_2(\calS)$ for every $\calS \in \calI$; but also $f_2(\calS_1^\star)\geq \min_{\calS \in \calI} f_2(\calS)$, and, as a result, indeed, $\min_{\calS \in \calI} f_1(\calS)\geq \min_{\calS \in \calI} f_2(\calS)$.  Changing the symbol of the dummy variable $\calS$ to $\bar{\calB}_T$, we get ineq.~\eqref{aux40:1}.}

Denote now the right-hand-side of ineq.~\eqref{aux40:1} by $z(\hat{\calP}_{1:T-1})$. Since  $\delta_T=\alpha_T-\beta_T$, and for $\bar{\calP}_T$ in  ineq.~\eqref{aux40:1} it is $\bar{\calP}_T\subseteq \calV_T\setminus\bar{\calB}_{T}$ and $|\bar{\calP}_T|\leq \delta_T$, then it equivalently holds:
\begin{align}
&z(\hat{\calP}_{1:T-1})=\nonumber\\
&\min_{\bar{\calB}_{T}\subseteq \calV_T, |\bar{\calB}_{T}|\leq \beta_T}\;\;\max_{\bar{\calA}_T\subseteq \calV_T, |\bar{\calA}_T|\leq \alpha_T}f(\hat{\calP}_{1:T-1},\bar{\calA}_T\setminus \bar{\calB}_{T}).\label{aux40:2}
\end{align}
Let in ineq.~\eqref{aux40:2} $w(\bar{\calA}_T\setminus \bar{\calB}_{T})\triangleq f(\hat{\calP}_{1:T-1},\bar{\calA}_T\setminus \bar{\calB}_{T})$; we prove that the following holds true:
\begin{equation}
z(\hat{\calP}_{1:T-1})\geq \max_{\bar{\calA}_T\subseteq \calV_T, |\bar{\calA}_T|\leq \alpha_T}\;\;\min_{\bar{\calB}_{T}\subseteq \calV_T, |\bar{\calB}_{T}|\leq \beta_T}w(\bar{\calA}_T\setminus \bar{\calB}_{T}).\label{aux40:3}
\end{equation}
{Particularly,  for any $\hat{\calA}_T\subseteq \calV_T$ such that $|\hat{\calA}_T|\;\leq \alpha_T$, and any $\hat{\calS}_{T,1}\subseteq \calV_T$ such that $|\hat{\calS}_{T,1}|\;\leq \beta_T$, it is
\begin{equation}\label{aux_ineq_000}
\max_{\bar{\calA}_T\subseteq \calV_T, |\bar{\calA}_T|\leq \alpha_T} w(\bar{\calA}_T\setminus \hat{\calS}_{T,1})\geq w(\hat{\calA}_T\setminus \hat{\calS}_{T,1}).
\end{equation}
From ineq.~\eqref{aux_ineq_000}, {following the same reasoning as for the derivation of ineq.~\eqref{aux40:1} from ineq.~\eqref{aux400:00}, considering in ineq.~\eqref{aux_ineq_000}  the $\hat{\calS}_{T,1}$ to be the free variable}, we get
\begin{align}
\min_{\bar{\calB}_{T}\subseteq \calV_T, |\bar{\calB}_{T}|\leq \beta_T}\;\;\max_{\bar{\calA}_T\subseteq \calV_T, |\bar{\calA}_T|\leq \alpha_T} w(\bar{\calA}_T\setminus \bar{\calB}_{T})&\geq \nonumber\\
&\hspace{-3cm}\min_{\bar{\calB}_{T}\subseteq \calV_T, |\bar{\calB}_{T}|\leq \beta_T}w(\hat{\calA}_T\setminus \bar{\calB}_{T}).\label{aux40:4}
\end{align}
Now, ineq.~\eqref{aux40:4} implies ineq.~\eqref{aux40:3}.  The reason: \eqref{aux40:4} holds for any $\bar{\calA}_T\subseteq \calV_T$ such that $|\bar{\calA}_T|\;\leq \alpha_T$, while the left-hand-side of \eqref{aux40:4} is equal to $z(\hat{\calP}_{1:T-1})$, which is independent of  $\hat{\calA}_T$; therefore, if we maximize the right-hand-side of \eqref{aux40:4} with respect to $\bar{\calA}_T$, then indeed we get \eqref{aux40:3}.}
  
All in all, due to ineq.~\eqref{aux40:3},  ineq.~\eqref{aux40:1} becomes:
\begin{align}
&\min_{\bar{\calB}_{T}\subseteq \calV_T, |\bar{\calB}_{T}|\leq \beta_T}h(\calS_{1,1},\ldots,\calS_{T-1,1},\bar{\calB}_{T})\geq\nonumber\\
&\max_{\bar{\calA}_T\subseteq \calV_T, |\bar{\calA}_T|\leq \alpha_T}\;\;\min_{\bar{\calB}_{T}\subseteq \calV_T, |\bar{\calB}_{T}|\leq \beta_T} f(\hat{\calP}_{1:T-1},\bar{\calA}_T\setminus \bar{\calB}_{T}).\label{aux40:5}
\end{align}
The left-hand-side of ineq.~\eqref{aux40:5} is a function of $\calS_{1,1},\ldots,\calS_{T-1,1}$; denote it as $h'(\calS_{1,1},\ldots,\calS_{T-1,1})$.  Similarly, the right-hand-side of ineq.~\eqref{aux40:5} is a function of $\hat{\calP}_{1:T-1}$; denote it as $f'(\hat{\calP}_{1:T-1})$. Given these notations, ineq.~\eqref{aux40:5} is equivalently written as
\begin{align}
&h'(\calS_{1,1},\ldots,\calS_{T-1,1})\geq f'(\hat{\calP}_{1:T-1}),\label{aux40:6}
\end{align}
which has the same form as ineq.~\eqref{aux40:0}.  Therefore, by following the same steps as those we used from ineq.~\eqref{aux40:0} and onward, we similarly get
\begin{align}
&\min_{\bar{\calB}_{T-1}\subseteq \calV_{T-1}, |\bar{\calB}_{T-1}|\leq \beta_{T-1}}h'(\calS_{1,1},\ldots,\calS_{T-2,1},\bar{\calB}_{T-1})\geq\nonumber\\
&\max_{\bar{\calA}_{T-1}\subseteq \calV_{T-1}, |\bar{\calA}_{T-1}|\leq \alpha_{T-1}}\;\;\min_{\bar{\calB}_{T-1}\subseteq \calV_{T-1}, |\bar{\calB}_{T-1}|\leq \beta_{T-1}}\nonumber\\ &\hspace{3cm}f'(\hat{\calP}_{1:T-2},
\bar{\calA}_{T-1}\setminus \bar{\calB}_{T-1}),\label{aux40:7}
\end{align}
{which, given the definitions of $h'(\cdot)$ and $f'(\cdot)$,  is equivalent to
\begin{align*}
&\min_{\bar{\calB}_{T-1}\subseteq \calV_{T-1}, |\bar{\calB}_{T-1}|\leq \beta_{T-1}}\;\; \min_{\bar{\calB}_{T}\subseteq \calV_T, |\bar{\calB}_{T}|\leq \beta_T}\\
&\qquad\qquad\qquad\qquad h(\calS_{1,1},\ldots,\calS_{T-2,1},\bar{\calB}_{T-1},\bar{\calB}_{T})\geq\nonumber\\
&\max_{\bar{\calA}_{T-1}\subseteq \calV_{T-1}, |\bar{\calA}_{T-1}|\leq \alpha_{T-1}}\;\;\min_{\bar{\calB}_{T-1}\subseteq \calV_{T-1}, |\bar{\calB}_{T-1}|\leq \beta_{T-1}}\nonumber\\
&\hspace{1.4cm}\max_{\bar{\calA}_T\subseteq \calV_T, |\bar{\calA}_T|\leq \alpha_T}\;\;\min_{\bar{\calB}_{T}\subseteq \calV_T, |\bar{\calB}_{T}|\leq \beta_T}\\ &\hspace{2.5cm} f(\hat{\calP}_{1:T-2},\bar{\calA}_{T-1}\setminus \bar{\calB}_{T-1}, \bar{\calA}_T\setminus \bar{\calB}_{T}).
\end{align*}
Eq.~\eqref{aux40:7}  has the same form as ineq.~\eqref{aux40:5}.}  Therefore, repeating the same steps as above for another $T-2$ times  {(starting now from \eqref{aux40:7} instead of \eqref{aux40:5})}, we get 
\begin{align}
&\min_{\bar{\calB}_{1}\subseteq \calV_1, |\bar{\calB}_{1}|\leq \beta_1}\cdots\min_{\bar{\calB}_{T}\subseteq \calV_T, |\bar{\calB}_{T}|\leq \beta_T}h(\bar{\calB}_{1:T})\geq\nonumber\\
&\max_{\bar{\calA}_{1}\subseteq \calV_{1}, |\bar{\calA}_{1}|\leq \alpha_{1}}\;\min_{\bar{\calB}_{1}\subseteq \calV_{1}, |\bar{\calB}_{1}|\leq \beta_{1}}\!\!\cdots\!\!\max_{\bar{\calA}_{T}\subseteq \calV_{T}, |\bar{\calA}_{T}|\leq \alpha_{T}}\;\min_{\bar{\calB}_{1}\subseteq \calV_{T}, |\bar{\calB}_{T}|\leq \beta_{T}}\nonumber\\ &\hspace{5cm}f(\bar{\calA}_{1}\setminus \bar{\calB}_{1},\ldots,
\bar{\calA}_{T}\setminus \bar{\calB}_{T}),\label{aux40:40}
\end{align}
which implies ineq.~\eqref{eq:toprovefrom_max_to_minmax} since the right-hand-side of ineq.~\eqref{aux40:40} is equal to the right-hand-side of ineq.~\eqref{eq:toprovefrom_max_to_minmax}, while for the left-hand-side of ineq.~\eqref{aux40:40} the following holds:
\belowdisplayskip=-9pt\begin{equation*}
\min_{\bar{\calB}_{1}\subseteq \calV_1, |\bar{\calB}_{1}|\leq \beta_1}\cdots\min_{\bar{\calB}_{T}\subseteq \calV_T, |\bar{\calB}_{T}|\leq \beta_T}h(\bar{\calB}_{1:T})\leq f(\calP_{1:T}).
\end{equation*}
\hfill $\blacksquare$

\section*{Appendix B: Proof of Proposition~\ref{prop:runtime}}

We compute the running time of \alg's line~3  and lines 5-8. \mbox{Line 3} needs $|\calV_t|\tau_f+|\calV_t|\log(|\calV_t|)+|\calV_t|+O(\log(|\calV_t|))$ time: it asks for $|\calV_t|$ evaluations of $f$, and their sorting, which takes $|\calV_t|\log(|\calV_t|)+|\calV_t|+O(\log(|\calV_t|))$ time (using, e.g., the merge sort algorithm).  Lines 5-8 need $(\alpha_t-\beta_t)[|\calV_t|\tau_f+|\calV_t|]$ time: the while loop is repeated $\alpha_t-\beta_t$ times, and during each loop at most $|\calV_t|$ evaluations of $f$ are needed (line~5), plus at most $|\calV_t|$ steps for a maximal element
to be found (line~6).  Hence, \alg runs at each $t$ in $(\alpha_t-\beta_t)[|\calV_t|\tau_f+|\calV_t|]+|\calV_t|\log(|\calV_t|)+|\calV_t|+O(\log(|\calV_t|))=O(|\calV_t|(\alpha_t-\beta_t)\tau_f)$~time. 

\section*{Appendix C: Proof of Theorem~\ref{th:alg_rob_sub_max_performance}}

\begin{figure}[t]
\def \setAone{ (0,0) circle (1cm) }
\def \setBone{ (.5,0) circle (0.4cm)}
\def \setAtwo{ (2.5,0) circle (1cm) }
\def \setBtwo{ (3.0,0) circle (0.4cm)}
\def \myrectangle{ (-1.5, -1.5) rectangle (4, 1.5) }
\begin{center}
\begin{tikzpicture}
\draw \myrectangle node[below left]{$\calV_t$};
\draw \setAone node[left]{$\mathcal{S}_{t,1}$};
\draw \setBone node[]{$\calB_{t,1}^\star$};
\draw \setAtwo node[left]{$\mathcal{S}_{t,2}$};
\draw \setBtwo node[]{$\calB_{t,2}^\star$};
\end{tikzpicture}
\end{center}
\caption{Venn diagram, where the sets $\mathcal{S}_{t,1}, \mathcal{S}_{t,2},\calB_{t,1}^\star, \calB_{t,2}^\star$ are as follows: per \alg, $\mathcal{S}_{t,1}$  and $\mathcal{S}_{t,2}$ are such that $\calA_t=\mathcal{S}_{t,1}\cup \mathcal{S}_{t,2}$.  Additionally, due to their construction,  $\mathcal{S}_{t,1}\cap \mathcal{S}_{t,2}=\emptyset$. Next, 
$\calB_{t,1}^\star$ and $\calB_{t,2}^\star$ are such that  $\calB_{t,1}^\star=\calB^\star_{1:T}\cap\mathcal{S}_{t,1}$, and $\calB_2^\star=\calB^\star_{1:T}\cap\mathcal{S}_{t,2}$; therefore, $\calB_{t,1}^\star\cap \calB_{t,2}^\star=\emptyset$ and $\calB^\star_{1:T}=(\calB_{1,1}^\star\cup \calB_{1,2}^\star)\cup \cdots\cup(\calB_{T,1}^\star\cup \calB_{T,2}^\star)$.  
}\label{fig:venn_diagram_for_proof}
\end{figure}
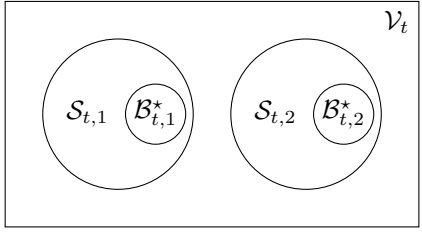

We first prove ineq.~\eqref{ineq:bound_non_sub} and then~\eqref{ineq:bound_sub}.
We use the notation:
\begin{itemize}
\item  $\calS_{t,1}^+\triangleq \calS_{t,1}\setminus \calB^\star_{t}$, i.e., $\calS_{t,1}^+$ is the remaining set after the optimal (worst-case) removal $\calB^\star_{t}$;
\item $\calS_{t,2}^+\triangleq \calS_{t,2}\setminus \calB^\star_{t}$;
\item $\calP_{1:T}$ be a solution to eq.~\eqref{eq:local_greedy}.
\end{itemize}
\paragraph*{Proof of ineq.~\eqref{ineq:bound_non_sub}} For $T>1$, we have:
\belowdisplayskip=8pt\begin{align}
&\!\!\!f(\calA_{1:T}\setminus \calB^\star_{1:T})\nonumber\\
&=f(\calS_{1,1}^+\cup \calS_{1,2}^+,\ldots,\calS_{T,1}^+\cup \calS_{T,2}^+)\label{ineq2:aux_14}\\
&\geq (1-c_f)\sum_{t=1}^T\sum_{v\in\calS_{t,1}^+\cup \calS_{t,2}^+}f(v)\label{ineq2:aux_15}\\
&\geq (1-c_f)\sum_{t=1}^T\sum_{v\in\calS_{t,2}}f(v)\label{ineq2:aux_16}\\
&\geq (1-c_f)^2f(\calS_{1,2},\ldots,\calS_{T,2})\label{ineq2:aux_17}\\
&\geq (1-c_f)^5f(\calP_{1:T})\label{ineq2:aux_18}\\
&\geq (1-c_f)^5f^\star,\label{ineq2:aux_19}
\end{align}
where eq.~\eqref{ineq2:aux_14} follows from the definitions of $\calS_{t,1}^+$ and $\calS_{t,2}^+$; ineq.~\eqref{ineq2:aux_15} follows from ineq.~\eqref{ineq2:aux_14}, due to Lemma~\ref{lem:curvature}; ineq.~\eqref{ineq2:aux_16} follows from ineq.~\eqref{ineq2:aux_15}, because: for all  $v \in \calS_{t,1}^+$ and $v' \in \calS_{t,2}\setminus \calS_{t,2}^+$ it is $f(v)\geq f(v')$,
 and $\calS_{t,2}=(\calS_{t,2}\setminus \calS_{t,2}^+)\cup \calS_{t,2}^+$; ineq.~\eqref{ineq2:aux_17} follows from ineq.~\eqref{ineq2:aux_16} due to Corollary~\ref{cor:ineq_from_lemmata}; ineq.~\eqref{ineq2:aux_18} follows from ineq.~\eqref{ineq2:aux_17} due to Corollary~\ref{lem:adapt_res_greedy_vs_opt_local_greedy}; finally, ineq.~\eqref{ineq2:aux_19} follows from ineq.~\eqref{ineq2:aux_18} due to Lemma~\ref{lem:from_max_to_minmax}. 
 
 For $T=1$, the proof follows the same steps up to ineq.~\eqref{ineq2:aux_17}, at which point $f(\calS_{1,2})\geq (1-c_f)f(\calP_1)$ instead, due to Lemma~\ref{lem:local_greedy_vs_opt} (since $\calS_{1,2}=\calM_1$).
 \hfill $\blacksquare$

\paragraph*{Proof of ineq.~\eqref{ineq:bound_sub}}
For $T>1$ we follow similar steps:
\begin{align}
&\!\!\!f(\calA_{1:T}\setminus \calB^\star_{1:T})\nonumber\\
&=f(\calS_{1,1}^+\cup \calS_{1,2}^+,\ldots,\calS_{T,1}^+\cup \calS_{T,2}^+)\label{ineq5:aux_14}\\
&\geq (1-\kappa_f)\sum_{t=1}^T\sum_{v\in\calS_{t,1}^+\cup \calS_{t,2}^+}f(v)\label{ineq5:aux_15}\\
&\geq (1-\kappa_f)\sum_{t=1}^T\sum_{v\in\calS_{t,2}}f(v)\label{ineq5:aux_16}
\end{align}
\begin{align}
&\geq (1-\kappa_f)f(\calS_{1,2},\ldots,\calS_{T,2})\label{ineq5:aux_17}\\
&\geq (1-\kappa_f)^4f(\calP_{1:T})\label{ineq5:aux_18}\\
&\geq (1-\kappa_f)^4 f^\star,\label{ineq5:aux_19}
\end{align}
where eq.~\eqref{ineq5:aux_14} follows from the definitions of~$\calS_{t,1}^+$ and $\calS_{t,2}^+$; ineq.~\eqref{ineq5:aux_15} follows from ineq.~\eqref{ineq5:aux_14} due to Lemma~\ref{lem:curvature2} and the fact that $c_f=\kappa_f$ for $f$ being submodular; ineq.~\eqref{ineq5:aux_16} follows from ineq.~\eqref{ineq5:aux_15} because for all $v \in \calS_{t,1}^+$ and $v' \in \calS_{t,2}\setminus \calS_{t,2}^+$ it is $f(v)\geq f(v')$, while $\calS_{t,2}=(\calS_{t,2}\setminus \calS_{t,2}^+)\cup \calS_{t,2}^+$; ineq.~\eqref{ineq5:aux_17} follows from ineq.~\eqref{ineq5:aux_16} because  $f$ is submodular and, as~a result, $f(\mathcal{S})+f(\mathcal{S}')\geq f(\mathcal{S}\cup \mathcal{S}')$, for any $\mathcal{S}, \mathcal{S}'\subseteq \calV$~\cite[Proposition 2.1]{nemhauser78analysis};  ineq.~\eqref{ineq5:aux_18} follows from ineq.~\eqref{ineq5:aux_17} due to Corollary~\ref{lem:adapt_res_greedy_vs_opt_local_greedy}, along with the fact that since $f$ is monotone submodular it is $c_f=\kappa_f$; finally, ineq.~\eqref{ineq5:aux_19} follows from ineq.~\eqref{ineq5:aux_18} due to Lemma~\ref{lem:from_max_to_minmax}. 

For $T=1$, the proof follows the same steps up to ineq.~\eqref{ineq5:aux_17}, at which point $f(\calS_{1,2})\geq 1/\kappa_f(1-e^{-\kappa_f})f(\calP_1)$, due to \cite[Theorem~5.4]{conforti1984curvature}.
\hfill $\blacksquare$

\section*{Appendix D: Proof of Theorem~\ref{th:aposteriori}}

To prove ineq.~\eqref{ineq:a_posteriori_bound_sub}, we have
\begin{align}
&\!\!\!f(\calA_{1:t}\setminus \calB^\star_{1:t})\nonumber\\
&=f(\calM_{1:t})\frac{f(\calA_{1:t}\setminus \calB^\star_{1:t})}{f(\calM_{1:t})}\nonumber\\
&\geq \left\{\begin{array}{lr}
\frac{1-{e^{-\kappa_f}}}{\kappa_f}
\frac{f(\calA_{1}\setminus \mathcal{B}^\star_{1})}{f(\calM_{1})}f(\calP_{1}), & t=1;\\
\frac{1}{1+\kappa_f}
\frac{f(\calA_{1:t}\setminus \mathcal{B}^\star_{1:t})}{f(\calM_{1:t})}f(\calP_{1:t}), & t>1,
\end{array} \right. \label{eq:aaa}
\end{align}

\begin{align}
&\geq \left\{\begin{array}{lr}
\frac{1-{e^{-\kappa_f}}}{\kappa_f}
\frac{f(\calA_{1}\setminus \mathcal{B}^\star_{1})}{f(\calM_{1})}f^\star_1, & t=1;\\
\frac{1}{1+\kappa_f}
\frac{f(\calA_{1:t}\setminus \mathcal{B}^\star_{1:t})}{f(\calM_{1:t})}f^\star_t, & t>1,
\end{array} \right. \label{eq:aaa2}
\end{align}
where ineq.~\eqref{eq:aaa} holds since \cite[Theorem~5.4]{conforti1984curvature} implies $f(\calM_1)\geq 1/\kappa_f(1-e^{-\kappa_f})f(\calP_1)$, while
\cite[Theorem~2.3]{conforti1984curvature} implies $f(\calM_{1:t})\geq 1/(1+\kappa_f)f(\calP_{1:t})$.  Finally, ineq.~\eqref{eq:aaa2} is proved following the same steps as in Lemma~\ref{lem:from_max_to_minmax}'s proof.

The proof of ineq.~\eqref{ineq:a_posteriori_bound_non_sub} follows similar steps as above but it is based instead on Lemma~\ref{lem:local_greedy_vs_opt}.

\section*{Appendix E: Proof of Theorem~\ref{th:tightness_optimality}}

It can be verified that for $t=1$ eq.~\eqref{ineq:a_posteriori_bound_sub} is tight for any $\beta_t\leq \alpha_t$ for the families of functions in \cite[Theorem~5.4]{conforti1984curvature}, and for $t>1$ it is tight for the families of functions in \cite[Theorem~2.12]{conforti1984curvature}.  Similarly, it can be verified eq.~\eqref{ineq:a_posteriori_bound_non_sub} is optimal for the families of functions in \cite[Theorem~8]{sviridenko2017optimal} for $\alpha_t=|\calV_t|^{1/2}$ and any $\beta_t\leq \alpha_t-|\calV_t|^{1/3}$. 

\section*{Acknowledgements}
We thank Rakesh Vohra of the University of Pennsylvania, and Luca Carlone of the Massachusetts Institute of Technology for inspiring comments and discussions. 

\bibliographystyle{IEEEtran}
\bibliography{references,security_references}

\begin{biography}[{\includegraphics[width=1in,height=1.25in,keepaspectratio]{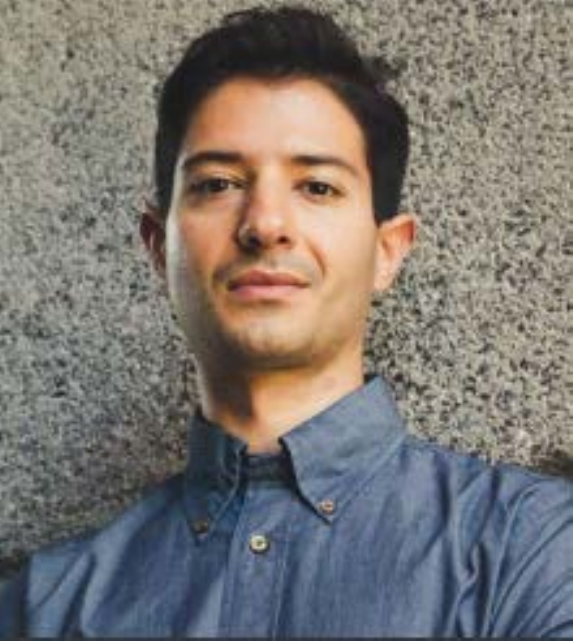}}]{Vasileios Tzoumas} received his Ph.D. in Electrical and Systems Engineering at the University of Pennsylvania (2018). He holds a Master of Arts in Statistics from the Wharton School of Business at the University of Pennsylvania (2016); a Master of Science in Electrical Engineering from the University of Pennsylvania (2016); and a diploma in Electrical and Computer Engineering from the National Technical University of Athens (2012).
		Vasileios is as an Assistant Professor in the Department of Aerospace Engineering, University of Michigan, Ann Arbor. 
		Previously, he was at the Massachusetts Institute of Technology (MIT), in the Department of Aeronautics and Astronautics, and in the Laboratory for Information and Decision Systems (LIDS), were he was a research scientist (2019-2020), and a post-doctoral associate (2018-2019). Vasileios was a visiting Ph.D. student at the Institute for Data, Systems, and Society (IDSS) at MIT during 2017.
		Vasileios works on control, learning, and perception, as well as combinatorial and distributed optimization, with applications to robotics, cyber-physical systems, and self-reconfigurable aerospace systems.  He aims for trustworthy collaborative autonomy.  His work includes foundational results on robust and adaptive combinatorial optimization, with applications to multi-robot information gathering for resiliency against robot failures and adversarial removals.
		Vasileios is a recipient of the Best Paper Award in Robot Vision at the 2020 IEEE International Conference on Robotics and Automation (ICRA), and was a Best Student Paper Award Finalist at the 2017 IEEE Conference in Decision and Control (CDC).
\end{biography}

\begin{biography}[{\includegraphics[width=1in,height=1.25in,clip,keepaspectratio]{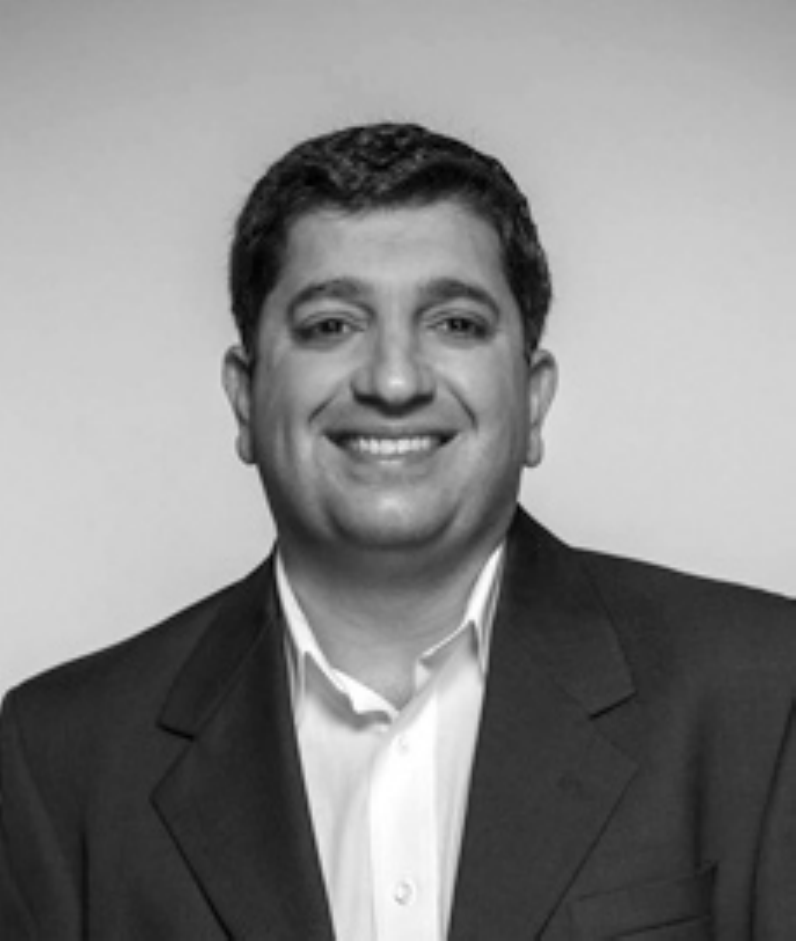}}] {Ali Jadbabaie} (S'99-M'08-SM'13-F'15) is the JR East Professor of Engineering and Associate Director of the Institute for Data, Systems and Society at MIT, where he is also on the faculty of the department of civil and environmental engineering and a principal investigator in the Laboratory for Information and Decision Systems (LIDS). He is the director of the Sociotechnical Systems Research Center, one of MIT's 13 laboratories. He received his Bachelors (with high honors) from Sharif University of Technology in Tehran, Iran, a Masters degree in electrical and computer engineering from the University of New Mexico, and his Ph.D.~in control and dynamical systems from the California Institute of Technology. He was a postdoctoral scholar at Yale University before joining the faculty at Penn in July 2002. Prior to joining MIT faculty, he was the Alfred Fitler Moore a Professor of Network Science and held secondary appointments in computer and information science and operations, information and decisions in the Wharton School. He was the inaugural editor-in-chief of IEEE Transactions on Network Science and Engineering, a new interdisciplinary journal sponsored by several IEEE societies. He is a recipient of a National Science Foundation Career Award, an Office of Naval Research Young Investigator Award, the O. Hugo Schuck Best Paper Award from the American Automatic Control Council, and the George S. Axelby Best Paper Award from the IEEE Control Systems Society. His students have been winners and finalists of student best paper awards at various ACC and CDC conferences. He is an IEEE fellow and a recipient of the Vannevar Bush Fellowship from the office of Secretary of Defense. His current research interests include the interplay of dynamic systems and networks with specific emphasis on multi-agent coordination and control, distributed optimization, network science, and network economics.
\end{biography}

\begin{biography}[{\includegraphics[width=1in,height=1.25in,clip,keepaspectratio]{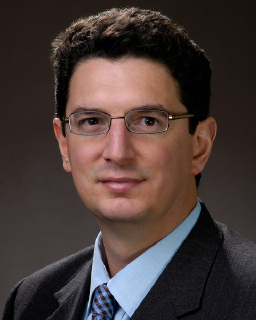}}]{George J.~Pappas} (S'90-M'91-SM'04-F'09) received the Ph.D.~degree in electrical engineering and computer sciences from the University of California, Berkeley, CA, USA, in  1998. He is currently the Joseph Moore Professor and Chair of the Department of Electrical and Systems Engineering, University of Pennsylvania, Philadelphia, PA, USA. He  also holds a secondary appointment with the Department of Computer and Information Sciences and the Department of Mechanical Engineering and Applied Mechanics. He is a Member of the GRASP Lab and the PRECISE Center. He had previously served as the Deputy Dean for Research with the School of Engineering and Applied Science. His research interests include control theory and, in particular, hybrid systems, embedded systems, cyberphysical systems, and hierarchical and distributed control systems, with applications to unmanned aerial vehicles, distributed robotics, green buildings, and biomolecular networks. Dr. Pappas has received various awards, such as the Antonio Ruberti Young Researcher Prize, the George S. Axelby Award, the Hugo Schuck Best Paper Award, the George H. Heilmeier Award, the National Science Foundation PECASE award and numerous best student papers awards.
\end{biography}

\end{document}